\documentclass{article}
\pdfoutput=1

\usepackage{arxiv}

\usepackage[utf8]{inputenc} 
\usepackage[T1]{fontenc}    
\usepackage[hypertexnames=false]{hyperref}       
\usepackage{url}            
\usepackage{amsfonts}       
\usepackage{nicefrac}       
\usepackage{microtype}      

\usepackage{lipsum}         
\usepackage{graphicx}
\usepackage{natbib}
\usepackage{doi}

\usepackage{mathtools,amsthm,amssymb}
\usepackage{algorithm}
\usepackage{algpseudocodex}
\usepackage{xcolor}
\allowdisplaybreaks[4]
\usepackage{cleveref}       
\usepackage{etoolbox}
\usepackage{nicematrix}
\usepackage{enumitem}
\usepackage{threeparttablex}
\usepackage{varwidth}
\usepackage{colortbl}
\usepackage{booktabs}

\hypersetup{
  colorlinks=true,
	linkcolor=red!75!black,
	citecolor=blue!50!black
}

\title{On Linear Convergence in Smooth Convex-Concave Bilinearly-Coupled Saddle-Point Optimization: Lower Bounds and Optimal Algorithms}

\date{November 21, 2024}

\newif\ifuniqueAffiliation
\uniqueAffiliationtrue

\ifuniqueAffiliation 
\author{
	Dmitry Kovalev \\
	Yandex Research \\
	\texttt{dakovalev1@gmail.com}
	\And
	Ekaterina Borodich\\
	Moscow Institute of Physics and Technology\\
	\texttt{borodich.ed@phystech.edu}
}
\else
\usepackage{authblk}

\author[1]{%
	Dmitry Kovalev\thanks{\texttt{dakovalev1@gmail.com}}%
}
\affil[1]{Yandex Research}
\fi
\renewcommand{\headeright}{}
\renewcommand{\shorttitle}{}

\newtheorem{theorem}{Theorem}
\newtheorem{lemma}{Lemma}
\newtheorem{corollary}{Corollary}

\newtheorem{assumption}{Assumption}

\newtheorem{definition}{Definition}
\crefname{assumption}{assumption}{assumptions}

\DeclarePairedDelimiter{\abs}{\lvert}{\rvert}
\DeclarePairedDelimiter{\norm}{\|}{\|}
\DeclarePairedDelimiter{\sqn}{\|}{\|^2}
\DeclarePairedDelimiter{\ceil}{\lceil}{\rceil}
\DeclarePairedDelimiter{\floor}{\lfloor}{\rfloor}

\def\<#1,#2>{\langle #1,#2\rangle}

\DeclareMathOperator{\dom}{dom}
\DeclareMathOperator{\diag}{diag}
\DeclareMathOperator{\range}{range}
\DeclareMathOperator{\prox}{prox}
\DeclareMathOperator{\proj}{proj}

\DeclareMathOperator{\spanset}{span}

\DeclareMathOperator*{\argmin}{arg\,min}
\DeclareMathOperator*{\Argmin}{Arg\,min}

\DeclareMathOperator*{\Argmax}{Arg\,max}

\newcommand{\bg}{\mathrm{D}}
\newcommand{\R}{\mathbb{R}}
\newcommand{\Sym}{\mathbb{S}}
\newcommand{\ones}{\mathbf{1}}
\newcommand{\zeros}{\mathbf{0}}

\newcommand{\lminp}{\lambda_{\min}^+}
\newcommand{\smax}{\sigma_{\max}}
\newcommand{\smin}{\sigma_{\min}}
\newcommand{\sminp}{\sigma_{\min}^+}

\newcommand{\cA}{\mathcal{A}}

\newcommand{\cC}{\mathcal{C}}

\newcommand{\cK}{\mathcal{K}}
\newcommand{\cL}{\mathcal{L}}
\newcommand{\cM}{\mathcal{M}}

\newcommand{\cO}{\mathcal{O}}

\newcommand{\cR}{\mathcal{R}}
\newcommand{\cS}{\mathcal{S}}

\newcommand{\cX}{\mathcal{X}}
\newcommand{\cY}{\mathcal{Y}}
\newcommand{\cZ}{\mathcal{Z}}

\newcommand{\mB}{\mathbf{B}}

\newcommand{\mE}{\mathbf{E}}
\newcommand{\mF}{\mathbf{F}}
\newcommand{\mG}{\mathbf{G}}

\newcommand{\mI}{\mathbf{I}}
\newcommand{\mJ}{\mathbf{J}}

\newcommand{\mO}{\mathbf{O}}
\newcommand{\mP}{\mathbf{P}}
\newcommand{\mQ}{\mathbf{Q}}

\newcommand{\mW}{\mathbf{W}}

\makeatletter
\newcommand{\annotatehypertarget}[1]{\Hy@raisedlink{\hypertarget{#1}{}}}
\makeatother
\newcounter{annotatecount}
\newcounter{annotateidx}
\newcounter{annotatejdx}
\newcounter{annotatelabelcount}
\newcounter{annotateglobalindex}
\newcommand{\atran}[2]{\stepcounter{annotatecount}\overset{\text{\annotatehypertarget{\alph{annotatecount}\theannotateglobalindex}{(\hyperlink{desc\alph{annotatecount}\theannotateglobalindex}{\alph{annotatecount}})}}}{#1}\csgdef{annotatedescription\theannotatecount}{#2}}
\newcommand{\aeq}[1]{\atran{=}{#1}}
\newcommand{\aleq}[1]{\atran{\leq}{#1}}
\newcommand{\ageq}[1]{\atran{\geq}{#1}}
\newcommand{\agreater}[1]{\atran{>}{#1}}

\newcommand{\asucceq}[1]{\atran{\succeq}{#1}}
\newcommand{\apreceq}[1]{\atran{\preceq}{#1}}

\newcommand{\annotateinitused}{\setcounter{annotateidx}{0}\whileboolexpr{test{\ifnumless{\theannotateidx}{\theannotatecount}}}{\stepcounter{annotateidx}\csgdef{aused\theannotateidx}{0}}}
\newcommand{\annotategetlabels}{\setcounter{annotatejdx}{0}\setcounter{annotatelabelcount}{0}\whileboolexpr{test{\ifnumless{\theannotatejdx}{\theannotatecount}}}{\stepcounter{annotatejdx}\ifcsequal{annotatedescription\theannotateidx}{annotatedescription\theannotatejdx}{\csgdef{aused\theannotatejdx}{1}\stepcounter{annotatelabelcount}\csedef{annotatelabel\theannotatelabelcount}{\alph{annotatejdx}}}{}}}
\newcommand{\annotateprintlabels}{\setcounter{annotatejdx}{0}\whileboolexpr{test{\ifnumless{\theannotatejdx}{\theannotatelabelcount}}}{\stepcounter{annotatejdx}\ifnumequal{\theannotatejdx}{\theannotatelabelcount}{\ifnumequal{\theannotatejdx}{1}{}{~and~}}{}\annotatehypertarget{desc\csuse{annotatelabel\theannotatejdx}\theannotateglobalindex}{(\hyperlink{\csuse{annotatelabel\theannotatejdx}\theannotateglobalindex}{\csuse{annotatelabel\theannotatejdx}})}\ifnumless{\theannotatejdx}{\theannotatelabelcount}{\ifnumless{\theannotatejdx+1}{\theannotatelabelcount}{,~}{}}{}}}
\newcommand{\annotate}{\annotateinitused\setcounter{annotateidx}{0}\whileboolexpr{test{\ifnumless{\theannotateidx}{\theannotatecount}}}{\stepcounter{annotateidx}\ifcsstring{aused\theannotateidx}{0}{\ifnumequal{\theannotateidx}{1}{}{;~}\annotategetlabels\annotateprintlabels~\csuse{annotatedescription\theannotateidx}}{}}\setcounter{annotatecount}{0}\stepcounter{annotateglobalindex}}

\newcommand{\sX}{\cX}
\newcommand{\sY}{\cY}
\newcommand{\sZ}{\cZ}
\newcommand{\rng}[2]{\{#1,\ldots,#2\}}

\newcommand{\mind}[1]{\hspace{#1}&\hspace{-#1}}

\newcommand{\basis}[2]{\mathbf{e}_{#2}^{#1}}
\DeclareMathOperator{\dist}{dist}
\newcommand{\distsol}{\cR^2_{\delta_x\delta_y}}

\newcommand{\oz}{\overline{z}}

\newcommand{\hz}{\hat{z}}

\newcommand{\condx}{\kappa_x}
\newcommand{\condy}{\kappa_y}
\newcommand{\condxya}{\kappa_{xy}}
\newcommand{\mem}{\cM}
\newcommand{\memx}{\cM_x}
\newcommand{\memy}{\cM_y}
\newcommand{\memf}{\cM_{f}}
\newcommand{\memg}{\cM_{g}}
\newcommand{\memxb}{\cM_{\mB^\top}}
\newcommand{\memyb}{\cM_{\mB}}
\newcommand{\moreau}[2]{M_{#2#1}}
\newcommand{\const}{\mathrm{const}}

\begin{document}
\maketitle

\begin{abstract}
	We revisit the smooth convex-concave bilinearly-coupled saddle-point problem of the form $\min_x\max_y f(x) + \<y,\mB x> - g(y)$. In the highly specific case where each of the functions $f(x)$ and $g(y)$ is either affine or strongly convex, there exist lower bounds on the number of gradient evaluations and matrix-vector multiplications required to solve the problem, as well as matching optimal algorithms. A notable aspect of these algorithms is that they are able to attain linear convergence, i.e., the number of iterations required to solve the problem is proportional to $\log(1/\epsilon)$. However, the class of bilinearly-coupled saddle-point problems for which linear convergence is possible is much wider and can involve smooth non-strongly convex functions $f(x)$ and $g(y)$. Therefore, {\em we develop the first lower complexity bounds and matching optimal linearly converging algorithms for this problem class}. Our lower complexity bounds are much more general, but they cover and unify the existing results in the literature. On the other hand, our algorithm implements the separation of complexities, which, for the first time, enables the simultaneous achievement of both optimal gradient evaluation and matrix-vector multiplication complexities, resulting in the best theoretical performance to date.
\end{abstract}


\section{Introduction}\label{sec:intro}

In this paper, we consider the following saddle-point optimization problem with a bilinear coupling function:
\begin{equation}\label{eq:main}
	\adjustlimits\min_{x \in \sX}\max_{y \in \sY} \left[F(x, y) = f(x) + \<y,\mB x> - g(y)\right],
\end{equation}
where $\sX = \R^{d_x}$ and $\sY = \R^{d_y}$ are finite-dimensional Euclidean spaces, $\mB \in \R^{d_y\times d_x}$ is a coupling matrix, and $f(x)\colon \sX \to \R$ and $g(y)\colon \sY \to \R$ are continuous functions. We aim to solve problem~\eqref{eq:main} in the fundamental setting where both functions $f(x)$ and $g(y)$ are convex and smooth.\footnote{A function is called smooth if it is differentiable and has a Lipschitz-continuous gradient. See \Cref{sec:pre} for an equivalent formal definition.}

Saddle-point problems of the form~\eqref{eq:main} appear in various fields such as economics \citep{von1947theory}, game theory \citep{roughgarden2010algorithmic}, and statistics \citep{berger2013statistical}. Moreover, these problems have a wide range of applications in machine learning, including supervised learning \citep{zhang2017stochastic,wang2017exploiting,xiao2019dscovr}, reinforcement learning \citep{du2017stochastic}, computer vision \citep{chambolle2011first}, robust optimization \citep{ben2002robust,liu2017primal}, distributed optimization \citep{lan2020communication,scaman2018optimal,kovalev2021lower,yarmoshik2024decentralized,kovalev2024lower}, and the training of generative adversarial networks \citep{mescheder2017numerics,nagarajan2017gradient}.

\subsection{First-Order Methods and Linear Convergence}

The majority of machine learning applications of problem~\eqref{eq:main} involve high-dimensional spaces $\sX$ and $\sY$. In this scenario, the most widely used and often the only scalable optimization algorithms are {\em first-order methods}. These methods implement an iterative process to find an approximate solution to the problem using the evaluation of the gradients of the functions $f(x)$ and $g(y)$, as well as matrix-vector multiplication with the matrices $\mB$ and $\mB^\top$. More specifically, they perform iterative updates of the current estimate of the solution until it converges to the exact solution up to a given accuracy. One of the main goals of our paper is to develop efficient first-order optimization methods for solving problem~\eqref{eq:main}.

In this paper, we are interested in first-order methods for solving problem~\eqref{eq:main} that are able to achieve {\em linear convergence}. That is, we are interested in algorithms that can find an $\epsilon$-approximate solution to the problem using at most $\cO(\cK \cdot \log (1/\epsilon))$ gradient evaluations and matrix-vector multiplications, where $\cO(\cdot)$ hides universal constants, $\epsilon > 0$ is an arbitrary precision, and $\cK \geq 1$ is a constant that possibly depends on the internal properties of the problem such as condition numbers, etc.

In this work, we also intend to consider problem classes where linear convergence is possible in principle. A typical and one of the most fundamental examples of such a class is problem~\eqref{eq:main} with strongly convex functions $f(x)$ and $g(y)$. There are plenty of linearly converging first-order optimization methods in this {\em strongly-convex-strongly-concave} setting, which include the gradient descent ascent \citep{zhang2022near}, the extragradient method \citep{korpelevich1976extragradient}, and the optimistic gradient method \citep{gidel2018variational}. Moreover, there is an array of algorithms that enjoy improved, or accelerated, convergence rates \citep{kovalev2022accelerated,thekumparampil2022lifted,jin2022sharper,du2022optimal,li2023nesterov} with the help of the Nesterov momentum trick \citep{nesterov2013introductory}.
Another fundamental problem class where linear convergence is possible is the class of bilinear min-max games, which is a special case of problem~\eqref{eq:main} with affine functions $f(x)$ and $g(y)$. Such problems can be solved using the corresponding algorithms with linear convergence rates \citep{azizian2020accelerating,li2022convergence}. Finally, \citet{kovalev2022accelerated} developed an algorithm for solving the general smooth convex-concave problem~\eqref{eq:main} and provided a set of sufficient conditions under which the proposed algorithm attains linear convergence.

\subsection{Optimal Algorithms and Complexity Separation}

In this paper, we are concerned with the task of developing {\em optimal algorithms}, which is one of the ultimate goals in optimization research. This task can be divided into two key parts. The first part involves finding lower bounds on the {\em oracle complexity} of solving the optimization problem, i.e., the number of oracle calls, such as gradient evaluations or matrix-vector multiplications, required to find an approximate solution to the problem. The second part is to find optimization algorithms that match these lower bounds. Such algorithms are called optimal because their oracle complexity cannot be improved due to the lower complexity bounds. For example, in the case of bilinear min-max games, lower bounds were proposed by \citet{ibrahim2020linear} and matching optimal algorithms were developed by \citet{azizian2020accelerating,li2022convergence}.

Unfortunately, apart from the case of bilinear min-max games, the question of finding optimal algorithms for solving problem~\eqref{eq:main} is far from being resolved, even in the fundamental {\em strongly-convex-strongly-concave} setting. Although separate lower bounds on the gradient evaluation and matrix-vector multiplication complexities have already been developed by \citet{nesterov2013introductory} and \citet{zhang2022lower}, respectively, the existing state-of-the-art algorithms \citep{kovalev2022accelerated,thekumparampil2022lifted,jin2022sharper,du2022optimal,li2023nesterov} cannot simultaneously reach these bounds. The main reason is that these algorithms perform the same number of evaluations of the gradients $\nabla f(x)$ and $\nabla g(y)$ and matrix-vector multiplications with the matrices $\mB$ and $\mB^\top$ at each iteration while solving the problem, whereas the lower bounds on these numbers \citep{nesterov2013introductory,zhang2022lower} can be significantly different. Thus, to reach the desired lower bounds, an optimal algorithm would have to implement {\em the separation of complexities} by skipping gradient evaluations and/or matrix-vector multiplications from time to time. \citet{borodich2023optimal,alkousa2020accelerated,sadiev2022communication,lan2021mirror} attempted to develop efficient first-order methods with the complexity separation for solving the problem; however, these algorithms are not able to achieve optimal complexities by a substantial margin, see \Cref{tab} for details.

The situation is even worse in other cases, such as the {\em strongly-convex-concave} or {\em convex-strongly-concave} settings, where only one of the functions $f(x)$ or $g(y)$ is strongly convex, or the {\em convex-concave setting}, where neither of the functions is strongly convex. To the best of our knowledge, there are no lower complexity bounds that would cover these cases, with the exception of the highly specific cases of bilinear min-max games \citep{ibrahim2020linear} and affinely constrained minimization \citep{salim2022optimal}. Therefore, the question remains unresolved as to whether the current state-of-the-art linearly converging algorithms for this setting \citep{kovalev2022accelerated,sadiev2022communication} are optimal or not.

\subsection{Main Contributions}

\definecolor{colorx}{rgb}{0,0.2,1}
\definecolor{colory}{rgb}{0,1,0}
\definecolor{colorxy}{rgb}{1,0,0}
\newcommand{\hlmath}[2]{\setlength\fboxsep{0.1em}\colorbox{#1!35}{$\displaystyle #2$}}
\newcommand{\hltextbg}[2]{\setlength\fboxsep{0.05em}\colorbox{#1!35}{#2}}
\newcommand{\hltext}[2]{{\color{#1}{#2}}}

The above discussion reveals significant gaps in the current theoretical understanding of smooth convex-concave saddle-point problems with bilinear coupling. In particular, the existing lower complexity bounds are insufficient, and the state-of-the-art optimization algorithms are limited. Summarizing these gaps leads to the following open research question:
\begin{center}
	{\em Is it possible to develop an optimal linearly converging first-order optimization method for solving the smooth convex-concave bilinearly-coupled saddle-point problem~\eqref{eq:main}?}
\end{center}
We provide a positive answer to this question and present the following key contributions:
\begin{itemize}
	\item [\bf (i)] We describe the class of smooth convex-concave saddle-point problems of the form~\eqref{eq:main} for which it is possible to achieve linear convergence. We establish {\em the first lower complexity bounds} for this class. In particular, we show that to find an $\epsilon$-approximate solution to problem~\eqref{eq:main}, any first-order optimization method requires at least {$\tilde\cO(\condx)$} evaluations of the gradient {$\nabla f(x)$}, {$\tilde\cO(\condy)$} evaluations of the gradient {$\nabla g(y)$}, and {$\tilde\cO(\condxya)$} matrix-vector multiplications with the matrices {$\mB$} and {$\mB^\top$}, where $\condx$, $\condy$, and $\condxya$ denote certain condition numbers associated with functions $f(x)$, $g(y)$, and matrix $\mB$.\footnote{Here, $\tilde\cO(\cdot)$ hides the logarithmic factor $\log(1/\epsilon)$, and universal (and possibly additive) constants. The precise definitions of $\condx$, $\condy$, and $\condxya$ are provided in \Cref{sec:pre}.}

	\item [\bf (ii)] We show that our lower complexity bounds are tight. That is, we develop {\em the first optimal algorithm} that matches these lower bounds. This algorithm implements the complexity separation, allowing us to simultaneously achieve both optimal gradient evaluation and matrix-vector multiplication complexities. To the best of our knowledge, such a result has never been established in the literature, even in the fundamental strongly-convex-strongly-concave setting.

	\item[\bf (iii)] As a side contribution, we develop a new algorithm for solving a class of composite monotone variational inequalities. Just like the current state-of-the-art method of \citet{lan2021mirror}, our algorithm implements the separation of complexities, but enjoys substantially improved convergence rates and works in a much broader range of settings. Refer to \Cref{sec:opt_alg} and \Cref{sec:disc:alg} for details.
\end{itemize}

Our lower complexity bounds are much more general than the existing lower bounds for the special cases of strongly-convex-strongly-concave \citep{nesterov2013introductory,zhang2022lower}, bilinear \citep{ibrahim2020linear}, and affinely constrained \citep{salim2022optimal} optimization. On the other hand, our lower bounds recover and provide unification of these existing results. Besides, our optimal algorithm shows the best theoretical performance ``on the market'', which, to the best of our knowledge, outclasses all existing methods in the literature, with the exception of the algorithms of \citet{azizian2020accelerating,li2022convergence} and \citet{salim2022optimal,kovalev2020optimal}, which are already optimal in the aforementioned specific cases of bilinear and affinely constrained optimization, respectively.

\section{Preliminaries}\label{sec:pre}

\subsection{Main Definitions and Assumptions}\label{sec:ass}

In this paper, we use the notation described in \Cref{sec:notation}.
Further, in this section, we provide a formal description of the assumptions that we impose on problem~\eqref{eq:main}. First, we define the (strong) convexity and smoothness properties of a differentiable function.
\begin{definition}\label{def:convex}
	A differentiable function $h(x) \colon \R^d \to \R$ is called $\mu$-strongly convex for $\mu \geq 0$ if the following inequality holds for all $x,x' \in \R^d$:
	\begin{equation}
		\bg_h(x,x') \geq \tfrac{1}{2}\cdot\mu\sqn{x - x'}.
	\end{equation}
	A differentiable function $h(x) \colon \R^d \to \R$ is called convex if the same inequality holds with $\mu = 0$.
\end{definition}
\begin{definition}\label{def:smooth}
	A differentiable convex function $h(x) \colon \R^d \to \R$ is called $L$-smooth for $L\geq 0$ if the following inequality holds for all $x_1,x_2 \in \R^d$:
	\begin{equation}
		\abs{\bg_h(x,x')} \leq \tfrac{1}{2}\cdot L\sqn{x-x'}.
	\end{equation}
\end{definition}
Next, we formalize the (strong) convexity and smoothness assumptions that we impose on functions $f(x)$ and $g(y)$ as \Cref{ass:x,ass:y}. Note that we allow the strong convexity constants $\mu_x$ and $\mu_y$ to be equal to zero, thus covering the case of non-strongly convex functions $f(x)$ and $g(y)$.
\begin{assumption}\label{ass:x}
	Function $f(x)\colon \sX \to \R$ is $\mu_x$-strongly convex and $L_x$-smooth for $L_x > \mu_x \geq 0$.
\end{assumption}
\begin{assumption}\label{ass:y}
	Function $g(y)\colon \sY \to \R$ is $\mu_y$-strongly convex and $L_y$-smooth for $L_y > \mu_y \geq 0$.
\end{assumption}
Finally, the next \Cref{ass:xy} describes the spectral properties of the coupling matrix $\mB$.
\begin{assumption}\label{ass:xy}
	There exist constants $L_{xy} > \mu_{xy},\mu_{yx} \geq 0$, such that
	\begin{align*}
		\mu_{xy}^2
		 & \leq
		\begin{cases}
			\lambda_{\min}^+(\mB^\top\mB) & \text{if }\nabla f(x) \in \range \mB^\top \text{ for all } x \in \sX \\
			\lambda_{\min}(\mB^\top\mB)   & \text{otherwise}
		\end{cases},
		\\
		\mu_{yx}^2
		 & \leq
		\begin{cases}
			\lambda_{\min}^+(\mB\mB^\top) & \text{if }\nabla g(y) \in \range \mB \text{ for all } y \in \sY
			\\ \lambda_{\min}(\mB\mB^\top) & \text{otherwise}
		\end{cases},
		\\
		L_{xy}^2
		 & \geq
		\lambda_{\max}(\mB^\top\mB) = \lambda_{\max}(\mB\mB^\top),
	\end{align*}
	Additionally, we assume that if $\mu_{xy} > 0$ and $\mu_{yx} > 0$, then $\mu_{xy} = \mu_{yx}$.
\end{assumption}
Further, to shorten the notation, we gather all the parameters defined in \Cref{ass:x,ass:y,ass:xy} into a single vector $\pi = (L_x,L_y,L_{xy},\mu_x,\mu_y,\mu_{xy},\mu_{yx}) \in \Pi$, where $\Pi \subset \R_+^7$ denotes the parameter set.

\subsection{Key Assumption for Linear Convergence}\label{sec:linear}

As discussed in \Cref{sec:intro}, we are interested in algorithms for solving problem~\eqref{eq:main} that exhibit linear convergence. We introduce the key \Cref{ass:xyxy}, which will enable us to establish linear lower complexity bounds and devise optimal linearly converging algorithms.
\begin{assumption}\label{ass:xyxy}
	Under \Cref{ass:x,ass:y,ass:xy}, the following inequality holds:
	\begin{equation}\label{eq:xyxy}
		\min\{\delta_x,\delta_y\} > 0,
		\quad\text{where}\quad
		\delta_x = \mu_{x} + \mu_{xy}^2/L_y
		\quad\text{and}\quad
		\delta_y = \mu_y + \mu_{yx}^2/L_x.
	\end{equation}
\end{assumption}
To better understand this assumption, consider the standard primal and dual reformulations of problem~\eqref{eq:main}, which are given as follows:
\begin{equation}\label{eq:pd}
	\min_{x \in \sX}\left[P(x) = f(x) + g^*(\mB x)\right],
	\qquad\quad
	\max_{y \in \sY}\left[D(y) = -g(y) - f^*(-\mB^\top y)\right].
\end{equation}
One can show that the primal objective function $P(x)$ and the dual objective function $-D(y)$ satisfy the {\em quadratic growth} condition \citep{anitescu2000degenerate,karimi2016linear} with constants $\delta_x$ and $\delta_y$, respectively. This fact provides a good starting point for understanding why linear convergence is plausible under \Cref{ass:xyxy}. On the other hand, \citet{kovalev2022accelerated} showed that this assumption is sufficient for developing a linearly converging algorithm. Moreover, in \Cref{sec:lb}, we obtain \Cref{thm:lb_sublinear}, which implies that \Cref{ass:xyxy} is also necessary for achieving linear convergence, thus making it both a {\em necessary and sufficient} condition.

We also need to characterize the linear convergence rates of the first-order methods that we consider in this paper. Such rates are typically expressed via the {\em condition numbers} associated with a given optimization problem. Consequently, we define the following condition numbers for problem~\eqref{eq:main}:
\begin{equation}
	\label{eq:kappa}
	\begin{aligned}
		\condx   & = L_x/\delta_x,                
		         &
		\condy   & = L_y/\delta_y,                
		         &
		\condxya & = L_{xy}^2/(\delta_x\delta_y). 
	\end{aligned}
\end{equation}
The condition numbers $\condx$ and $\condy$ correspond to the functions $f(x)$ and $g(y)$, respectively. These can be seen as extensions of the standard condition numbers $L_{x}/\mu_x$ and $L_y/\mu_y$, which are commonly used in smooth and strongly convex optimization \citep{nesterov2013introductory}. Similarly, the condition number $\condxya$ associated with the bilinear coupling term is a generalization of the standard condition number ${L_{xy}^2}/{(\mu_x\mu_y)}$, which is widespread in strongly-convex-strongly-concave saddle-point optimization \citep{zhang2022lower,ibrahim2020linear}.

Further, we would like to ensure that the condition numbers defined in \cref{eq:kappa} are lower-bounded by some small universal constants. This is achieved by the following additional \Cref{ass:Pi} on the parameter set $\Pi$. It allows us to avoid addressing some corner cases where the condition numbers are small, which are neither theoretically nor practically interesting. It should be noted that \Cref{ass:Pi} does not impose any fundamental restrictions;\footnote{In particular, it is always possible to increase the smoothness constants $L_x$, $L_y$, and $L_{xy}$ to satisfy \Cref{ass:Pi}.} it is merely introduced to streamline our complex theoretical findings.
\begin{assumption}\label{ass:Pi}
	For all $\pi \in \Pi$ the following additional constraints are satisfied:
	\begin{equation}
		L_{x} > 4\mu_x,\;\;
		L_{y} > 4\mu_y,\;\;
		L_{xy} > 18\max\{\mu_{xy},\mu_{yx},\sqrt{\mu_x\mu_y}\},\;\;
		\sqrt{L_xL_y} > 4\max\{\mu_{xy},\mu_{yx}\}.
	\end{equation}
\end{assumption}

\subsection{Structure of the Solution Set}

In this paper, we denote the solution set of the saddle-point problem~\eqref{eq:main} as $\cS \subset \sX \times \sY$. Under \Cref{ass:x,ass:y,ass:xy}, $(x^*,y^*) \in \cS$ if and only if the following first-order optimality conditions hold:
\begin{equation}\label{eq:opt}
	\nabla f(x^*) + \mB^\top y^* = 0,
	\quad
	\nabla g(y^*) - \mB x^* = 0.
\end{equation}
Moreover, under \Cref{ass:xyxy}, the solution set is always non-empty and has an affine structure, as indicated by \Cref{lem:sol}.
\begin{lemma}\label{lem:sol}
	Under \Cref{ass:x,ass:y,ass:xy,ass:xyxy}, the solution set $\cS$ of problem~\eqref{eq:main} is nonempty and is given as
	\begin{equation}
		\textstyle
		\cS = \cS_x \times \cS_y,
		\quad\text{where}\quad
		\cS_x = \Argmin_{x \in \sX} P(x),
		\quad
		\cS_y = \Argmax_{y \in \sY} D(y).
	\end{equation}
	Moreover, the primal and dual solution sets $\cS_x \subset \sX$ and $\cS_y \subset \sY$ have the following affine structure:
	\begin{equation}
		\cS_x = x^* + \begin{cases}
			\{0\}    & \mu_x > 0        \\
			\ker \mB & \text{otherwise}
		\end{cases},
		\qquad
		\cS_y = y^* + \begin{cases}
			\{0\}         & \mu_y > 0        \\
			\ker \mB^\top & \text{otherwise}
		\end{cases},
	\end{equation}
	where $(x^*,y^*) \in \cS$ is an arbitrary solution to problem~\eqref{eq:main}.
\end{lemma}
We also define a weighted squared distance function $\distsol(x,y)$ as follows:
\begin{equation}\label{eq:distsol}
	\distsol(x,y) = \delta_x \dist(x;\cS_x) + \delta_y\dist(y;\cS_y).
\end{equation}
We are going to use this function to measure the quality of a given approximate solution to problem~\eqref{eq:main} in both lower complexity bounds and the convergence analysis of optimal algorithms.

\section{Lower Complexity Bounds}\label{sec:lb}

\subsection{First-Order Saddle-Point Optimization Methods}
In this section, we present lower bounds on the number of gradient evaluations and matrix-vector multiplications required to solve problem~\eqref{eq:main}. These lower bounds apply to a specific class of algorithms that we refer to as {\em first-order saddle-point optimization methods}. A formal description of this class is provided in \Cref{def:alg}. This definition is mostly inspired by the common {\em linear span assumption} \citep{nesterov2013introductory,zhang2022lower,ibrahim2020linear}. However, the standard existing definitions focus only on iteration complexity. This is insufficient in our case, as we need to derive more specific lower bounds on the numbers of gradient evaluations and matrix-vector multiplications. Therefore, in \Cref{def:alg}, we introduce a continuous execution time parameter $\tau \geq 0$ and assume that the evaluation of the gradients $\nabla f(x)$ and $\nabla g(y)$ takes time $\tau_f$ and $\tau_g$, respectively, while matrix-vector multiplication with matrices $\mB$ and $\mB^\top$ takes time $\tau_{\mB}$. A similar approach was previously used in distributed optimization by \citet{scaman2017optimal,scaman2018optimal,kovalev2024lower}, where they had to ensure a distinction between communication and local computation complexities.

\begin{definition}\label{def:alg}
	An algorithm is called a first-order saddle-point optimization method with gradient computation times $\tau_f,\tau_g > 0$, and matrix-vector computation time $\tau_\mB > 0$, if it satisfies the following constraints:
	\begin{enumerate}
		\item[\bf (i)] {\bf Memory.} At any time $\tau \geq 0$, the algorithm maintains a memory, which is represented by a set $\mem(\tau) = \memx(\tau)\times \memy(\tau)$, where $\memx(\tau)\subset\sX$ and $\memy(\tau) \subset \sY$. The memory can be updated by computing the gradients $\nabla f(x)$ and $\nabla g(y)$, and by performing matrix-vector multiplications with matrices $\mB$ and $\mB^\top$. This is represented by the following inclusions:
		      \begin{equation}
			      \memx(\tau)
			      \subset
			      \memf(\tau) \cup \memxb(\tau),\quad
			      \memy(\tau)
			      \subset
			      \memg(\tau) \cup \memyb(\tau),
		      \end{equation}
		      where sets $\memf(\tau)$, $\memg(\tau)$, $\memxb(\tau)$, and $\memyb(\tau)$ are defined below.
		\item[\bf (ii)] {\bf Gradient computation.} At any time $\tau \geq 0$, the algorithm can update the memory by computing the gradients $\nabla f(x)$ and $\nabla g(y)$, which take time $\tau_f$ and $\tau_g$, respectively. That is, for all $\tau \geq 0$, sets $\memf(\tau) \subset \sX$ and $\memg(\tau) \subset \sY$ are defined as follows:
		      \begin{equation}
			      \begin{aligned}
				      \memf(\tau) & = \begin{cases}
					                      \spanset(\{x, \nabla f(x) : x \in \memx(\tau - \tau_f)\}) & \tau \geq \tau_f \\
					                      \varnothing                                               & \tau < \tau_f
				                      \end{cases}, \\
				      \memg(\tau) & = \begin{cases}
					                      \spanset(\{y, \nabla g(y) : y \in \memy(\tau - \tau_g)\}) & \tau \geq \tau_g \\
					                      \varnothing                                               & \tau < \tau_g
				                      \end{cases}.
			      \end{aligned}
		      \end{equation}
		\item[\bf (iii)] {\bf Matrix-vector multiplication.} At any time $\tau \geq 0$, the algorithm can update the memory by performing matrix-vector multiplication with matrices $\mB$ and $\mB^\top$, which takes time $\tau_\mB$. That is, for all $\tau \geq 0$, sets $\memxb(\tau) \subset \sX$ and $\memyb(\tau) \subset \sY$ are defined as follows:
		      \begin{equation}
			      \begin{aligned}
				      \memxb(\tau) & = \begin{cases}
					                       \spanset(\{x, \mB^\top y : (x,y) \in \mem(\tau - \tau_\mB)\}) & \tau \geq \tau_\mB \\
					                       \varnothing                                                   & \tau < \tau_\mB
				                       \end{cases}, \\
				      \memyb(\tau) & = \begin{cases}
					                       \spanset(\{\mB x, y  : (x,y) \in \mem(\tau - \tau_\mB)\}) & \tau \geq \tau_\mB \\
					                       \varnothing                                               & \tau < \tau_\mB
				                       \end{cases}.
			      \end{aligned}
		      \end{equation}
		\item[\bf (iv)] {\bf Initialization and output.} At time $\tau=0$, the algorithm must initialize the memory with the zero vector, that is, $\memx(0) = \{0\}, \memy(0) = \{0\}$. At any time $\tau \geq 0$, the algorithm must specify a single output vector from the memory, $(x_o(\tau), y_o(\tau)) \in \mem(\tau)$.
	\end{enumerate}
\end{definition}

\subsection{Lower Bounds}\label{sec:lb_thm}

In this section, we present our lower complexity bounds. We start with \Cref{thm:lb_sublinear}, which shows that it is not possible to obtain a linearly converging algorithm for solving problem~\eqref{eq:main} if \Cref{ass:xyxy} does not hold. The proof can be found in \Cref{sec:proof:lb_sublinear}.
This theorem indicates that there exists a specific "hard" instance of problem~\eqref{eq:main}, such that any first-order saddle-point optimization method fails to converge in terms of the distance to the solution set and converges sublinearly in terms of the primal-dual gap. It is important to clarify that the main purpose of \Cref{thm:lb_sublinear} is to demonstrate the impossibility of attaining linear convergence in general if \Cref{ass:xyxy} does not hold, rather than to provide tight lower complexity bounds for this setting. Thus, we leave the investigation of the general case of problem~\eqref{eq:main} under \Cref{ass:x,ass:y,ass:xy} for future work.
\begin{theorem}\label{thm:lb_sublinear}
	Let $\pi \in \Pi$, $R_x > 0$, and $\epsilon > 0$ be arbitrary parameters, distance, and precision, respectively. Suppose that \Cref{ass:xyxy} does not hold, i.e., without loss of generality, $\delta_x = 0$. There exists a problem~\eqref{eq:main} satisfying \Cref{ass:x,ass:y,ass:xy} with parameters $\pi$, such that $\dist(0;\cS_x) = R_x$, and for any first-order saddle-point optimization method and execution time $\tau > 0$, the following inequality holds:
	\begin{equation}
		\dist^2(x_o(\tau);\cS_x) > \tfrac{1}{8}R_x^2.
	\end{equation}
	Moreover, to reach precision on the primal-dual gap $P(x_o(\tau)) - D(y_o(\tau)) < \epsilon$ by any first-order saddle-point optimization method, the execution time $\tau$ must satisfy the following inequality:
	\begin{equation}
		\tau \geq \Omega\left(\tau_f \cdot {\dist(0;\cS_x)\sqrt{L_x/\epsilon}}\right).
	\end{equation}
\end{theorem}

Now, we are ready to present lower complexity bounds for problem~\eqref{eq:main} under \Cref{ass:x,ass:y,ass:xy,ass:xyxy,ass:Pi} in \Cref{thm:lb}. The proof can be found in \Cref{sec:proof:lb}.
The lower bound on the total execution time $\tau$ in \cref{eq:lb} contains the terms \hltextbg{colorx}{$\sqrt{\condx}\log \tfrac{cR^2}{\epsilon}$}, \hltextbg{colory}{$\sqrt{\condy}\log \tfrac{cR^2}{\epsilon}$}, and \hltextbg{colorxy}{$\sqrt{\condxya}\log \tfrac{cR^2}{\epsilon}$}. These terms can be respectively interpreted as gradient evaluation complexities with respect to the gradients \hltext{colorx}{$\nabla f(x)$} and \hltext{colory!80!black}{$\nabla g(y)$}, and matrix-vector multiplication complexity with respect to the matrices \hltext{colorxy}{$\mB$} and \hltext{colorxy}{$\mB^\top$}, as they are respectively multiplied by the corresponding times \hltext{colorx}{$\tau_f$}, \hltext{colory!80!black}{$\tau_g$}, and \hltext{colorxy}{$\tau_{\mB}$}. In addition, these complexities are proportional to the logarithmic factor $\log(1/\epsilon)$, making them linear as we previously discussed. Furthermore, in \Cref{sec:alg}, we will prove the tightness of the lower bound by developing a matching optimal algorithm.
\begin{theorem}\label{thm:lb}
	Under \Cref{ass:Pi}, let $\pi \in \Pi$, $R > 0$, and $\epsilon > 0$ be arbitrary parameters, distance, and precision, respectively. Suppose that \Cref{ass:xyxy} hold. There exists a problem~\eqref{eq:main} satisfying \Cref{ass:x,ass:y,ass:xy} with parameters $\pi$, such that $\distsol(0,0) = R^2$, and to reach precision $\distsol(x_o(\tau),y_o(\tau)) < \epsilon$ by any first-order saddle-point optimization method, the execution time $\tau$ must satisfy the following inequality:
	\begin{equation}\label{eq:lb}
		\tau \geq \Omega\left(
		\tau_f \cdot \hlmath{colorx}{\sqrt{\condx}\log \tfrac{cR^2}{\epsilon}}
		+
		\tau_g \cdot \hlmath{colory}{\sqrt{\condy}\log \tfrac{cR^2}{\epsilon}}
		+
		\tau_{\mB} \cdot \hlmath{colorxy}{\sqrt{\condxya}\log \tfrac{cR^2}{\epsilon}}
		\right),
	\end{equation}
	where $c > 0$ is a universal constant.
\end{theorem}
The result in \Cref{thm:lb} has two important merits. First, this lower bound is tight, which we prove by developing a matching optimal algorithm in \Cref{sec:alg}. Second, by making an appropriate restriction of the parameter set $\Pi$, we can recover the existing lower complexity bounds for the important and fundamental special cases of strongly-convex-strongly-concave saddle-point optimization \citep{zhang2022lower,nesterov2013introductory}, bilinear saddle-point optimization \citep{ibrahim2020linear}, and strongly convex minimization with affine constraints \citep{salim2022optimal}. On the other hand, our result applies to an arbitrary choice of parameters $\pi \in \Pi$. Therefore, \Cref{thm:lb} and our definition of the condition numbers $\condx$, $\condy$, and $\condxya$ in \Cref{sec:pre} provide unification and substantial generalization of the existing results. See \Cref{sec:disc:lb} for additional discussion.

\section{Optimal Algorithm}\label{sec:alg}
\newcommand{\pit}[4]{p_{#1}^{#3;#4}(#2)}
\newcommand{\pitpr}[4]{\hat{p}_{#1}^{#3;#4}(#2)}
\newcommand{\zit}[2]{z^{#1}_{#2}}
\newcommand{\hzit}[2]{\oz^{#1}_{#2}}
\newcommand{\git}[3]{\Delta_{#1}^{#2;#3}}
\newcommand{\gitq}[2]{\Delta_{Q_#1}^{#1;#2}}
\newcommand{\Hit}[3]{H_{#1}^{#2;#3}}
\newcommand{\Lit}[3]{L_{#1}^{#2;#3}}
\newcommand{\Mit}[3]{M_{#1}^{#2;#3}}
\newcommand{\procname}{RecursiveProcedure}
\newcommand{\zin}{z_{\text{in}}}
\newcommand{\xin}{x_{\text{in}}}
\newcommand{\yin}{y_{\text{in}}}
\newcommand{\zout}{z_{\text{out}}}
\newcommand{\xout}{x_{\text{out}}}
\newcommand{\yout}{y_{\text{out}}}

\subsection{Monotone Variational Inequalities}
\newcommand{\Zcon}{\cC_z}
In this section, we develop an optimal algorithm for solving problem~\eqref{eq:main}. To do this, we consider a more general monotone variational inequality problem, which is given as follows:
\begin{equation}\label{eq:vi}
	\text{find}\;\;
	z^* \in \Zcon
	\;\;\text{such that}\;\;
	p(z^*) - p(z) + \<Q(z),z^* - z> \leq 0
	\;\;\text{for all}\;\;
	z \in \Zcon,
\end{equation}
where $\Zcon$ is a closed and convex subset of the finite-dimensional Euclidean space $\sZ = \R^{d_z}$, and differentiable convex function $p(z)\colon \sZ \to \R$ and continuous monotone operator $Q(z) \colon \sZ \to \sZ$ have the following finite-sum structures:
\begin{equation}\label{eq:vi_sum}
	p(z) = \sum_{i=1}^{n} p_i(z),\quad
	Q(z) = \sum_{i=1}^{n} Q_i(z),\quad
	\text{where}\quad
	p_i(z)\colon \sZ \to \R,\quad
	Q_i(z) \colon \sZ \to \sZ.
\end{equation}
Vector $z^*$ defined in \cref{eq:vi} is often called a {\em weak} solution to the monotone variational inequality. In the setting of this paper, it is equivalent to the {\em strong} solution\footnote{Vector $z^* \in \Zcon$ is a strong solution to the variational inequality if $p(z^*) - p(z) + \<Q(z^*),z^* - z> \leq 0$ for all $z \in \Zcon$.}; refer to \citet{kinderlehrer2000introduction} for details.

Further, we assume that the gradients $\nabla p_i(z)$ and operators $Q_i(z)$ are monotone and Lipschitz with respect to the norm $\sqn{\cdot}_\mP$, where $\mP \in \Sym^{d_z}_{++}$. These assumptions are commonly used in the literature and are formalized through the following \Cref{def:monotone,def:Lipschitz} and \Cref{ass:p,ass:Q}. Note that \Cref{ass:p} implies that each function $p_i(z)$ is convex and smooth.
\begin{definition}\label{def:monotone}
	An operator $G(x)\colon \R^d \to \R^d$ is called $\mu$-strongly monotone with respect to the norm $\norm{\cdot}_\mP$ for $\mu \geq 0$ if the following inequality holds for all $x,x' \in \R^d$:
	\begin{equation}
		\<G(x) - G(x'),x - x'> \geq \mu\sqn{x-x'}_\mP.
	\end{equation}
	An operator $G(x)\colon \R^d \to \R^d$ is called monotone if the same inequality holds with $\mu = 0$.
\end{definition}
\begin{definition}\label{def:Lipschitz}
	An operator $G(x)\colon \R^d \to \R^d$ is called $M$-Lipschitz with respect to the norm $\norm{\cdot}_\mP$ for $M \geq 0$ if the following inequality holds for all $x,x' \in \R^d$:
	\begin{equation}
		\norm{G(x) - G(x')}_{\mP^{-1}} \leq M\norm{x-x'}_{\mP}.
	\end{equation}
\end{definition}
\begin{assumption}\label{ass:p}
	For all $1 \leq i  \leq n$, the gradient $\nabla p_i(z)$ is monotone and $L_i$-Lipschitz w.r.t.~$\norm{\cdot}_\mP$.
\end{assumption}
\begin{assumption}\label{ass:Q}
	For all $1 \leq i  \leq n$, operator $Q_i(z)$ is monotone and $M_i$-Lipschitz w.r.t.~$\norm{\cdot}_\mP$.
\end{assumption}


\subsection{Optimal Sliding Algorithm for Monotone Variational Inequalities}
\begin{algorithm}[t]
    \caption{}
    \label{alg}
    \begin{algorithmic}[1]
        \State {\bf input:} $\zin \in \Zcon$
        \State {\bf parameters:}
        $\{\alpha_t\}_{t=0}^{\infty}\subset \R_{++}$,\;
        $\{L_i\}_{i=1}^n, \{M_i\}_{i=1}^n \subset \R_{+}$,\;
        $\{T_k\}_{k=1}^n \in \{1,2,\ldots\}$,\;
        $\mP \in \Sym_{++}^{d_z}$
        \For{$i=1,\ldots,n$}
        \State $\zit{i}{0,\ldots,0} = \zin$
        \label{line:z_in}
        \State $p_i^0(z) \equiv p_i(z)$
        \label{line:p_in}
        \EndFor
        \State $\zout = \Call{\procname}{1}$
        \label{line:z_out}
        \State {\bf output:} $\zout \in \sZ$
        \LComment{Auxiliary Recursive Procedure:}
        \Procedure{\procname}{$k$, $t_1,\ldots,t_{k-1}$}
        \If{$k = n+1$}
        \State \Return $\argmin_{z \in \Zcon} \sum_{i=1}^n \pit{i}{z}{n}{t_1,\ldots,t_n}$
        \label{line:argmin}
        \Else
        \State $\hzit{k}{0} = \zit{k}{t_1,\ldots,t_{k-1},0}$
        \For{$t_k = 0,\ldots,T_{k}-1$}
        \For{$i=1,\ldots,n$}
        \State $\pitpr{i}{z}{k}{t_1,\ldots,t_k} \equiv
            \begin{cases}
                \alpha_{t_k}^{-1}\pit{i}{\alpha_{t_k}z + (1-\alpha_{t_k})\hzit{k}{t_k}}{k-1}{t_1,\ldots,t_{k-1}}
                 & i \geq k \\
                \pit{i}{z}{k-1}{t_1,\ldots,t_{k-1}}
                 & i < k
            \end{cases}$
        \label{line:p_prime}
        \State $\Lit{i}{k}{t_1,\ldots,t_k} = L_i \cdot \prod_{l=1}^{k}\alpha_{t_l}$
        \label{line:L}
        \State $\Mit{i}{k}{t_1,\ldots,t_k} = M_i\cdot \prod_{l=1}^{k}(\alpha_{t_l} / \alpha_{T_l-1})$
        \label{line:M}
        \State $\Hit{i}{k}{t_1,\ldots,t_k} = \Lit{i}{k}{t_1,\ldots,t_k} + \Mit{i}{k}{t_1,\ldots,t_k}$
        \label{line:H}
        \State $\git{i}{k}{t_1,\ldots,t_k} = \nabla \pitpr{i}{\zit{k}{t_1,\ldots,t_k}}{k}{t_1,\ldots,t_k} + Q_i(\zit{k}{t_1,\ldots,t_k})$
        \label{line:grad}
        \State $\pit{i}{z}{k}{t_1,\ldots,t_k} \equiv
            \begin{cases}
                \frac{\Hit{k}{k}{t_1,\ldots,t_k}}{2}\sqn{z - \zit{k}{t_1,\ldots,t_k}}_{\mP} + \<z, \git{k}{k}{t_1,\ldots,t_k}>
                 & i = k    \\
                \pitpr{i}{z}{k}{t_1,\ldots,t_k}
                 & i \neq k
            \end{cases}$
        \label{line:p}
        \EndFor
        \State $\zit{k}{t_1,\ldots,t_{k-1},t_k+1/2} = \Call{\procname}{k+1,t_1,\ldots,t_k}$
        \label{line:call}
        \State $\hzit{k}{t_k + 1} = \alpha_{t_k}\zit{k}{t_1,\ldots,t_{k-1},t_k+1/2} + (1-\alpha_{t_k})\hzit{k}{t_k}$
        \label{line:hz}
        \State $\gitq{k}{t_1,\ldots,t_k} = Q_k(\zit{k}{t_1,\ldots,t_k}) - Q_k(\zit{k}{t_1,\ldots,t_{k-1},t_k+1/2})$
        \label{line:delta_q}
        \State $\zit{k}{t_1,\ldots,t_{k-1},t_k+1} = \zit{k}{t_1,\ldots,t_{k-1},t_k+1/2} + (\Hit{k}{k}{t_1,\ldots,t_k}\mP)^{-1}\gitq{k}{t_1,\ldots,t_k}$
        \label{line:z}
        \For{$l = k+1,\ldots,n$}
        \State $\zit{l}{t_1,\ldots,t_{k-1},t_k+1,0,\ldots,0} = \zit{l}{t_1,\ldots,t_{k-1},t_k,T_{k+1},0,\ldots,0}$
        \label{line:z_table}
        \EndFor
        \EndFor
        \State \Return $\hzit{k}{T_k}$
        \label{line:return}
        \EndIf
        \EndProcedure
    \end{algorithmic}
\end{algorithm}

Now, we are ready to present our new algorithm for solving the variational inequality problem~\eqref{eq:vi}. One of the key ideas behind the development of this algorithm is our new perspective on the celebrated accelerated gradient method of \citet{nesterov2013introductory}. In particular, a single step of this algorithm, applied to minimizing an $L$-smooth convex function $h(z)\colon \sZ \to \R$, can be seen as applying a single step of the standard gradient method to the function $h_t(z)$ with the fixed stepsize $1/L$, where $h_t(z)\colon \sZ \to \R$ is defined as follows:
\begin{equation}\label{eq:h_t}
	h_t(z) = \alpha_t^{-2}h(\alpha_t z + (1-\alpha_t)\oz^t),\quad\text{where}\quad \oz^t \in \sZ.
\end{equation}
Indeed, the stepsize $1/L$ is suitable, since one can show that function $h_t(z)$ is $L$-smooth as well. Hence, using the standard recursion for the gradient descent, for all $z \in \sZ$, we obtain the inequality
\begin{equation}
	\tfrac{1}{2}L\sqn{z^{t+1} - z} + h_t(z^{t+1}) \leq \tfrac{1}{2}L\sqn{z^{t} - z} + h_t(z).
\end{equation}
Next, we can define $\oz^{t+1} = \alpha_t z^{t+1} + (1-\alpha_t)\oz^t$, and use the definition of function $h_t(z)$ in \cref{eq:h_t}, the convexity of function $h(z)$, and the definition of $\alpha_t$ in \cref{eq:alpha_t}.\footnote{From the definition of $\alpha_t$ in \cref{eq:alpha_t}, it follows that $\alpha_t \in (0,1]$ and $\alpha_t^{-2} = \alpha_{t-1}^{-2} + \alpha_t^{-1}$.} This gives the following recursion:
\begin{equation}
	\tfrac{1}{2}L\sqn{z^{t+1} - z} +  \alpha_t^{-2}[h(\oz^{t+1}) - h(z)] \leq \tfrac{1}{2}L\sqn{z^{t} - z} + \alpha_{t-1}^{-2}[h(\oz^t) - h(z)],
\end{equation}
which implies the desired optimal rate $h(\oz^t) - \min_z h(z) = \cO(LR^2/t^2)$, where $R > 0$ is the initial distance to the solution. Overall, the derivations above offer a vast simplification compared to the standard proof of \citet{nesterov2013introductory}.

Inspired by the {\em sliding} algorithm of \citet{lan2021mirror,lan2016accelerated,kovalev2022optimal}, we apply a series of transformations of the form~\eqref{eq:h_t} to functions $p_i(z)$ in a recursive fashion. This leads, subject to some additional details, to \Cref{alg} for solving problem~\eqref{eq:vi}. Moreover, using the considerations above, we obtain the key theoretical result in \Cref{thm:alg}. The proof can be found in \Cref{sec:proof:alg}.
\begin{theorem}\label{thm:alg}
	Let \Cref{ass:p,ass:Q} hold, where $M_i,L_i \geq 0$ and  $M_i + L_i > 0$. Let $\alpha_t$ be defined recursively as follows:
	\begin{equation}\label{eq:alpha_t}
		\alpha_0 = 1,\qquad
		\alpha_{t+1} = 2\cdot\big(1 + \sqrt{\rule{0pt}{8pt}\smash{1 + 4/\alpha_t^2}}\big)^{-1}
		\;\;\text{for}\;\; t\geq 1.
	\end{equation}
	Then $\zout \in \Zcon$ and the following inequality holds for all $z \in \Zcon$:
	\begin{equation}\label{eq:alg}
		p(\zout) - p(z) + \<Q(z),\zout - z>
		\leq
		\sum_{i=1}^n\left(\frac{4^{i}L_i}{\prod_{j=1}^{i}T_j^2} + \frac{2^{i}M_i}{\prod_{j=1}^{i}T_j}\right)\cdot \tfrac{1}{2}\sqn{\zin - z}_\mP.
	\end{equation}
\end{theorem}
Furthermore, we can reorder functions $p_i(z)$ and operators $Q_i(z)$ in ascending order of the values of the Lipschitz constants $L_i$ and $M_i$, which leads to the complexity result in \Cref{cor:alg}. The proof can be found in \Cref{sec:proof:alg2}.
\begin{corollary}\label{cor:alg}
	Under the conditions of \Cref{thm:alg}, to ensure the following inequality
	\begin{equation}\label{eq:vi_gap}
		p(\zout) - p(z) + \<Q(z),\zout - z> \leq \epsilon n \sqn{\zin - z}_\mP
	\end{equation}
	for all $z \in \Zcon$ and $\epsilon > 0$, it is sufficient to perform no more than
	\begin{equation}
		6^n \cdot \max\big\{\sqrt{L_i/\epsilon}, M_i / \epsilon,1\big\}
	\end{equation}
	computations of the gradient $\nabla p_i(z)$ and operator $Q_i(z)$ for $1 \leq i \leq n$.
\end{corollary}
Using \Cref{cor:alg}, we can show that \Cref{alg} achieves the optimal complexity separation for solving the variational inequality problem~\eqref{eq:vi} as long as $n = \cO(1)$. Consequently, \Cref{alg} theoretically outperforms the existing state-of-the-art algorithm of \citet{lan2021mirror}, which is designed for the case $n=2$ with additional restrictions.\footnote{The algorithm of \citet{lan2021mirror} works in the case $n=2$, where $Q_1(z) \equiv \zeros_{d_z}$ and $p_2(z) \equiv 0$.} See \Cref{sec:disc:alg} for details.

\subsection{Application to the Main Saddle-Point Problem}\label{sec:opt_alg}

In this section, we show how to adapt \Cref{alg} to solve the main problem~\eqref{eq:main} and reach the lower complexity bounds in \Cref{thm:lb}. To do this, we consider a special instance of problem~\eqref{eq:vi}, where $n = 3$, $\sZ = \Zcon = \sX\times \sY$, operators $Q_i(z) = Q_i(x,y)$ are defined as follows:
\begin{equation}\label{eq:Q}
	Q_1(x,y) = 0,\qquad
	Q_2(x,y) = 0,\qquad
	Q_3(x,y) =
	\begin{bNiceMatrix}[c]
		\mO_{d_x} & \mB^\top  \\
		-\mB      & \mO_{d_y}
	\end{bNiceMatrix}
	\begin{bNiceMatrix}
		x \\ y
	\end{bNiceMatrix},
\end{equation}
and functions $p_i(z) = p_i(x,y)$ are defined as follows:
\begin{equation}\label{eq:p}
	\begin{aligned}
		p_1(x,y) & = f(x), \\
		p_2(x,y) & = g(y),
	\end{aligned}
	\qquad
	p_3(x,y) = \frac{\beta_x}{2}\sqn{\mB x - \nabla g(\yin)} + \frac{\beta_y}{2}\sqn{\mB^\top y + \nabla f(\xin)},
\end{equation}
where $\zin = (\xin,\yin) \in \sZ$ and $\beta_x,\beta_y \geq 0$ are defined as follows:
\begin{equation}\label{eq:beta_xy}
	\beta_x = 1/(4L_y),\qquad
	\beta_y = 1/(4L_x).
\end{equation}
We also define matrix $\mP \in \Sym_{++}^{d_z}$ as the following diagonal matrix:
\begin{equation}\label{eq:P}
	\mP = \diag(\delta_x \mI_{d_x},\delta_y \mI_{d_y}).
\end{equation}
We apply \Cref{alg} to solve this problem instance and, using \Cref{thm:alg} and \Cref{cor:alg}, we obtain the following result in \Cref{thm:opt_alg}. The proof can be found in \Cref{sec:proof:opt_alg}.
\begin{theorem}\label{thm:opt_alg}
	Under \Cref{ass:x,ass:y,ass:xy,ass:xyxy,ass:Pi},
	functions $p_i(z)$, operators $Q_i(z)$, and matrix $\mP$ defined in \cref{eq:P,,eq:p,,eq:Q} satisfy the conditions of \Cref{thm:alg} with the following parameters:
	\begin{equation}\label{eq:LM}
		L_1 = \condx,\quad
		L_2 = \condy,\quad
		L_3 = \condxya,\quad
		M_1 = M_2 = 0,\quad
		M_3 = \sqrt{\condxya}.
	\end{equation}
	Moreover, the input $\zin = (\xin,\yin)$ and the output $\zout=(\xout,\yout)$ of \Cref{alg} satisfy the inequality $\Psi(\zout) \leq \tfrac{2}{3}\Psi(\zin)$ as long as the numbers of inner iterations $\{T_i\}_{i=1}^3$ are chosen according to \Cref{cor:alg}. Here, the Lyapunov function $\Psi(x,y)$ is defined as follows:
	\begin{equation}\label{eq:Psi}
		\Psi(z) = \Psi(x,y) = \distsol(x,y) + 12\bg_f(x,x^*) + 12\bg_g(y,y^*),
	\end{equation}
	where $z^* = (x^*,y^*) =  \proj_\cS(\zin) = \proj_\cS(\zout)$.
\end{theorem}
\Cref{thm:opt_alg} implies that we can reduce the value of the Lyapunov function $\Psi(x,y)$ defined in \cref{eq:Psi} by a constant factor with a single run of \Cref{alg}. Hence, we can apply the standard restarting technique to this algorithm and obtain the complexity result in \Cref{cor:opt_alg}. The proof can be found in \Cref{sec:proof:opt_alg2}.
\begin{corollary}\label{cor:opt_alg}
	Under the conditions of \Cref{thm:opt_alg}, to reach precision $\distsol(x,y) \leq \epsilon$, it is sufficient to perform \hltextbg{colorx}{$\cO\big(\sqrt{\condx}\log \tfrac{cR^2}{\epsilon}\big)$}, \hltextbg{colory}{$\cO\big(\sqrt{\condy}\log \tfrac{cR^2}{\epsilon}\big)$}, and \hltextbg{colorxy}{$\cO\big(\sqrt{\condxya}\log \tfrac{cR^2}{\epsilon}\big)$} computations of the gradients \hltext{colorx}{$\nabla f(x)$} and \hltext{colory!80!black}{$\nabla g(y)$}, and matrix-vector multiplications with the matrices \hltext{colorxy}{$\mB$} and \hltext{colorxy}{$\mB^\top$}, respectively. Here, $R^2 = \distsol(0,0)$ is the initial distance, $c = 1 + 12\condx+12\condy$, and $\epsilon \in (0,R^2)$.
\end{corollary}
The complexity result in \Cref{cor:opt_alg} matches the lower complexity bounds in \Cref{thm:lb} up to universal and/or additive constants. Hence, this result is optimal. Moreover, to the best of our knowledge, this result theoretically outperforms all existing state-of-the-art algorithms, including the algorithms of \citet{kovalev2022accelerated,li2023nesterov,jin2022sharper,du2022optimal,thekumparampil2022lifted,borodich2023optimal,alkousa2020accelerated,sadiev2022communication,chambolle2011first}. See \Cref{sec:disc:opt_alg} for additional discussion.

\section{Additional Discussions}

\subsection{Lower Complexity Bounds}\label{sec:disc:lb}
As mentioned in \Cref{sec:lb_thm}, the lower complexity bound in \Cref{thm:lb} recovers several important problem classes, which are special instances of problem~\eqref{eq:main}. We can recover this problem classes by imposing additional constraints on the parameter set $\Pi$.
\begin{itemize}
	\item[\bf (i)] The class of smooth strongly-convex-strongly-concave saddle-point optimization problems corresponds to the constraint
	      \begin{equation}
		      \mu_{xy} = \mu_{yx} = 0.
	      \end{equation}
	      In this case, the lower bound in \cref{eq:lb} becomes the following:
	      \begin{equation}
		      \tilde{\Omega}\left(\tau_f\cdot\sqrt{\frac{L_x}{\mu_x}}  + \tau_g\cdot \sqrt{\frac{L_y}{\mu_y}} + \tau_\mB\cdot \frac{L_{xy}}{\sqrt{\mu_x\mu_y}}\right),
	      \end{equation}
	      where $\tilde{\Omega}(\cdot)$ hides universal constants and logarithmic factors.
	      This result recovers the existing lower complexity bounds of \citet{zhang2022lower,nesterov2013introductory}.
	\item[\bf (ii)]
	      The class of bilinear saddle-point optimization problems is obtained by choosing
	      \begin{equation}
		      L_x = L_y = 0,\quad
		      \mu_x = \mu_y = 0,\quad
		      \tau_f = \tau_g = 0,\quad
		      \mu_{xy} = \mu_{yx} > 0,
	      \end{equation}
	      and the lower bound in \cref{eq:lb} turns into the following:
	      \begin{equation}\label{eq:lb_bilinear}
		      \tilde{\Omega}\left(\tau_\mB\cdot \frac{L_{xy}}{\mu_{xy}}\right).
	      \end{equation}
	      This result recovers the existing lower complexity bound of \citet{ibrahim2020linear}. Note that strictly speaking, we cannot choose $L_x = L_y = 0$ due to \Cref{ass:Pi}. However, this is not an issue, because \Cref{ass:Pi} allows us to choose arbitrary $L_x, L_y > 0$ such that $\sqrt{L_xL_y} = 5\mu_{xy}$ and still obtain the lower bound~\eqref{eq:lb_bilinear} from \cref{eq:lb}. In addition, as mentioned in \Cref{sec:pre}, this assumption is not a fundamental restriction but is rather used to avoid covering uninteresting corner cases in our theoretical proofs.
	\item[\bf (iii)] The class of smooth strongly convex optimization problems with affine constraints is obtained by choosing
	      \begin{equation}
		      L_y = \mu_y = 0,\quad
		      \mu_{xy} = 0,\quad
		      \tau_g = 0.
	      \end{equation}
	      In this case, the lower bound in \cref{eq:lb} becomes the following:
	      \begin{equation}\label{eq:lb_affine}
		      \tilde{\Omega}\left(\tau_f\cdot\sqrt{\frac{L_x}{\mu_x}} + \tau_\mB\cdot \frac{L_{xy}}{\mu_{yx}}\sqrt{\frac{L_x}{\mu_x}}\right),
	      \end{equation}
	      which recovers the existing result of \citet{salim2022optimal}. Similarly to the previous case {\bf (ii)}, we can choose $L_y = {17\mu_{yx}^2}/{L_x}$ instead of $L_y = 0$ to satisfy \Cref{ass:Pi} and still obtain the lower bound~\eqref{eq:lb_affine}.

	      It is worth mentioning the work of \citet{ouyang2021lower}, who offer sublinear lower complexity bounds for this problem class. However, their result does not contradict ours since they consider the case $\mu_{yx} = \delta_y = 0$, i.e., \Cref{ass:xyxy} does not hold. It is also important to highlight that affinely constrained optimization problems, where $\mu_{yx} = 0$, hold limited interest. Indeed, in this setting, it is typically assumed that $\mu_{yx}^2 = \lminp(\mB\mB^\top)$, which, by definition, is always nonzero.
\end{itemize}

\subsection{Algorithm for Solving the Variational Inequality Problem}\label{sec:disc:alg}

We compare \Cref{alg} for solving problem~\eqref{eq:vi} with the algorithm of \citet{lan2021mirror} in the case $n=2$, where function $p_2(z)$ and operator $Q_1(z)$ are zero. Note that this is the main problem setting used by \citet{lan2021mirror}. Let $R \geq 0$ be the following distance parameter associated with the constraint set $\Zcon$:
\begin{equation}
	R = \sup_{z \in \Zcon}\norm{z - \zin}_\mP.
\end{equation}
We compare the numbers of evaluations of the gradient $\nabla p_1(z)$ and operator $Q_2(z)$ required by both algorithms to find a vector $\zout \in \Zcon$ that satisfies the following accuracy criterion:
\begin{equation}
	\sup_{z \in \Zcon}p(\zout) - p(z) + \<Q(z),\zout - z> \leq \epsilon,
\end{equation}
where $\epsilon > 0$ is an arbitrary precision. Note that the parameter $R$ is finite only if the constraint set is bounded. However, we can easily tackle this issue by following the standard approach and replacing the constraint set $\Zcon$ with its intersection with the ball $\{ z \in \sZ : \norm{z-\zin}_{\mP} \leq D\}$, where $D > 0$ is a positive parameter. Refer, for instance, to \citet{nesterov2007dual}.

\begin{table}[t]
    \centering
    \caption{Comparison of \Cref{alg} with the algorithm of \citet{lan2021mirror}.}
    \label{tab2}
    \begin{NiceTabular}[cell-space-limits=0.1em]{|c|c|c|}
        \CodeBefore
        \rowcolor{gray!20}{4}
        \Body
        \hline
        \Block{2-1}{\bf Algorithm}
         & \Block{1-2}{\bf Complexity}
         &
        \\
        \cline{2-3}
         & $\nabla p_1(z)$
         & $Q_2(z)$
        \\
        \hline
        \citet{lan2021mirror}
         & $\displaystyle \cO\left(\sqrt{\frac{L_1R^2}{\epsilon}}\right)$
         & $\displaystyle \cO\left(\sqrt{\frac{L_1R^2}{\epsilon}} + \frac{M_2R^2}{\epsilon}\right)$
        \\
        \hline
        \begin{varwidth}{\textwidth}
            \centering
            \Cref{alg}\\
            (\Cref{cor:alg})
        \end{varwidth}
         & $\displaystyle \cO\left(\sqrt{\frac{L_1R^2}{\epsilon}}\right)$
         & $\displaystyle \cO\left(\frac{M_2R^2}{\epsilon}\right)$
        \\
        \hline
    \end{NiceTabular}
\end{table}

The comparison of \Cref{alg} with the algorithm of \citet{lan2021mirror} is summarized in \Cref{tab2}. One can observe that the theoretical complexities of these algorithms coincide up to universal constants when $\sqrt{L_1/\epsilon} \leq M_2 / \epsilon$. However, \Cref{alg} can significantly outperform the algorithm of \citet{lan2021mirror} in the case where $\sqrt{L_1/\epsilon} \gg M_2 / \epsilon$. It is important to highlight that this case is necessary to consider, as it plays an essential role in the task of developing an optimal algorithm for solving the main problem~\eqref{eq:main}. In addition, the algorithm of \citet{lan2021mirror} only works in the case $n=2$ with the additional restrictions described above. On the other hand, using the result in \Cref{cor:alg}, it is easy to verify that our \Cref{alg} can achieve the optimal complexity separation as long as $n = \cO(1)$.

\subsection{Optimal Algorithm for Solving the Main Problem}\label{sec:disc:opt_alg}

{\bf Comparison with existing results.}
The theoretical complexity of \Cref{alg} with restarting, applied to solve the smooth bilinearly-coupled saddle-point optimization problem~\eqref{eq:main}, is established in \Cref{cor:opt_alg} and is proven to be optimal due to the lower complexity bounds in \Cref{thm:lb}. We compare this result with the theoretical complexities of the existing state-of-the-art linearly converging first-order methods. These include the algorithms for the strongly-convex-strongly-concave case \citep{kovalev2022accelerated,li2023nesterov,jin2022sharper,du2022optimal,thekumparampil2022lifted,borodich2023optimal,chambolle2011first,alkousa2020accelerated}, the strongly-convex-concave case \citep{kovalev2022accelerated,sadiev2022communication}, and the convex-concave case \citep{kovalev2022accelerated}. This comparison is summarized in \Cref{tab}. One can observe that our optimal result is substantially better compared to the existing algorithms. It is also worth highlighting that the complexity of our algorithm matches the complexities of the algorithms of \citet{salim2022optimal} and \citet{azizian2020accelerating,li2022convergence}, which are optimal in the case of affinely constrained minimization and bilinear saddle-point optimization, respectively, as discussed in \Cref{sec:disc:lb}.

{\bf Auxiliary Variational Inequality Subproblem in \Cref{sec:opt_alg}.}
From the optimality conditions~\eqref{eq:opt}, it is easy to observe that the main problem~\eqref{eq:main} is equivalent to the following variational inequality problem, i.e., finding $z^* \in \sZ$ such that
\begin{equation}
	(p_1(z^*) + p_2(z^*)) - (p_1(z) + p_2(z)) + \<Q_3(z),z^* - z> \leq 0
	\;\;\;\text{for all}\;\;\;
	z \in \sZ,
\end{equation}
where functions $p_1(z),p_2(z)$ and operator $Q_3(z)$ are defined in \cref{eq:p,eq:Q}. This problem matches the special instance of the monotone variational inequality problem defined in \Cref{sec:opt_alg}, with the only difference being the addition of the quadratic function $p_3(z)$ defined in \cref{eq:p}. The additional quadratic regularization terms in the function $p_3(z)$ help to achieve the optimal linear convergence rates in all cases where $\delta_x > 0$ and $\delta_y > 0$, even when $\mu_x = 0$ and/or $\mu_y = 0$. Moreover, these terms do not break the convergence to the solution $z^*$ of the original problem~\eqref{eq:main}. Indeed, it is easy to show that $\nabla p_3(z^*)$ converges to zero as long as $\zin$ converges to $z^*$. Refer to the proof of \Cref{thm:opt_alg} in \Cref{sec:proof:alg} for more details.

\newcommand{\tabopt}[1]{
    \begin{varwidth}{\textwidth}
        \centering
        {\bf Optimal}\tabularnote{Lower complexity bounds are established in \Cref{thm:lb}. These bounds are matched by \Cref{alg}, which is established by \Cref{thm:opt_alg} and \Cref{cor:opt_alg}.}\\
        {\bf(this paper)}#1
    \end{varwidth}
}
\begin{table}[t]
    \centering
    \caption{Comparison of the optimal complexity of \Cref{alg} developed in this paper (\Cref{thm:lb}, \Cref{cor:opt_alg}) with the existing state-of-the-art linearly-converging algorithms in the strongly-convex-strongly-concave, strongly-convex-concave, and convex-concave settings.}
    \label{tab}
    \begin{NiceTabular}[cell-space-limits=0.2em,notes/detect-duplicates,notes/style={\color{red}(\arabic{#1})}]{|c|c|c|c|}
        \hline
        \Block{2-1}{\bf Algorithm}
         & \Block{1-3}{{\bf Complexity}\tabularnote{For brevity, we omit universal constants and logarithmic factors such as $\log \frac{1}{\epsilon}$.}}
        \\\cline{2-4}
         & $\nabla f(x)$\cellcolor{colorx!25}
         & $\nabla g(y)$\cellcolor{colory!25}
         & $\mB$ and $\mB^\top$\cellcolor{colorxy!25}
        \\\hline
        \Block{1-4}{\bf Strongly-convex-strongly-concave case ($\mu_x,\mu_y > 0$ and $\mu_{xy} = \mu_{yx} = 0$)}
        \\\hline
        \begin{varwidth}{\textwidth}
            \centering
            \citet{kovalev2022accelerated,li2023nesterov}
            \\
            \citet{jin2022sharper,du2022optimal}
            \\
            \citet{thekumparampil2022lifted}
        \end{varwidth}
         & \Block{1-3}{$\sqrt{\frac{L_x}{\mu_x}} + \sqrt{\frac{L_y}{\mu_y}} + \frac{L_{xy}}{\sqrt{\mu_x\mu_y}}$}
        \\\hline
        \citet{borodich2023optimal}
         & \Block{1-2}{$\sqrt{\frac{L_x}{\mu_x}} + \sqrt{\frac{L_y}{\mu_y}}$}
         &
         & $\sqrt{\frac{L_x}{\mu_x}} + \sqrt{\frac{L_y}{\mu_y}} + \frac{L_{xy}}{\sqrt{\mu_x\mu_y}}$
        \\\hline
        \citet{chambolle2011first}
         & \Block{1-2}{N/A\tabularnote{Requires computation of the proximal operators of functions $f(x)$ and/or $g(y)$.}}
         &
         & $\frac{L_{xy}}{\sqrt{\mu_x\mu_y}}$
        \\\hline
        \citet{alkousa2020accelerated}
         & $\sqrt{\frac{L_x}{\mu_x}}$
         & $\frac{L_{xy}\sqrt{L_y}}{\sqrt{\mu_x}\mu_y}$
         & $\frac{\sqrt{L_{xy}^3}}{\sqrt{\mu_x}\mu_y}$
        \\\hline
        \tabopt{\tabularnote{Here, the lower bounds were also established by \citet{zhang2022lower,nesterov2013introductory}.}}
         & \Block[fill=colorx!25]{1-1}{$\sqrt{\frac{L_x}{\mu_x}}$}
         & \Block[fill=colory!25]{1-1}{$\sqrt{\frac{L_y}{\mu_y}}$}
         & \Block[fill=colorxy!25]{1-1}{$\frac{L_{xy}}{\sqrt{\mu_x\mu_y}}$}
        \\\hline
        \Block{1-4}{{\bf Strongly-convex-concave case ($\mu_x,\mu_{yx} > 0$ and $\mu_{xy} = \mu_{y} = 0$)}\tabularnote{This case is symetric to the convex-strongly-concave case, which we omit for brevity.}}
        \\\hline
        \citet{kovalev2022accelerated}
         & \Block{1-3}{$\frac{L_{xy}}{\mu_{yx}}\sqrt{\frac{L_x}{\mu_x}} + \frac{\sqrt{L_xL_y}}{\mu_{yx}} + \frac{L_{xy}^2}{\mu_{yx}^2}$}
        \\\hline
        \citet{sadiev2022communication}
         & $\frac{L_{xy}}{\mu_{yx}}\sqrt[4]{\frac{L_x^{3}}{\mu_x^{3}}} + \frac{L_{xy}^2}{\mu_{yx}^2}\sqrt[4]{\frac{L_x}{\mu_x}}$
         & N/A\tabularnote{Requires computation of the proximal operators of functions $f(x)$ and/or $g(y)$.}
         & $\frac{L_{xy}}{\mu_{yx}}\sqrt{\frac{L_x}{\mu_x}} + \frac{L_{xy}^2}{\mu_{yx}^2}$
        \\\hline
        \tabopt{}
         & \Block[fill=colorx!25]{1-1}{$\sqrt{\frac{L_x}{\mu_x}}$}
         & \Block[fill=colory!25]{1-1}{$\frac{\sqrt{L_xL_y}}{\mu_{yx}}$}
         & \Block[fill=colorxy!25]{1-1}{$\frac{L_{xy}}{\mu_{yx}}\sqrt{\frac{L_x}{\mu_x}}$}
        \\\hline
        \Block{1-4}{\bf Convex-concave case ($\mu_{xy}=\mu_{yx} > 0$ and $\mu_{x} = \mu_{y} = 0$)}
        \\\hline
        \citet{kovalev2022accelerated}
         & \Block{1-3}{$\frac{L_{xy}\sqrt{L_xL_y}}{\mu_{xy}^2} + \frac{L_{xy}^2}{\mu_{xy}^2}$}
        \\\hline
        \tabopt{}
         & \Block[fill=colorx!25]{1-1}{$\frac{\sqrt{L_xL_y}}{\mu_{xy}}$}
         & \Block[fill=colory!25]{1-1}{$\frac{\sqrt{L_xL_y}}{\mu_{xy}}$}
         & \Block[fill=colorxy!25]{1-1}{$\frac{L_{xy}\sqrt{L_xL_y}}{\mu_{xy}^2}$}
        \\\hline
    \end{NiceTabular}
\end{table}

\newpage

\bibliographystyle{apalike}
\bibliography{references}

\appendix
\part*{Appendix}

\section{Notation}\label{sec:notation}

In this paper, we are going to use the following notations:
$\Sym_{\vphantom{+}}^p$ and $\Sym_{++}^p$ denote the sets of $p\times p$ symmetric and symmetric positive definite matrices, respectively;
$\mI_p$ denotes the $p\times p$ identity matrix,
$\mJ_{p\times q}$ and $\mO_{p\times q}$ denote the $p\times q$ all-ones and all-zeros matrices, respectively,
$\mJ_p = \mJ_{p\times p}$ and $\mO_p = \mO_{p\times p}$;
$\basis{p}{j} \in \R^p$ denotes the $j$-th unit basis vector,
$\ones_{p} = (1,\ldots,1) \in \R^p$,
$\zeros_{p} = (0,\ldots,0) \in \R^p$.
In addition,
$\norm{\cdot}$ denotes the standard Euclidean norm of a vector, and $\<\cdot,\cdot>$ denotes the standard scalar product of two vectors,
$\norm{\cdot}_{\mP} = \norm{\mP^{\frac{1}{2}}(\cdot)}$ and $\<\cdot,\cdot>_\mP = \<\mP(\cdot),\cdot>$ denote the weighted Euclidean norm and scalar product, respectively, where $\mP \in \Sym_{++}^p$;
$\lambda_{\min}(\cdot)$, $\lambda_{\min}^+(\cdot)$, and $\lambda_{\max}(\cdot)$ denote the smallest, smallest positive, and largest eigenvalues of a matrix, respectively;
$\smax(\cdot)$ and $\sminp(\cdot)$ denote the largest and smallest positive singular values of a given matrix.

For a nonempty closed convex set $\cA \subset \R^d$ and a vector $x \in \R^d$, we define the standard distance function $\dist(x;\cA)$ as follows:
\begin{equation}
	\dist(x;\cA) = \argmin_{x'\in\cA}\norm{x-x'}.
\end{equation}
For a differentiable function $h(x)$, we denote the Bregman divergence associated with $h(x)$ as $\bg_h(x,x')$, which is defined as follows:
\begin{equation}
	\bg_h(x,x') = h(x) - h(x') - \<\nabla h(x'), x - x'>.
\end{equation}
For a proper, closed, and convex function $h(x)$, we denote its Fenchel conjugate as $h^*(x)$, its Moreau envelope as $\moreau{h}{\lambda}(x)$, and its proximal operator as $\prox_{\lambda h}(x)$. These are respectively defined as follows:
\begin{equation}
	\begin{aligned}
		h^*(x)                 & = \sup_{x'} \big(\<x,x'> - h(x')\big),
		\\
		\moreau{h}{\lambda}(x) & = \min_{x'} \left(h(x') + \frac{1}{2\lambda}\sqn{x'-x}\right),
		\\
		\prox_{\lambda h}(x)   & = \argmin_{x'} \left(h(x') + \frac{1}{2\lambda}\sqn{x'-x}\right).
	\end{aligned}
\end{equation}

\section{Proof of Lemma~\ref{lem:sol}}

We define linear spaces $\cL_x \subset \sX$ and $\cL_y \subset \sY$ as follows:
\begin{equation}\label{eq:cLxy}
    \cL_x = \begin{cases}
        \{0\}    & \mu_x > 0        \\
        \ker \mB & \text{otherwise}
    \end{cases},
    \qquad
    \cL_y = \begin{cases}
        \{0\}         & \mu_y > 0        \\
        \ker \mB^\top & \text{otherwise}
    \end{cases}.
\end{equation}
One can show that the following identity holds:
\begin{equation}\label{eq:F_diff}
    F(x+d_x, y+d_y) = F(x,y)
    \quad\text{for all}\quad
    (x,y) \in \sX \times \sY
    \quad\text{and}\quad
    (d_x,d_y) \in \cL_x \times \cL_y.
\end{equation}
Indeed, for the saddle part, we obviously have $\<y+d_y,\mB(x+d_x)> = \<y,\mB x>$.
Furthermore, we can show that $f(x+d_x) = f(x)$ and $g(y+d_y) = g(y)$. Indeed, \Cref{ass:xy} and \cref{eq:cLxy} imply $\nabla f(x) \in \cL_x^\perp$ and $\nabla g(y) \in \cL_y^\perp$. Hence, we obtain
\begin{align*}
    f(x + d_x) - f(x) & = \int_{0}^{1}\<\nabla f(x+d_x\cdot t),d_x> \mathrm{d}t = 0, \\
    g(y + d_y) - g(y) & = \int_{0}^{1}\<\nabla g(y+d_y\cdot t),d_y> \mathrm{d}t = 0,
\end{align*}
which conludes the proof of \cref{eq:F_diff}. In addition, it is easy to show that
\begin{equation}
    \label{eq:dom_dual}
    \dom f^*(\cdot) \subset \cL_x^\perp
    \quad\text{and}\quad
    \dom g^*(\cdot) \subset \cL_y^\perp.
\end{equation}

Consider the following saddle-point problem:
\begin{equation}\label{eq:main2}
    \min_{x \in \cL_x^\perp}\max_{y \in \cL_y^\perp} F(x,y).
\end{equation}
We can show that this problem has a unique solution $(x^*,y^*) \in \cL_x^\perp  \times \cL_y^\perp$, which, together with \cref{eq:F_diff}, implies \Cref{lem:sol}. Let us further prove this statement.

One can show that function $P(x)$ is strongly convex on $\cL_x^\perp$. Indeed, if $\mu_x >0$ this statement is obvious. Otherwise, \Cref{ass:xyxy} implies $\mu_{xy} >0$, which in turn implies the strong convexity of function $g^*(\mB x)$ on $\cL_x^\perp$ thanks to the strong convexity of function $g^*(y)$. The strong convexity of function $g^*(y)$ is implied by the smoothness of function $g(y)$.

Next, we show that $\dom P(\cdot) \neq \varnothing$, which immediately implies $\dom P(\cdot) \cap \cL_x^\perp \neq \varnothing$, thanks to \cref{eq:F_diff}. Indeed, if $\mu_y > 0$ function $g^*(y)$ is smooth, which implies $\dom P(\cdot) = \sX$. Otherwise, \Cref{ass:xyxy} implies $\mu_{yx} > 0$, which in turn implies $\nabla g(y) \in \range \mB$ for all $y \in \sY$. Hence, there exists $x \in \cL_x^\perp$ such that $\mB x = \nabla g(y)$ for some $y \in \sY$. On the other hand, $\nabla g(y) \in \dom g^*(\cdot)$, which implies $\mB x \in \dom g^*(\cdot)$ and $x \in \dom P(\cdot)$.

The strong convexity of $P(x)$ on $\cL_x^\perp$ and the fact that $\dom P(\cdot) \cap \cL_x^\perp \neq \varnothing$ imply that there exists a unique solution $x^* \in \cL_x^\perp$ to the following problem:
\begin{equation}
    \min_{x \in \cL_x^\perp} P(x).
\end{equation}
Similarly, there exists a unique solution $y^* \in \cL_y^\perp$ to the following problem:
\begin{equation}
    \label{eq:dual2}
    \max_{y \in \cL_y^\perp} D(y).
\end{equation}
Moreover, vectors $x^*$ and $y^*$ are solutions to the primal and dual problems in \cref{eq:pd}, respectively, thanks to \cref{eq:F_diff}.

Let $h(x) = g^*(\mB x)$. Vector $x^*$ is a solution to the problem $\min_{x \in \sX} [f(x) + h(x)]$. Hence, standard theory implies $-\nabla f(x^*) \in \partial h(x^*)$, or
\begin{equation}
    h(x) \geq h(x^*) - \<\nabla f(x^*),x - x^*>
    \quad\text{for all}\quad
    x \in \sX.
\end{equation}
From this inequality, for arbitrary $x \in x^* + \ker \mB$, we obtain
\begin{align*}
    0 \leq  \<\nabla f(x^*),x - x^*>,
\end{align*}
which implies $\nabla f(x^*) \in \range \mB^\top$. Hence, there exists $y \in \cL_y^\perp$ such that $\nabla f(x^*) = -\mB^\top y$, which for all $x \in \sX$, implies
\begin{align*}
    h(x) \geq h(x^*) + \<y,\mB x - \mB x^*>.
\end{align*}
Hence, for all $z \in \range \mB$, we obtain
\begin{align*}
    g^*(z) \geq g^*(\mB x^*) + \<y,z - \mB x^*>.
\end{align*}
In addition, this inequality holds for all $z \in \sY$ due to \cref{eq:dom_dual}. Hence, $y \in \partial g^*(\mB x^*)$, which implies $\mB x^* = \nabla g(y)$. We also have $x^* \in \partial f^*(-\mB^\top y)$. Hence, $y$ is a solution of problem~\eqref{eq:dual2}, which implies $y = y^*$. It remains to observe, that $(x^*,y^*)$ satisfies the first-order optimality conditions~\eqref{eq:opt}. Hence, the strong duality holds in both problems~\eqref{eq:main} and~\eqref{eq:main2}, which concludes the proof. \qed
\newpage
\section{Proof of Theorem~\ref{thm:lb_sublinear}}\label{sec:proof:lb_sublinear}

Note that the condition $\delta_x = 0$ implies $\mu_x=0$ and $\mu_{xy} = 0$.
Consider the following special instance of problem~\eqref{eq:main}:
\begin{equation}
    \adjustlimits\min_{u_x \in \R^{d_1}}\min_{v_x \in \R^{d_2}}\max_{y \in \R^{d_2}}
    f(u_x) + \frac{L_x}{2}\sqn{v_x} + \mu_{yx}\<y,v_x> - \frac{\mu_y}{2}\sqn{y},
\end{equation}
where $f(x)\colon \R^{d_1}\to \R$ is the $L_x$-smooth function proposed by \citet[Theorem~2.1.7]{nesterov2013introductory}. This problem has a single solution $(u_x^*,0,0) \in \R^{d_1}\times\R^{d_2}\times\R^{d_2}$, where $u_x^* = \argmin_{u_x \in \R^{d_1}} f(u_x)$. Moreover, the primal-dual gap is lower-bounded as follows:
\begin{align*}
    P(u_x,v_x) - D(y) \geq f(u_x) - f(u_x^*).
\end{align*}
Thus, the statement of \Cref{thm:lb_sublinear} trivially follows from Theorem~2.1.7 of \citet{nesterov2013introductory}.\qed
\newcommand{\hx}{\hat{x}}
\newcommand{\hy}{\hat{y}}

\section{Proof of Theorem~\ref{thm:lb}}\label{sec:proof:lb}

The proof of \Cref{thm:lb} relies on the following \Cref{lem:lb_fg,,lem:lb_B1,,lem:lb_B2}.
The proof of \Cref{lem:lb_fg} is available in \cref{sec:proof:lb_fg}.
The proof of \Cref{lem:lb_B1} is available in \cref{sec:proof:lb_B1}.
The proof of \Cref{lem:lb_B2} is available in \cref{sec:proof:lb_B2}.
\begin{lemma}\label{lem:lb_fg}
    Under conditions of \Cref{thm:lb}, the execution time can be lower-bounded as follows:
    \begin{equation}
        \tau \geq \Omega\left(
        \tau_f\cdot\sqrt{\condx}\log \tfrac{cR^2}{\epsilon}
        +\tau_g\cdot\sqrt{\condy}\log \tfrac{cR^2}{\epsilon}
        \right).
    \end{equation}
\end{lemma}
\begin{lemma}\label{lem:lb_B1}
    Under conditions of \Cref{thm:lb}, let $\mu_y, \mu_{yx} > 0$ and $\mu_{yx}^2 \geq \mu_x\mu_y$.
    Then the execution time can be lower-bounded as follows:
    \begin{equation}
        \tau \geq \Omega\left(
        \tau_\mB\cdot\sqrt{\condxya}\log \tfrac{cR^2}{\epsilon}
        \right).
    \end{equation}
\end{lemma}
\begin{lemma}\label{lem:lb_B2}
    Under conditions of \Cref{thm:lb}, let $\mu_x,\mu_y > 0$ and $\mu_x\mu_y \geq \max\{\mu_{xy}^2,\mu_{yx}^2\}$. Then the execution time can be lower-bounded as follows:
    \begin{equation}
        \tau \geq \Omega\left(
        \tau_\mB\cdot \sqrt{\condxya}\log \tfrac{cR^2}{\epsilon}
        \right).
    \end{equation}
\end{lemma}

It remains to obtain the lower bound
$\Omega\big(\tau_\mB\cdot \sqrt{\condxya}\log \tfrac{cR^2}{\epsilon}\big)$ without the additional assumptions that were made in \Cref{lem:lb_B1,lem:lb_B2}. It can be done by considering the following special cases:
\begin{itemize}
    \item[\bf (i)] {\bf Case $\mu_x = \mu_y = 0$.} In this case, we have $\mu_{xy} = \mu_{yx} > 0$ due to \Cref{ass:xy,ass:xyxy}. We can replace $\mu_y = 0$ with a very small value $\mu_y > 0$ and apply \Cref{lem:lb_B1} to obtain the desired result.
    \item[\bf (ii)] {\bf Case $\mu_x = 0$ and $\mu_y > 0$.}
          \begin{enumerate}
              \item [\bf (ii.a)] {\bf Case $\mu_{yx} > 0$.} We can apply \Cref{lem:lb_B1}.
              \item [\bf (ii.b)] {\bf Case $\mu_{yx} = 0$.} This case is symmetric to case {\bf (iii.b)}.
          \end{enumerate}
    \item[\bf (iii)] {\bf Case $\mu_x > 0$ and $\mu_y = 0$.} 
          \begin{enumerate}
              \item [\bf (iii.a)] {\bf Case $\mu_{xy} > 0$.} This case is symmetric to case {\bf (ii.a)}.
              \item [\bf (iii.b)] {\bf Case $\mu_{xy} = 0$.} In this case, we have $\mu_{yx} > 0$ due to \Cref{ass:xyxy}. We can replace $\mu_y = 0$ with a very small value $\mu_y > 0$ and apply \Cref{lem:lb_B1} to obtain the desired result.
          \end{enumerate}
    \item[\bf (iv)] {\bf Case $\mu_x > 0$ and $\mu_y > 0$.}
          \begin{itemize}
              \item[\bf (iv.a)] {\bf Case $\mu_x\mu_y \geq \max\{\mu_{xy}^2,\mu_{yx}^2\}$.} We can apply \Cref{lem:lb_B2}.
              \item[\bf (iv.b)] {\bf Case $\mu_x\mu_y < \mu_{yx}^2$.} We can apply \Cref{lem:lb_B1}.
              \item[\bf (iv.c)] {\bf Case $\mu_x\mu_y < \mu_{xy}^2$.} This case is symmetric to case {\bf (iv.b)}.
          \end{itemize}
\end{itemize}
This concludes the proof.\qed




\newpage

\section{Proof of Lemma~\ref{lem:lb_fg}}\label{sec:proof:lb_fg}

{\bf Case $\mu_{xy} > 0$.}
We consider a special instance of problem~\eqref{eq:main}, where $\sX = \sY = \R^d$, functions $f(x)$ and $g(y)$ and matrix $\mB$ are defined as follows:
\begin{equation}
    \begin{aligned}
        f(x) & =\frac{\mu_x}{2}\sqn{x} + \frac{L_x - \mu_x}{2}\sqn{\mF x} - A \<\basis{d}{1},x>, \\
        g(y) & =\frac{L_y}{2}\sqn{y}, \qquad \mB = \mu_{xy}\mI_d,
    \end{aligned}
\end{equation}
where matrix $\mF \in \R^{(d-1)\times d}$ is defined as follows:
\begin{equation}
    \mF = \frac{1}{2}\begin{bNiceMatrix}
        1 & -1                   \\
          & \Ddots & \Ddots      \\
          &        & 1      & -1
    \end{bNiceMatrix}.
\end{equation}
This problem instance has a unique solution $(x^*,y^*) \in \sX \times \sY$, which is given as follows:
\begin{equation}
    \begin{aligned}
        x^*  = \argmin_{x \in \sX} \frac{\delta_x}{2}\sqn{x} + \frac{L_x - \mu_x}{2}\sqn{\mF x} - A \<\basis{d}{1},x>,\qquad
        y^*  = \frac{\mu_{xy}}{L_y} x^*.
    \end{aligned}
\end{equation}
The rest of the proof is similar to the proofs of \Cref{lem:lb_B1,lem:lb_B2}.

{\bf Case $\mu_{xy} = 0$.} In this case, we assume $\mu_{yx} > 0$, otherwise we can use the proof of the previous case.
We consider a special instance of problem~\eqref{eq:main}, where $\sX = \R^{d+1}$ and $\sY = \R$, functions $f(x)$ and $g(y)$ and matrix $\mB$ are defined as follows:
\begin{equation}
    \begin{aligned}
        f(x) & = f(u_x,v_x) = \frac{\mu_x}{2}\sqn{x} + \frac{L_x - \mu_x}{2}\sqn{\mF u_x} - A \<\basis{d}{1},u_x>, \\
        g(y) & =\frac{L_y}{2}\sqn{y},
        \qquad \mB = \mu_{yx}\begin{bNiceMatrix}
                                 0 & \Cdots & 0 & 1
                             \end{bNiceMatrix},
    \end{aligned}
\end{equation}
where $x = (u_x,v_x)$, $u_x \in \R^d$, $v_x \in \R$, and matrix $\mF \in \R^{(d-1)\times d}$ is defined as follows:
\begin{equation}
    \mF = \frac{1}{2}\begin{bNiceMatrix}
        1 & -1                   \\
          & \Ddots & \Ddots      \\
          &        & 1      & -1
    \end{bNiceMatrix}.
\end{equation}
This problem instance has a unique solution $(x^*,y^*) = (u_x^*,v_x^*,y^*) \in \sX \times \sY$, which is given as follows:
\begin{equation}
    \begin{aligned}
        u_x^*  = \argmin_{x \in \sX} \frac{\mu_x}{2}\sqn{u_x} + \frac{L_x - \mu_x}{2}\sqn{\mF u_x} - A \<\basis{d}{1},u_x>,\qquad
        v_x^* = y^* = 0.
    \end{aligned}
\end{equation}
The rest of the proof is similar to the proofs of \Cref{lem:lb_B1,lem:lb_B2}.\qed

\newpage

\section{Proof of Lemma~\ref{lem:lb_B1}}\label{sec:proof:lb_B1}
\newcommand{\tdelta}{\tilde{\delta}}

We consider the following ``hard'' instance of problem~\eqref{eq:main}:
\begin{itemize}
    \item[\bf(i)] Linear spaces $\sX$ and $\sY$ are defined as $\sX = (\R^d)^{n_x}$ and $\sY = (\R^d)^{n_y}$, where
          \begin{equation}
              n_x = 3n
              \quad\text{and}\quad
              n_y = \begin{cases}
                  3n & \mu_{xy} > 0 \\3n-1&\mu_{xy} = 0
              \end{cases},
              \quad
              n \in \{2,3,\ldots\}.
          \end{equation}
    \item[\bf(ii)] Function $f(x)$ is defined as follows:
          \begin{equation}\label{eq:f}
              f(x) = \sum_{i=1}^{n_x}f_i(x_i),
          \end{equation}
          where we use the notation $x = (x_1,\ldots,x_{3n}) \in (\R^d)^{3n}$,
          and functions $f_i(x_i)\colon\R^d \to \R$ are defined as follows:
          \begin{equation}
              f_i(z) = \begin{cases}
                  \frac{1}{2}\mu_x\sqn{z} + \tfrac{1}{2}(L_x-\tdelta_x)\sqn{\mF_1 z}
                   & i \in \rng{1}{n}     \\
                  \frac{1}{2}\tdelta_x\sqn{z}
                   & i \in \rng{n+1}{2n}  \\
                  \frac{1}{2}\tdelta_x\sqn{z} + \frac{1}{2}(L_x-\tdelta_x)\sqn{\mF_2 z} - A\<\basis{d}{1},z>
                   & i \in \rng{2n+1}{3n} \\
              \end{cases},
          \end{equation}
          where $A \in \R$ will be determined later, $\tdelta_x > 0$ is defined as follows:
          \begin{equation}\label{eq:tdelta_x}
              \tdelta_{x} = \mu_x + {4\mu_{xy}^2}/{L_y},
          \end{equation}
          and matrices $\mF_1 \in \R^{2\floor{d/2} \times d}$ and $\mF_2 \in \R^{2\floor{(d-1)/2}\times d}$ are defined as follows:
          \begin{equation}\label{eq:F}
              \begin{aligned}
                  \mF_1 & = \frac{1}{\sqrt{2}}\begin{bNiceMatrix}[margin]
                                                  1      & -1 & 0 & \Cdots &   &        \\
                                                  0      & 0  & 1 & -1     & 0 & \Cdots \\
                                                  \Vdots &    &   &        &   &
                                              \end{bNiceMatrix},     \\
                  \mF_2 & = \frac{1}{\sqrt{2}}\begin{bNiceMatrix}[margin]
                                                  0      & 1 & -1 & 0 & \Cdots              \\
                                                  0      & 0 & 0  & 1 & -1     & 0 & \Cdots \\
                                                  \Vdots &   &    &   &        &   &
                                              \end{bNiceMatrix}.
              \end{aligned}
          \end{equation}
    \item[\bf(iii)] Function $g(y)$ is defined as follows:
          \begin{equation}\label{eq:g}
              g(y) = \sum_{i=1}^{n_y}g_i(y_i),
          \end{equation}
          where we use the notation $y = (y_1,\ldots,y_{n_y}) \in (\R^d)^{n_y}$, and
          functions $g_i(y_i)\colon\R^d \to \R$ are defined as follows:
          \begin{equation}
              g_i(y_i) = \begin{cases}
                  \frac{1}{2}\tilde{L}_y\sqn{y_i} & i = 1              \\
                  \frac{1}{2}\mu_y\sqn{y_i}       & i \in \rng{2}{n_y}
              \end{cases},
              \quad\text{where}\quad
              \tilde{L}_y = \begin{cases}
                  L_y   & \mu_{xy} > 0 \\
                  \mu_y & \mu_{xy} = 0
              \end{cases}.
          \end{equation}
    \item[\bf(iv)] Matrix $\mB \in \R^{n_yd\times n_xd}$ is defined as follows:
          \begin{equation}\label{eq:B}
              \mB = \begin{cases}
                  \mE \otimes \mI_d  & \mu_{xy} > 0 \\
                  \mE' \otimes \mI_d & \mu_{xy} = 0
              \end{cases}
          \end{equation}
          where matrix $\mE \in \R^{3n\times 3n}$ is defined as follows:
          \begin{equation}\label{eq:E}
              \mE = \begin{bNiceArray}[margin,cell-space-limits=2pt]{ccc|cccc|ccc}
                  \gamma&\Cdots^{n\text{ times}}&\gamma\\
                  \hline
                  \beta&&&-\beta\\
                  &\Ddots^{n\text{ times}}&&\Vdots\\
                  &&\beta&-\beta\\
                  \hline
                  &&&\alpha&-\alpha\\
                  &&&&\Ddots&\Ddots^{n-1\text{ times}}\\
                  &&&&&\alpha&-\alpha\\
                  \hline
                  &&&&&&-\beta&\beta\\
                  &&&&&&\Vdots&&\Ddots^{n\text{ times}}\\
                  &&&&&&-\beta&&&\beta\\
              \end{bNiceArray},
          \end{equation}
          and matrix $\mE' \in \R^{(3n-1)\times 3n}$ is defined as follows:
          \begin{equation}
              \label{eq:E'}
              \mE' = \begin{bNiceArray}[margin,cell-space-limits=2pt]{ccc|cccc|ccc}
                  \beta&&&-\beta\\
                  &\Ddots^{n\text{ times}}&&\Vdots\\
                  &&\beta&-\beta\\
                  \hline
                  &&&\alpha&-\alpha\\
                  &&&&\Ddots&\Ddots^{n-1\text{ times}}\\
                  &&&&&\alpha&-\alpha\\
                  \hline
                  &&&&&&-\beta&\beta\\
                  &&&&&&\Vdots&&\Ddots^{n\text{ times}}\\
                  &&&&&&-\beta&&&\beta\\
              \end{bNiceArray},
          \end{equation}
          where $\alpha,\beta,\gamma > 0$ are defined as follows:
          \begin{equation}\label{eq:alpha_beta_gamma}
              \alpha = \frac{L_{xy}}{2},\quad
              \beta = \frac{L_{xy}}{n},\quad
              \gamma = \frac{2\mu_{xy}}{\sqrt{n}}.
          \end{equation}
\end{itemize}
One can verify that the problem described above satisfies \Cref{ass:x,ass:y,ass:xy}.
Indeed, each function $g_i(z)$ is obviously $L_y$-smooth and $\mu_y$-strongly convex, and each function $f_i(z)$ is $L_x$-smooth and $\mu_x$-strongly convex due to the fact that $\smax^2(\mF_1)= \smax^2(\mF_2) =1$, and $\mu_x \leq \tdelta_x\leq L_x$, where the latter inequality is implied by \Cref{ass:Pi} as follows:
\begin{equation}\label{eq:tdelta_x_upper}
    \tdelta_x = \mu_x + {4\mu_{xy}^2}/{L_y} < \tfrac{1}{4}L_x + \tfrac{1}{4}L_x = \tfrac{1}{2}L_x.
\end{equation}
Moreover, we establish the following \Cref{lem:E}, which describes the spectral properties of matrix $\mB$. The proof is available in \Cref{sec:proof:E}.
\begin{lemma}\label{lem:E}
    Let  $n \in \{2,3,\ldots\}$ and $\alpha,\beta>0$. Then for $\gamma > 0$, the singular values of matrix $\mE$ defined in \cref{eq:E} can be bounded as follows:
    \begin{equation}
        \min\left\{\frac{n\gamma^2}{4},\frac{\beta^2}{36}, \frac{\alpha^2}{9n^2}\right\}
        \leq \smin^2(\mE) \leq \smax^2(\mE) \leq
        \max\left\{2n\gamma^2,2(n+1)\beta^2,4\alpha^2\right\},
    \end{equation}
    and for $\gamma = 0$, the singular values of matrix $\mE'$ defined in \cref{eq:E'} can be bounded as follows:
    \begin{equation}
        \min\left\{\frac{\beta^2}{36}, \frac{\alpha^2}{9n^2}\right\}
        \leq (\sminp(\mE'))^2 \leq \smax^2(\mE') \leq
        \max\left\{2(n+1)\beta^2,4\alpha^2\right\}.
    \end{equation}
\end{lemma}
Using \Cref{lem:E} and the definition of $\alpha,\beta,\gamma$ in \cref{eq:alpha_beta_gamma}, we can show that matrix $\mB$ satisfies \Cref{ass:xy} as long as $n$ is defined as follows:
\begin{equation}\label{eq:n}
    n = \floor*{\frac{L_{xy}}{6\mu_{yx}}}.
\end{equation}
Indeed, the definition of $n$ in \cref{eq:n}, the definitions of $\alpha,\beta,\gamma$ in \cref{eq:alpha_beta_gamma}, and \Cref{lem:E} imply
\begin{equation}
    \mu_{xy}^2=\mu_{yx}^2 \leq \smin^2(\mE) \leq \smax^2(\mE) \leq L_{xy}^2
\end{equation}
in the case $\mu_{xy} > 0$, and
\begin{equation}
    \mu_{yx}^2 \leq (\sminp(\mE'))^2 \leq \smax^2(\mE') \leq L_{xy}^2
\end{equation}
in the case $\mu_{xy} = 0$. Moreover, it is not hard to verify that $\range \mB^\top = \sX$ in the case $\mu_{xy} > 0$ and $\range \mB = \sY$ in both cases. Also note that $n \geq 3$ due to \Cref{ass:Pi}.


Next, we establish the following \Cref{lem:hard_sol}, which describes the solution to the problem defined above, the proof is available in \Cref{sec:proof:hard_sol}.
\begin{lemma}\label{lem:hard_sol}
    For all $d \in \{1,2,\ldots\}$, the instance of problem~\eqref{eq:main} defined above has a unique solution $(x^*,y^*) \in \sX \times \sY$. Moreover, there exists a vector $(x^\circ, y^\circ) \in  (\cL_{d,\rho})^{n_x} \times (\cL_{d,\rho})^{n_y}$ such that the following inequality holds:
    \begin{equation}
        \distsol(x^\circ,y^\circ) \leq C_\pi A^2 \rho^{2d},
    \end{equation}
    where $C_\pi > 0$ is some constant that possibly depends on the parameters $\pi \in \Pi$, but does not depend on $d$,
    $\cL_{d,\rho} \subset \R^{d}$ is a linear space which is defined for $\rho \in (0,1)$ as follows:
    \begin{equation}\label{eq:space_geometric}
        \cL_{d,\rho}
        = \range \begin{bNiceMatrix}[l,margin]
            1 & 0 & 0      & 0      & \Cdots                   \\
            1 & 0 & \rho^2 & 0      & \rho^4 & \Cdots          \\
            0 & 1 & 0      & \rho^2 & 0      & \rho^4 & \Cdots
        \end{bNiceMatrix}^\top,
    \end{equation}
    and $\rho \in (0,1)$ satisfies the following inequality:
    \begin{equation}
        \rho \geq
        \max\left\{1 - 89\cdot \frac{n}{\sqrt{\condxya}},
        \frac{1}{346}\right\}.
    \end{equation}
\end{lemma}

Let $(x^0,y^0) = (x_1^0,\ldots,x_{n_x}^0,y_1^0,\ldots,y_{n_y}^0) \in \sX \times \sY$ be defined as follows:
\begin{equation}\label{eq:x_init}
    (x^0,y^0) = \proj_{\spanset(\{\basis{d}{1}\})^{n_x+n_y}}((x^*,y^*)) = (\mI_{n_x + n_y} \otimes \mP)(x^*,y^*),
\end{equation}
where $\mP \in \R^{d \times d}$ is the orthogonal projection matrix onto the linear space $\spanset(\{\basis{d}{1}\}) \subset \R^d$, which is given as follows:
\begin{equation}
    \mP = \basis{d}{1}(\basis{d}{1})^\top.
\end{equation}
Then vector $(x^\circ, y^\circ)$ from \Cref{lem:hard_sol} satisfies the following relation:
\begin{equation}\label{eq:x_circ}
    \begin{aligned}
         & (\mI_{n_x+n_y}\otimes(\mI_d - \mP))(x^\circ,y^\circ)
        = (u_x,u_y)\otimes(0,1,0,\rho^2,\ldots) + (v_x,v_y)\otimes(0,0,1,0,\rho^2,\ldots),
        \\
         & \quad\begin{aligned}
                    \text{where}\quad
                    u_x
                     & = (u_{x,1},\ldots,u_{x,n_x}) \in \R^{n_x},
                    \quad
                    v_x = (v_{x,1},\ldots,v_{x,n_x}) \in \R^{n_x},
                    \\
                    u_y
                     & = (u_{y,1},\ldots,u_{y,n_y}) \in \R^{n_y},
                    \quad
                    v_y = (v_{y,1},\ldots,v_{y,n_y}) \in \R^{n_y}.
                \end{aligned}
    \end{aligned}
\end{equation}
Hence, we can obtain the following relation:
\begin{align*}
    \distsol(x^0,y^0)
     & =
    \delta_x\sqn{x^* - x^0}
    +\delta_y\sqn{y^* - y^0}
    \\&\aeq{uses the definition of $(x^0,y^0)$ in \cref{eq:x_init}}
    \delta_x\sqn{(\mI_{n_x}\otimes(\mI_d - \mP))x^*}
    +\delta_y\sqn{(\mI_{n_y}\otimes(\mI_d - \mP))y^*}
    \\&=
    \delta_x\sqn{(\mI_{n_x}\otimes(\mI_d - \mP))(x^* - x^\circ + x^\circ)}
    +\delta_y\sqn{(\mI_{n_y}\otimes(\mI_d - \mP))(y^* - y^\circ + y^\circ)}
    \\&\aleq{uses Young's inequality}
    2\delta_x\sqn{(\mI_{n_x}\otimes(\mI_d - \mP))(x^* - x^\circ)}
    +2\delta_y\sqn{(\mI_{n_y}\otimes(\mI_d - \mP))(y^* - y^\circ)}
    \\&
    +2\delta_x\sqn{(\mI_{n_x}\otimes(\mI_d - \mP))x^\circ}
    +2\delta_y\sqn{(\mI_{n_y}\otimes(\mI_d - \mP))y^\circ}
    \\&\leq
    2\delta_x\sqn{(\mI_{n_x}\otimes(\mI_d - \mP))x^\circ}
    +2\delta_y\sqn{(\mI_{n_y}\otimes(\mI_d - \mP))y^\circ}
    +2\distsol(x^\circ,y^\circ)
    \\&\aeq{uses \cref{eq:x_circ} and the properties of the Kronecker product}
    2(\delta_x\sqn{u_x}+\delta_y\sqn{u_y})\sqn{(0,1,0,\rho^2,\ldots)}
    \\&
    +2(\delta_x\sqn{v_x}+\delta_y\sqn{v_y})\sqn{(0,0,1,0,\rho^2,\ldots)}
    +2\distsol(x^\circ,y^\circ)
    \\&=
    2(\delta_x\sqn{u_x}+\delta_y\sqn{u_y})
    \sum_{j=0}^{\floor{d/2}-1}\rho^{4j}
    \\&
    +2(\delta_x\sqn{v_x}+\delta_y\sqn{v_y})
    \sum_{j=0}^{\floor{(d-1)/2}-1}\rho^{4j}
    +2\distsol(x^\circ,y^\circ)
    \\&=
    \frac{2(\delta_x\sqn{u_x}+\delta_y\sqn{u_y})(1 - \rho^{4\floor{d/2}})}{1-\rho^4}
    \\&
    +\frac{2(\delta_x\sqn{v_x}+\delta_y\sqn{v_y})(1 - \rho^{4\floor{(d-1)/2}})}{1-\rho^4}
    +2\distsol(x^\circ,y^\circ)
    \\&\leq
    \frac{2(\delta_x\sqn{(u_x,v_x)}+\delta_y\sqn{(u_y,v_y)})}{1-\rho^4}
    +2\distsol(x^\circ,y^\circ),
\end{align*}
where \annotate.

Further, we fix $k \in \rng{1}{d}$.
using the sparse structure of the matrices $\mF_1,\mF_2$ and $\mE$ defined in \cref{eq:F,eq:E}, respectively, and using the standard arguments \citep{nesterov2013introductory,ibrahim2020linear,zhang2022near,scaman2017optimal,scaman2018optimal,kovalev2024lower}, we can show that the output vectors $x_o(\tau) = (x_{o,1}(\tau),\ldots,x_{o,3n}(\tau)) \in \sX$ and $y_o(\tau) = (y_{o,1}(\tau),\ldots,y_{o,3n}(\tau)) \in \sY$ satisfy the following implication:
\begin{equation}\label{eq:span}
    \tau \leq D \cdot \tau_{\mB}n(k-1)
    \quad
    \Rightarrow
    \quad
    x_{o,i}(\tau),y_{o,j}(\tau) \in
    \spanset(\{\basis{d}{1},\ldots,\basis{d}{k}\})
\end{equation}
for all $i \in \rng{1}{n_x}$ and $j\in\rng{1}{n_y}$, where $D > 0$ is a universal constant.
The right-hand side of this implication implies the following:
\begin{align*}
    \distsol(x_o(\tau),y_o(\tau))
     & \ageq{uses Young's inequality}
    \tfrac{1}{2}\delta_x \sqn{x_o(\tau) - x^\circ}
    +\tfrac{1}{2}\delta_y \sqn{y_o(\tau) - y^\circ}
    -\distsol(x^\circ,y^\circ)
    \\&\ageq{uses \cref{eq:span} and the expression for $(x^\circ,y^\circ)$ in \cref{eq:x_circ}}
    \frac{(\delta_x\sqn{u_x}+\delta_y\sqn{u_y})}{2}
    \sum_{j=\floor{k/2}}^{\floor{d/2}-1}\rho^{4j}
    \\&
    +\frac{(\delta_x\sqn{v_x}+\delta_y\sqn{v_y})}{2}
    \sum_{j=\floor{(k-1)/2}}^{\floor{(d-1)/2}-1}\rho^{4j}
    -\distsol(x^\circ,y^\circ)
    \\&=
    \frac{(\delta_x\sqn{u_x}+\delta_y\sqn{u_y})(\rho^{4\floor{k/2}} - \rho^{4\floor{d/2}})}{2(1-\rho^4)}
    \\&
    +\frac{(\delta_x\sqn{v_x}+\delta_y\sqn{v_y})(\rho^{4\floor{(k-1)/2}} - \rho^{4\floor{(d-1)/2}})}{2(1-\rho^4)}
    -\distsol(x^\circ,y^\circ)
    \\&\geq
    \frac{(\delta_x\sqn{(u_x,v_x)}+\delta_y\sqn{(u_y,v_y)})(\rho^{2k} -  \rho^{2d-4})}{2(1-\rho^4)}
    -\distsol(x^\circ,y^\circ)
    \\&\ageq{uses the previously obtained upper bound on $\distsol(x^0,y^0)$}
    \frac{(\rho^{2k} -  \rho^{2d-4})}{4}\left(\distsol(x^0,y^0) - 2\distsol(x^\circ,y^\circ)\right)
    -\distsol(x^\circ,y^\circ)
    \\&=
    \frac{(\rho^{2k} -  \rho^{2d-4})}{4}\distsol(x^0,y^0)
    -\left(1 +  \frac{(\rho^{2k} -  \rho^{2d-4})}{2}\right)\distsol(x^\circ,y^\circ)
    \\&\ageq{uses \Cref{lem:hard_sol}}
    \frac{(\rho^{2k} -  \rho^{2d-4})}{4}\distsol(x^0,y^0)
    -\left(1 +  \frac{(\rho^{2k} -  \rho^{2d-4})}{2}\right)C_\pi A^2\rho^{2d},
\end{align*}
where \annotate.

Next, we establish the following \Cref{lem:init_dist}, the proof is available in \Cref{sec:proof:init_dist}.
\begin{lemma}\label{lem:init_dist}
    For all $d \in \{2,3,\ldots\}$, the unique solution to the instance of problem~\eqref{eq:main} defined above satisfies the following relation:
    \begin{equation}
        \distsol(x^0,y^0) = B_{\pi,d}A^2,
    \end{equation}
    where $B_{\pi,d} > 0$ is a constant that possibly depends on $d \in \{2,3,\ldots\}$ and the parameters $\pi \in \Pi$, Moreover there exists $\hat{d} \in \{2,3,\ldots\}$ such that the following inequality holds:
    \begin{equation}
        \min_{d\in \{\hat{d},\hat{d}+1,\ldots\}} B_{\pi,d} > 0
    \end{equation}
\end{lemma}
Using \Cref{lem:init_dist}, for $d \geq \hat{d}$, we can further lower-bound $\distsol(x_o(\tau),y_o(\tau))$ as follows:
\begin{equation}
    \distsol(x_o(\tau),y_o(\tau))
    \geq
    \frac{(\rho^{2k} -  \rho^{2d-4})}{4}B_{\pi,d}A^2
    -\left(1 +  \frac{(\rho^{2k} -  \rho^{2d-4})}{2}\right)C_\pi A^2\rho^{2d}.
\end{equation}
Next, we can choose $A = R / \sqrt{B_{\pi,d}}$ to ensure $\distsol(x^0,y^0) = R^2$ and obtain the following:
\begin{align*}
    \distsol(x_o(\tau),y_o(\tau))
     & \geq
    \frac{(\rho^{2k} -  \rho^{2d-4})}{4}R^2
    -\left(1 +  \frac{(\rho^{2k} -  \rho^{2d-4})}{2}\right)\frac{C_\pi R^2\rho^{2d}}{B_{\pi,d}}
    \\&\ageq{uses \Cref{lem:init_dist}}
    \frac{(\rho^{2k} -  \rho^{2d-4})}{4}R^2
    -\left(1 +  \frac{(\rho^{2k} -  \rho^{2d-4})}{2}\right)\frac{C_\pi R^2\rho^{2d}}{\min_{d\in \{\hat{d},\hat{d}+1,\ldots\}} B_{\pi,d}}
    \\&\ageq{is implied by choosing a large enough value of $d$}
    \frac{1}{5}\rho^{2k}R^2,
\end{align*}
where \annotate. The rest of the proof uses the lower bound on $\rho$ in \Cref{lem:hard_sol} and is almost identical to the final steps of the proof of \Cref{lem:lb_B2} in \Cref{sec:proof:lb_B2}.\qed

\subsection{Proofs of Auxiliary Lemmas}

\subsubsection{Proof of Lemma~\ref{lem:hard_sol}}\label{sec:proof:hard_sol}

In this proof, we consider the case $\mu_{xy} > 0$, since the case $\mu_{xy} = 0$ is almost identical.
Using the first-order optimality conditions~\eqref{eq:opt} in problem~\eqref{eq:main}, we obtain the following expression for the optimal dual variable $y^* \in \sY$:
\begin{equation}\label{eq:y}
    y^* = \nabla g^*(\mB x^*)
    \Rightarrow
    y^* = \left(\left(\begin{bNiceMatrix}
                {1}/{L_y}                             \\
                 & {1}/{\mu_y}                        \\
                 &             & \Ddots               \\
                 &             &        & {1}/{\mu_y} \\
            \end{bNiceMatrix} \mE\right)\otimes \mI_d\right) x^*,
\end{equation}
where the optimal primal variable is the solution to the primal minimization problem in \cref{eq:pd}:
\begin{equation}
    x^* = \argmin_{x\in\sX} f(x) + g^*(\mB x).
\end{equation}
Moreover, using the definition of functions $f(x)$ and $g(y)$ in \cref{eq:f,eq:g} and the definition of matrix $\mB$ in \cref{eq:B}, we can rewrite this problem as follows:
\begin{equation}
    \begin{aligned}
        \min_{x \in \sX}\;
         &
        \sum_{i=1}^n \left(f_1(x_i)+\frac{\beta^2}{2\mu_y}\sqn{x_i - x_{n+1}}\right)
        +\sum_{i=2n+1}^{3n}\left(f_{3n}(x_i)+\frac{\beta^2}{2\mu_y}\sqn{x_i - x_{2n}}\right)
        \\&
        +\frac{\gamma^2}{2L_y}\sqn{{\textstyle\frac{1}{n}\sum_{i=1}^n}x_i}
        +\sum_{i=n+1}^{2n} \frac{\tdelta_x}{2}\sqn{x_i}
        +\sum_{i=n+1}^{2n-1}\frac{\alpha^2}{2\mu_y}\sqn{x_{i+1} - x_i}.
    \end{aligned}
\end{equation}
It is also not hard to verify that the following inequality holds:
\begin{equation}
    \sum_{i=1}^n \left(f_1(x_i)+\frac{\beta^2}{2\mu_y}\sqn{x_i - x_{n+1}}\right)
    \geq
    n f_1({\textstyle\frac{1}{n}\sum_{i=1}^n}x_i) + \frac{n\beta^2}{2\mu_y}\sqn{{\textstyle\frac{1}{n}\sum_{i=1}^n}x_i - x_{n+1}},
\end{equation}
where equality is attained if and only if $x_1=\cdots=x_n$. Consequently, the problem can be further reformulated as follows:
\begin{equation}
    \min_{x \in \sX}\;
    n\moreau{h_1}{\frac{\mu_y}{\beta^2}}(x_{n+1}) + n\moreau{h_2}{\frac{\mu_y}{\beta^2}}(x_{2n})
    +\sum_{i=n+1}^{2n} \frac{\tdelta_x}{2}\sqn{x_i}
    +\sum_{i=n+1}^{2n-1}\frac{\alpha^2}{2\mu_y}\sqn{x_{i+1} - x_i},
\end{equation}
where functions $h_1(z),h_2(z)\colon \R^d \to \R$ are defined as follows:
\begin{equation}\label{eq:h}
    \begin{aligned}
        h_1(z) & = \frac{\tdelta_x}{2}\sqn{z} + \frac{L_x - \tdelta_x}{2}\sqn{\mF_1 z},                      \\
        h_2(z) & = \frac{\tdelta_x}{2}\sqn{z} + \frac{L_x - \tdelta_x}{2}\sqn{\mF_2 z} - A\<\basis{d}{1},z>,
    \end{aligned}
\end{equation}
and $\moreau{h_1}{\frac{\mu_y}{\beta^2}}(z)$ and $\moreau{h_2}{\frac{\mu_y}{\beta^2}}(z)$ are the corresponding Moreau envelopes. Moreover, the solution $x^*$ satisfies the following relations:
\begin{equation}\label{eq:x1}
    \begin{aligned}
        x_1^* = \cdots = x_n^*         & = \prox_{\frac{\mu_y}{\beta^2}h_1}(x_{n+1}^*) \\
        x_{2n+1}^* = \cdots = x_{3n}^* & = \prox_{\frac{\mu_y}{\beta^2}h_2}(x_{2n}^*)  \\
    \end{aligned}
\end{equation}
Further, we perform the minimization in the variables $x_{n+2},\ldots,x_{2n-1}$. Using the first-order optimality conditions, we obtain the following relations:
\begin{equation}
    \left(2 + \frac{\tdelta_x\mu_y}{\alpha^2}\right)x_i^* = x_{i-1}^* + x_{i+1}^*
    \quad\text{for}\quad
    i \in \rng{n+2}{2n-1}.
\end{equation}
This is nothing else but a linear recurrence, which is not hard to solve. Let $q > 0$ be the smallest root of the following characteristic polynomial:
\begin{equation}
    \label{eq:char}
    \left(2 + \frac{\tdelta_x\mu_y}{\alpha^2}\right)q = 1 + q^2,
\end{equation}
which is given as follows:
\begin{equation}\label{eq:q}
    q = \frac{\sqrt{4\alpha^2 + \tdelta_x\smash{\mu_y}}-\sqrt{\tdelta_x\smash{\mu_y}}}{\sqrt{4\alpha^2 + \tdelta_x\smash{\mu_y}}+\sqrt{\tdelta_x\smash{\mu_y}}}.
\end{equation}
Then $x_{n+2}^*,\ldots,x_{2n-1}^*$ can be expressed as follows:
\begin{equation}\label{eq:x2}
    x_{n+i}^* = \frac{x_{n+1}^*(q^{n-i} - q^{i-n}) + x_{2n}^*(q^{i-1} - q^{1-i})}{(q^{n-1} - q^{1-n})}.
\end{equation}
Moreover, one can observe the following:
\begin{align*}
    \mind{1em}
    \sum_{i=n+1}^{2n}\frac{\tdelta_x}{2}\sqn{x_i^*}
    +\sum_{i=n + 1}^{2n-1}\frac{\alpha^2}{2\mu_y}\sqn{x_{i+1}^* - x_{i}^*}
    \\&=
    \sum_{i=n+1}^{2n}\frac{\tdelta_x}{2}\sqn{x_i^*}
    +\sum_{i=n+1}^{2n-1}\frac{\alpha^2}{2\mu_y}\left(\sqn{x_i^*} + \sqn{x_{i+1}^*} - 2\<x_{i+1}^*,x_i^*>\right)
    \\&=
    \left(\frac{\tdelta_x}{2} + \frac{\alpha^2}{2\mu_y}\right)\left(\sqn{x_{n+1}^*}+\sqn{x_{2n}^*}\right)
    +\sum_{i=n+2}^{2n-1}\left(\frac{\tdelta_x}{2} + \frac{\alpha^2}{\mu_y}\right)\sqn{x_i^*}
    -\sum_{i=n+1}^{2n-1}\frac{\alpha^2}{\mu_y}\<x_{i+1}^*,x_i^*>
    \\&=
    \left(\frac{\tdelta_x}{2} + \frac{\alpha^2}{2\mu_y}\right)\left(\sqn{x_{n+1}^*}+\sqn{x_{2n}^*}\right)
    -\frac{\alpha^2}{2\mu_y}\left(\<x_{n+1}^*,x_{n+2}^*> + \<x_{2n}^*,x_{2n-1}^*>\right)
    \\&
    +\sum_{i=n+2}^{2n-1}\frac{\alpha^2}{2\mu_y}\left(\left(\frac{\tdelta_x\mu_y}{\alpha^2} + 2\right)\sqn{x_i^*} - \<x_i^*,x_{i+1}^* + x_{i-1}^*>\right)
    \\&\aeq{uses \cref{eq:char,eq:x2}}
    \left(\frac{\tdelta_x}{2} + \frac{\alpha^2}{2\mu_y}\right)\left(\sqn{x_{n+1}^*}+\sqn{x_{2n}^*}\right)
    -\frac{\alpha^2}{2\mu_y}\left(\<x_{n+1}^*,x_{n+2}^*> + \<x_{2n}^*,x_{2n-1}^*>\right)
    \\&\aeq{uses \cref{eq:x2}}
    \left(\frac{\tdelta_x}{2} + \frac{\alpha^2}{2\mu_y}\right)\left(\sqn{x_{n+1}^*}+\sqn{x_{2n}^*}\right)
    \\&
    -\frac{\alpha^2}{2\mu_y}\cdot\frac{(\sqn{x_{n+1}^*} + \sqn{x_{2n}^*})(q^{n-2} - q^{2-n}) + 2\<x_{n+1}^*,x_{2n}^*>(q - q^{-1})}{(q^{n-1} - q^{1-n})}
    \\&=
    \left(\frac{\tdelta_x}{2} + \frac{\alpha^2}{2\mu_y}\left(1 - \frac{(q^{n-2} - q^{2-n} + q - q^{-1})}{(q^{n-1} - q^{1-n})}\right)\right)\left(\sqn{x_{n+1}^*}+\sqn{x_{2n}^*}\right)
    \\&
    +\frac{\alpha^2}{2\mu_y}\cdot\frac{(q - q^{-1})}{(q^{n-1} - q^{1-n})}\sqn{x_{n+1}^*-x_{2n}^*}
    \\&=
    \left(\frac{\tdelta_x}{2} + \frac{\alpha^2}{2\mu_y} \cdot \frac{(1-q)(1-q^{n-2})}{(1+q^{n-1})}\right)
    \left(\sqn{x_{n+1}^*} + \sqn{x_{2n}^*}\right)
    \\&
    +\frac{\alpha^2}{2\mu_y}\cdot\frac{(1-q)(1+q)}{q(q^{1-n} - q^{n-1})}\sqn{x_{n+1}^*-x_{2n}^*}
    \\&=
    \left(\frac{\tdelta_x}{2} + \frac{\omega_n\alpha^2}{2\mu_y}\right)
    \left(\sqn{x_{n+1}^*} + \sqn{x_{2n}^*}\right)
    +\frac{\nu_n\alpha^2}{2\mu_y}\sqn{x_{n+1}^*-x_{2n}^*}
\end{align*}
where \annotate, and $\omega_n,\nu_n > 0$ are defined as follows:
\begin{equation}\label{eq:omega_nu}
    \omega_n = \frac{(1-q)(1-q^{n-2})}{(1+q^{n-1})},\quad
    \nu_n = \frac{(1-q)(1+q)}{q(q^{1-n} - q^{n-1})}.
\end{equation}
In addition, we can observe that the following relation holds:
\begin{equation}
    \frac{1}{2}\sqn{x_{n+1}^*-x_{2n}^*} = \min_{z \in \R^d} \left(\sqn{z - x_{n+1}^*} + \sqn{z - x_{2n}^*}\right).
\end{equation}
Therefore, the problem can be further reformulated as follows:
\begin{equation}
    \begin{aligned}
        \min_{x \in \sX}\min_{z \in \R^d}\;
         &
        \moreau{h_1}{\frac{\mu_y}{\beta^2}}(x_{n+1})
        +\left(\frac{\tdelta_x}{2n} + \frac{\omega_n\alpha^2}{2n\mu_y}\right)\sqn{x_{n+1}}
        +\frac{\nu_n\alpha^2}{n\mu_y}\sqn{x_{n+1}-z}
        \\&
        +\moreau{h_2}{\frac{\mu_y}{\beta^2}}(x_{2n})
        +\left(\frac{\tdelta_x}{2n} + \frac{\omega_n\alpha^2}{2n\mu_y}\right)\sqn{x_{2n}}
        +\frac{\nu_n\alpha^2}{n\mu_y}\sqn{x_{2n}-z}.
    \end{aligned}
\end{equation}
Let functions $h_1^{+}(z), h_2^{+}(z)\colon \R^d\to\R$ be defined as follows:
\begin{equation}
    h_j^{+}(z) = \moreau{h_j}{\frac{\mu_y}{\beta^2}}(z)
    +\left(\frac{\tdelta_x}{2n} + \frac{\omega_n\alpha^2}{2n\mu_y}\right)\sqn{z}
    ,\quad
    j = 1,2,
\end{equation}
and let functions $h_1^{++}(z), h_2^{++}(z)\colon \R^d\to\R$ be defined as follows:
\begin{equation}\label{eq:h++2}
    h_j^{++}(z) = \moreau{h_j^{+}}{\frac{n\mu_y}{2\nu_n\alpha^2}}(z)
    ,\quad
    j = 1,2.
\end{equation}
Then the latter problem reformulation can be rewritten as follows:
\begin{equation}\label{eq:z}
    z^* = \argmin_{z \in \R^d} h_1^{++}(z) + h_2^{++}(z),
\end{equation}
and the solution $x^* \in \sX$ satisfies the following relation:
\begin{equation}\label{eq:x3}
    x_{n+1}^* = \prox_{\frac{n\mu_y}{2\nu_n\alpha^2} h_1^+}(z^*),\quad
    x_{2n}^* = \prox_{\frac{n\mu_y}{2\nu_n\alpha^2} h_2^+}(z^*).
\end{equation}
Next, we establish the following \Cref{lem:Moreau}, which is used to obtain the explicit expressions for the Moreau envelopes. The proof is available in \Cref{sec:proof:Moreau}.
\begin{lemma}\label{lem:Moreau}
    Let function $h(z)\colon \R^d \to \R$ be defined as follows:
    \begin{equation}
        h(z) = \frac{\mu}{2}\sqn{z} + \frac{L-\mu}{2}\sqn{\mF z} - \<b,z>,
    \end{equation}
    where $L > \mu > 0$, and $\mF \in \R^{p\times d}$ and $b \in \R^d$ satisfy the following assumptions:
    \begin{equation}
        \mF\mF^\top = \mI_p,\quad
        \mF^\top\mF b = 0.
    \end{equation}
    Then for $\lambda > 0$, the Moreau envelope $\moreau{h}{\lambda}(z)$ is given as follows:
    \begin{equation}
        \moreau{h}{\lambda}(z) = \frac{\mu_\lambda}{2}\sqn{z} + \frac{L_\lambda-\mu_\lambda}{2}\sqn{\mF z} - B_\lambda\<b,z> - C_\lambda,
    \end{equation}
    where constants $L_\lambda > \mu_\lambda > 0$ and $B_\lambda, C_\lambda \in \R$ are defined as follows:
    \begin{equation}
        L_\lambda = \left(\lambda + 1/L\right)^{-1},\quad
        \mu_\lambda = \left(\lambda + 1/\mu\right)^{-1},\quad
        B_\lambda = (1 + \lambda\mu)^{-1},\quad
        C_{\lambda} = \frac{\lambda\sqn{b}}{2(1+\lambda\mu)}.
    \end{equation}
\end{lemma}
Using \Cref{lem:Moreau} and the definition of functions $h_1(z),h_2(z)$ in \cref{eq:h}, we can express functions $h_1^{+}(z), h_2^{+}(z)$ as follows:
\begin{equation}\label{eq:h+}
    \begin{aligned}
        h_1^+(z) & = \frac{\mu^+}{2}\sqn{z} + \frac{L^+ - \mu^+}{2}\sqn{\mF_1 z} + \const,
        \\
        h_2^+(z) & = \frac{\mu^+}{2}\sqn{z} + \frac{L^+ - \mu^+}{2}\sqn{\mF_2 z} - A^+\<\basis{d}{1},z> + \const,
    \end{aligned}
\end{equation}
and functions $h_1^{++}(z), h_2^{++}(z)$ can be expressed as follows:
\begin{equation}\label{eq:h++}
    \begin{aligned}
        h_1^{++}(z) & = \frac{\mu^{++}}{2}\sqn{z} + \frac{L^{++} - \mu^{++}}{2}\sqn{\mF_1 z} + \const,
        \\
        h_2^{++}(z) & = \frac{\mu^{++}}{2}\sqn{z} + \frac{L^{++} - \mu^{++}}{2}\sqn{\mF_2 z} - A^{++}\<\basis{d}{1},z> + \const,
    \end{aligned}
\end{equation}
where constants $L^+ > \mu^+ > 0$ are defined as follows:
\begin{equation}
    \label{eq:+}
    L^+       = \left(\frac{1}{L_x} + \frac{\mu_y}{\beta^2}\right)^{-1} + \frac{\tdelta_x}{n} + \frac{\omega_n\alpha^2}{n\mu_y},
    \quad
    \mu^+     = \left(\frac{1}{\tdelta_x} + \frac{\mu_y}{\beta^2}\right)^{-1} + \frac{\tdelta_x}{n} + \frac{\omega_n\alpha^2}{n\mu_y},
\end{equation}
constants $L^{++} > \mu^{++} > 0$ are defined as follows:
\begin{equation}
    \label{eq:++}
    L^{++}    = \left(\frac{1}{L^+} + \frac{n\mu_y}{2\nu_n\alpha^2}\right)^{-1},
    \quad
    \mu^{++}  = \left(\frac{1}{\mu^+} + \frac{n\mu_y}{2\nu_n\alpha^2}\right)^{-1},
\end{equation}
and constants $A^+,A^{++}\in R$ are defined as follows:
\begin{equation}\label{eq:A}
    A^+ = \frac{\beta^2A}{\beta^2 + \mu_y\tdelta_x},
    \quad
    A^{++} = \frac{2\nu_n\alpha^2 A^+}{2\nu_n\alpha^2 + n\mu_y\mu^+}.
\end{equation}
Next, we establish the following \Cref{lem:geometric}, the proof is available in \Cref{sec:proof:geometric}.
\begin{lemma}\label{lem:geometric}
    Let $\hz = (\hz_1,\ldots,\hz_d) \in \R^d$ be defined as follows:
    \begin{equation}
        \hz = \argmin_{z\in \R^d}\mu\sqn{z} + \frac{L - \mu}{2}\sqn{\mF z} - B\<\basis{d}{1},z>,
    \end{equation}
    where $L > \mu > 0$, $B \in \R$ and matrix $\mF \in \R^{(d-1)\times d}$ is defined as follows:
    \begin{equation}\label{eq:F2}
        \mF = \frac{1}{\sqrt{2}}\begin{bNiceMatrix}
            1 & -1     &        &    \\
              & \Ddots & \Ddots      \\
              &        & 1      & -1
        \end{bNiceMatrix}.
    \end{equation}
    Then there exists $z^\circ \in \spanset(\{(1,\rho,\ldots,\rho^{d-1})\})$ wuch that the following inequality holds:
    \begin{equation}
        \norm{z^\circ - \hz}
        \leq
        \frac{B\rho^d}{2\mu},
    \end{equation}
    where $\rho \in (0,1)$ is defined as follows:
    \begin{equation}
        \rho = \frac{\sqrt{L} - \sqrt{\mu}}{\sqrt{L} + \sqrt{\mu}}.
    \end{equation}
\end{lemma}
From \Cref{lem:geometric}, the definition of $z^*$ in \cref{eq:z}, and the definition of functions $h_j(z)$, $h_j^{+}(z)$, and $h_j^{++}(z)$ in \cref{eq:h,eq:h+,eq:h++}, it follows that there exists $z^\circ \in \spanset(\{(1,\rho,\ldots,\rho^{d-1})\})$ such that the following inequality holds:
\begin{equation}\label{eq:z_ub}
    \norm{z^\circ - z^*} \leq \frac{A^{++}\rho^d}{2\mu^{++}},
\end{equation}
where $\rho \in (0,1)$ is defined as follows:
\begin{equation}\label{eq:rho}
    \rho = \frac{\sqrt{L^{++}} - \sqrt{\mu^{++}}}{\sqrt{L^{++}} + \sqrt{\mu^{++}}}.
\end{equation}

Finally, we obtain the desired statement of \Cref{lem:hard_sol} with the help of the following \Cref{lem:geometric2}. The proof is available in \Cref{sec:proof:geometric2}.
\begin{lemma}\label{lem:geometric2}
    Let vector $(x^\circ,y^\circ) = (x^\circ_1,\ldots,x^\circ_{n_x},y^\circ_1,\ldots,y^\circ_{n_y}) \in \sX \times \sY$ be defined as follows:
    \begin{equation}
        x^\circ_i = \proj_{\cL_{d,\rho}}(x'_i),\quad
        y^\circ_i = \proj_{\cL_{d,\rho}}(y'_i),
    \end{equation}
    where vector $(x',y') = (x'_1,\ldots,x'_{n_x},y'_1,\ldots,y'_{n_y}) \in \sX \times \sY$ is defined as follows:
    \begin{equation}
        \label{eq:x_prime}
        \begin{aligned}
            x_i' & = \begin{cases}
                         \prox_{\frac{\mu_y}{\beta^2}h_1}(x_{n+1}')
                          & i \in \rng{1}{n}      \\
                         \prox_{\frac{n\mu_y}{2\nu_n\alpha^2} h_1^+}(z^\circ)
                          & i = n+1               \\
                         \frac{x_{n+1}'(q^{2n-i} - q^{i-2n}) + x_{2n}'(q^{i-(n+1)} - q^{(n+1)-i})}{(q^{n-1} - q^{1-n})}
                          & i \in \rng{n+2}{2n-1} \\
                         \prox_{\frac{n\mu_y}{2\nu_n\alpha^2} h_2^+}(z^\circ)
                          & i = 2n                \\
                         \prox_{\frac{\mu_y}{\beta^2}h_2}(x_{2n}')
                          & i \in \rng{2n+1}{3n}
                     \end{cases},
            \\
            y'   & = \nabla g^*(\mB x').
        \end{aligned}
    \end{equation}
    Then the following inequality holds:
    \begin{equation}
        \distsol(x^\circ,y^\circ)
        \leq
        C_\pi A^2\rho^{2d}.
    \end{equation}
    where $C_\pi > 0$ is some constant that possibly depends on the parameters $\pi \in \Pi$
\end{lemma}

It remains to lower-bound $\rho$. It is done with the help of the following \Cref{lem:q^n,,lem:omega_nu,,lem:++}, the proofs are available in \Cref{sec:proof:q^n,,sec:proof:omega_nu,,sec:proof:++}, respectively.
\begin{lemma}\label{lem:q^n}
    Under assumption $\mu_x\mu_y \leq \mu_{xy}^2$, the following ineqality holds:
    \begin{equation}
        q^{-n} \leq 2.
    \end{equation}
\end{lemma}
\begin{lemma}\label{lem:omega_nu}
    The following inequalities hold:
    \begin{equation}
        \begin{aligned}
            \omega_n & \leq \frac{(n-2)\tdelta_x\smash{\mu_y}}{\alpha^2}, \\
            \nu_n    & \geq \frac{1}{4(n-1)}.
        \end{aligned}
    \end{equation}
\end{lemma}
\begin{lemma}\label{lem:++}
    Constants $L^{++}$ and $\mu^{++}$ defined in \cref{eq:++} satisfy the following inequality:
    \begin{equation}
        \frac{L^{++}}{\mu^{++}} \geq 1 + \max\left\{
        \frac{1}{86},
        \frac{1}{55}\cdot\frac{\mu_{yx}^2}{\delta_x\delta_y}
        \right\}.
    \end{equation}
\end{lemma}
Using \Cref{lem:++}, we can lower-bound $\rho$ as follows:
\begin{align*}
    \rho
     & \aeq{uses the definition of $\rho$ in \cref{eq:rho}}
    1 - \frac{2}{\sqrt{\frac{L^{++}}{\mu^{++}}} + 1}
    \\&\ageq{uses \Cref{lem:++}}
    \max\left\{1 - \sqrt{220}\cdot\sqrt{\frac{\delta_x\delta_y}{\mu_{yx}^2}},
    \frac{\sqrt{87}-\sqrt{86}}{\sqrt{87}+\sqrt{86}}\right\}
    \\&\geq
    \max\left\{1 - \sqrt{220}\cdot\sqrt{\frac{\delta_x\delta_y}{\mu_{yx}^2}},
    \frac{1}{346}\right\}
    \\&\aeq{uses the definition of $\condxya$ in \cref{eq:kappa}}
    \max\left\{1 - \frac{\sqrt{220}}{\sqrt{\condxya}}\cdot\frac{L_{xy}}{\mu_{yx}},
    \frac{1}{346}\right\}
    \\&\ageq{uses the definition of $n$ in \cref{eq:n}}
    \max\left\{1 - 89\cdot \frac{n}{\sqrt{\condxya}},
    \frac{1}{346}\right\},
\end{align*}
where \annotate, which concludes the proof.\qed
\subsubsection{Proof of Lemma~\ref{lem:E}}\label{sec:proof:E}

Let matrices $\mW_i,\mW_i'\in\R^{i\times i}$ be defined for $i \in \rng{1}{n-1}$ as follows:
\begin{equation}
    \mW_i = \begin{bNiceMatrix}
        2  & -1 &        &        &    \\
        -1 & 2  &        &        &    \\
           &    & \Ddots & \Ddots &    \\
           &    & \Ddots & 2      & -1 \\
           &    &        & -1     & 2  \\
    \end{bNiceMatrix},
    \quad
    \mW_i' = \begin{bNiceMatrix}
        1  & -1 &        &        &    \\
        -1 & 2  &        &        &    \\
           &    & \Ddots & \Ddots &    \\
           &    & \Ddots & 2      & -1 \\
           &    &        & -1     & 2  \\
    \end{bNiceMatrix}.
\end{equation}
Then using the definition of matrix $\mE$ in \cref{eq:E}, we can write the matrix $\mE\mE^\top$ as the following block matrix:
\begin{equation}
    \mE\mE^\top =
    \begin{bNiceArray}[margin,cell-space-limits=2pt]{c|ccc|ccc|ccc}
        n\gamma^2&\beta\gamma&\Cdots^{n\text{ times}}&\beta\gamma\\
        \hline
        \beta\gamma&\Block{3-3}{\beta^2(\mI_n + \mJ_n)}&&&-\alpha\beta\\
        \Vdots^{n\text{ times}}&&&&\Vdots\\
        \beta\gamma&&&&-\alpha\beta\\
        \hline
        &-\alpha\beta&\Cdots&-\alpha\beta&\Block{3-3}{\alpha^2\mW_{n-1}}&&\\
        \\
        &&&&&&&\alpha\beta&\Cdots^{n\text{ times}}&\alpha\beta\\
        \hline
        &&&&&&\alpha\beta&\Block{3-3}{\beta^2(\mI_n + \mJ_n)}\\
        &&&&&&\Vdots^{n\text{ times}}\\
        &&&&&&\alpha\beta\\
    \end{bNiceArray}.
\end{equation}
Furthermore, let matrices $\mQ_i,\mQ_i' \in \R{(n+i)\times(n+i)}$ be defined for $i \in \rng{1}{n}$ as follows:
\begin{equation}
    \mQ_i =
    \begin{bNiceArray}[margin,cell-space-limits=2pt]{ccc|ccc}
        \Block{3-3}{\alpha^2\mW_i}&&\\
        \\
        &&&\alpha\beta&\Cdots^{n\text{ times}}&\alpha\beta\\
        \hline
        &&\alpha\beta&\Block{3-3}{\beta^2(\mI_n + \mJ_n)}\\
        &&\Vdots^{n\text{ times}}\\
        &&\alpha\beta\\
    \end{bNiceArray},
    \quad
    \mQ_i' =
    \begin{bNiceArray}[margin,cell-space-limits=2pt]{ccc|ccc}
        \Block{3-3}{\alpha^2\mW_i'}&&\\
        \\
        &&&\alpha\beta&\Cdots^{n\text{ times}}&\alpha\beta\\
        \hline
        &&\alpha\beta&\Block{3-3}{\beta^2(\mI_n + \mJ_n)}\\
        &&\Vdots^{n\text{ times}}\\
        &&\alpha\beta\\
    \end{bNiceArray}.
\end{equation}
It is not hard to verify that the following matrix inequality holds:
\begin{equation}\label{eq:Q_lb2}
    \mQ_{n-1}'
    \succeq
    \beta^2
    \begin{bNiceArray}[margin,cell-space-limits=2pt]{c|c}
        \mO_{n-1}\\
        \hline
        &\mI_{n}\\
    \end{bNiceArray}.
\end{equation}
Moreover, we can show that for $i \in \rng{1}{n}$, matrices $\mQ_{i},\mQ_{i}'$ satisfy the following inequalities:
\begin{equation}\label{eq:Q_lb}
    \begin{aligned}
        \mQ_i  & \frac{\alpha^2(n+i+1)}{n+i}(\basis{n+i}{1})(\basis{n+i}{1})^\top,  \\
        \mQ_i' & \succeq \frac{\alpha^2}{n+i}(\basis{n+i}{1})(\basis{n+i}{1})^\top.
    \end{aligned}
\end{equation}
Indeed, let us prove \cref{eq:Q_lb} by induction. The base case $i=1$ is trivial:
\begin{equation}
    \mQ_1' =
    \begin{bNiceArray}[margin,cell-space-limits=2pt]{c|ccc}
        \alpha^2&\alpha\beta&\Cdots&\alpha\beta\\
        \hline
        \alpha\beta&\Block{3-3}{\beta^2(\mI_n + \mJ_n)}\\
        \Vdots\\
        \alpha\beta\\
    \end{bNiceArray}
    \succeq
    \begin{bNiceArray}[margin,cell-space-limits=2pt]{c|ccc}
        \frac{\alpha^2}{n+1}\\
        \hline
        &0\\
        &&\Ddots\\
        &&&0\\
    \end{bNiceArray}
    =
    \frac{\alpha^2}{n+1}(\basis{n+1}{1})(\basis{n+1}{1})^\top.
\end{equation}
Furthermore, for an arbitrary index $i \in \rng{2}{n-1}$, we can show that
\begin{align*}
    \mQ_i'
    =
    \begin{bNiceArray}[margin,cell-space-limits=2pt]{c|cccc}
        \alpha^2&\alpha^2&0&\Cdots&0\\
        \hline
        \alpha^2&\Block{4-4}{\mQ_{i-1}}&&\\
        0\\
        \Vdots\\
        0\\
    \end{bNiceArray}
    \asucceq{uses the induction hypothesis~\eqref{eq:Q_lb} for the index $i-1$}
    \begin{bNiceArray}[margin,cell-space-limits=2pt]{c|cccc}
        \alpha^2&\alpha^2&0&\Cdots&0\\
        \hline
        \alpha^2&\frac{\alpha^2(n+i)}{n+i-1}&&\\
        0&&0\\
        \Vdots&&&\Ddots\\
        0&&&&0\\
    \end{bNiceArray}
    \asucceq{uses Young's inequality}
    \begin{bNiceArray}[margin]{c|ccc}
        \frac{\alpha^2}{n+i}\\
        \hline
        &0\\
        &&\Ddots\\
        &&&0\\
    \end{bNiceArray}.
\end{align*}
which is nothing else bu the induction hypothesis~\eqref{eq:Q_lb} for the index $i$, and where \annotate.
Next, we can obtain the following inequality for the matrix $\mE\mE^\top$:
\begin{align*}
    \mE\mE^\top
     & =
    \begin{bNiceArray}[margin,cell-space-limits=2pt]{c|ccc|ccc}
        n\gamma^2&\beta\gamma&\Cdots&\beta\gamma\\
        \hline
        \beta\gamma&\Block{3-3}{\beta^2(\mI_n + \mJ_n)}&&&-\alpha\beta\\
        \Vdots&&&&\Vdots\\
        \beta\gamma&&&&-\alpha\beta\\
        \hline
        &-\alpha\beta&\Cdots&-\alpha\beta&\Block{3-3}{\mQ_{n-1}}&&\\
        \\
        &&&&&&\\
    \end{bNiceArray}
    \succeq
    \begin{bNiceArray}[margin,cell-space-limits=2pt]{c|ccc|ccc}
        n\gamma^2&\beta\gamma&\Cdots&\beta\gamma\\
        \hline
        \beta\gamma&\Block{3-3}{\frac{\beta^2(n+1)}{n}\mJ_n}&&&-\alpha\beta\\
        \Vdots&&&&\Vdots\\
        \beta\gamma&&&&-\alpha\beta\\
        \hline
        &-\alpha\beta&\Cdots&-\alpha\beta&\Block{3-3}{\mQ_{n-1}}&&\\
        \\
        &&&&&&\\
    \end{bNiceArray}
    \\&\succeq
    \begin{bNiceArray}[margin,cell-space-limits=2pt]{c|ccc|ccc}
        0&\\
        \hline
        &\Block{3-3}{\beta^2\mJ_n}&&&-\alpha\beta\\
        &&&&\Vdots\\
        &&&&-\alpha\beta\\
        \hline
        &-\alpha\beta&\Cdots&-\alpha\beta&\Block{3-3}{\mQ_{n-1}}&&\\
        \\
        &&&&&&\\
    \end{bNiceArray}
    \succeq
    \begin{bNiceArray}[margin,cell-space-limits=2pt]{c|c|c}
        0&\\
        \hline
        &\mO_n\\
        \hline
        &&\mQ_{n-1}'\\
    \end{bNiceArray}
    \\&\succeq
    \frac{\alpha^2}{n+i}(\basis{3n}{2n+1-i})(\basis{3n}{2n+1-i})^\top.
\end{align*}
This implies the following inequality for the matrix $\mE\mE^\top$:
\begin{align*}
    \mE\mE^\top
     & \aeq{uses the facdt that $\sum_{i=1}^{n-1}\frac{2(n+i)}{3n(n-1)} = 1$}
    \sum_{i=1}^{n-1}\frac{2(n+i)}{3n(n-1)}\mE\mE^\top
    \succeq
    \sum_{i=1}^{n-1}\frac{2(n+i)}{3n(n-1)}
    \cdot\frac{\alpha^2}{n+i}(\basis{3n}{2n+1-i})(\basis{3n}{2n+1-i})^\top
    \\&=
    \frac{2\alpha^2}{3n(n-1)}
    \begin{bNiceArray}[margin,cell-space-limits=2pt]{c|c|c}
        \mO_{n+1}\\
        \hline
        &\mI_{n-1}\\
        \hline
        &&\mO_{n}
    \end{bNiceArray}.
\end{align*}
where \annotate.
Furthermore, using \cref{eq:Q_lb2}, we can obtain the following inequality for the matrix $\mE\mE^\top$:
\begin{align*}
    \mE\mE^\top
     & \succeq
    \begin{bNiceArray}[margin,cell-space-limits=2pt]{c|c|c}
        0&\\
        \hline
        &\mO_n\\
        \hline
        &&\mQ_{n-1}'\\
    \end{bNiceArray}
    \succeq
    \beta^2
    \begin{bNiceArray}[margin,cell-space-limits=2pt]{c|c}
        \mO_{2n}\\
        \hline
        &\mI_{n}\\
    \end{bNiceArray}.
\end{align*}
Next, we can obtain the following inequality for the matrix $\mE\mE^\top$:
\begin{align*}
    \mE\mE^\top
     & =
    \begin{bNiceArray}[margin,cell-space-limits=2pt]{c|ccc|ccc}
        n\gamma^2&\beta\gamma&\Cdots&\beta\gamma\\
        \hline
        \beta\gamma&\Block{3-3}{\beta^2(\mI_n + \mJ_n)}&&&-\alpha\beta\\
        \Vdots&&&&\Vdots\\
        \beta\gamma&&&&-\alpha\beta\\
        \hline
        &-\alpha\beta&\Cdots&-\alpha\beta&\Block{3-3}{\mQ_{n-1}}&&\\
        \\
        &&&&&&\\
    \end{bNiceArray}
    \\&
    \asucceq{uses \cref{eq:Q_lb}}
    \begin{bNiceArray}[margin,cell-space-limits=2pt]{c|ccc|c|c}
        n\gamma^2&\beta\gamma&\Cdots&\beta\gamma\\
        \hline
        \beta\gamma&\Block{3-3}{\beta^2(\mI_n + \mJ_n)}&&&-\alpha\beta\\
        \Vdots&&&&\Vdots\\
        \beta\gamma&&&&-\alpha\beta\\
        \hline
        &-\alpha\beta&\Cdots&-\alpha\beta&\frac{2n\alpha^2}{2n-1}\\
        \hline
        &&&&&\mO_{2n-2}\\
    \end{bNiceArray}
    \\&
    \asucceq{uses Young's inequality}
    \begin{bNiceArray}[margin,cell-space-limits=2pt]{c|c|c|c}
        0\\
        \hline
        &\beta^2\mI_n + \beta^2\left(1 - \frac{1}{n} - \frac{2n-1}{2n}\right)\mJ_n\\
        \hline
        &&0\\
        \hline
        &&&\mO_{2n-2}\\
    \end{bNiceArray}
    \\&=
    \begin{bNiceArray}[margin,cell-space-limits=2pt]{c|c|c|c}
        0\\
        \hline
        &\beta^2\mI_n - \frac{\beta^2}{2n}\mJ_n\\
        \hline
        &&0\\
        \hline
        &&&\mO_{2n-2}\\
    \end{bNiceArray}
    \\&\asucceq{uses the fact that $\mJ_n \preceq n \mI_n$}
    \begin{bNiceArray}[margin,cell-space-limits=2pt]{c|c|c|c}
        0\\
        \hline
        &\frac{\beta^2}{2}\mI_n\\
        \hline
        &&0\\
        \hline
        &&&\mO_{2n-2}\\
    \end{bNiceArray}.
\end{align*}
where \annotate.
Furthermore, we can obtain the following inequality for the matrix $\mE\mE^\top$:
\begin{align*}
    \mE\mE^\top
     & \asucceq{uses the ineqality obtained above}
    \begin{bNiceArray}[margin,cell-space-limits=2pt]{c|ccc|c|c}
        n\gamma^2&\beta\gamma&\Cdots&\beta\gamma\\
        \hline
        \beta\gamma&\Block{3-3}{\beta(\mI_n + \mJ_n)}&&&-\alpha\beta\\
        \Vdots&&&&\Vdots\\
        \beta\gamma&&&&-\alpha\beta\\
        \hline
        &-\alpha\beta&\Cdots&-\alpha\beta&\frac{2n\alpha^2}{2n-1}\\
        \hline
        &&&&&\mO_{2n-2}\\
    \end{bNiceArray}
    \\& \asucceq{uses the fact that $\mJ_n \preceq n \mI_n$}
    \begin{bNiceArray}[margin,cell-space-limits=2pt]{c|ccc|c|c}
        n\gamma^2&\beta\gamma&\Cdots&\beta\gamma\\
        \hline
        \beta\gamma&\Block{3-3}{\frac{\beta^2(n+1)}{n}\mJ_n}&&&-\alpha\beta\\
        \Vdots&&&&\Vdots\\
        \beta\gamma&&&&-\alpha\beta\\
        \hline
        &-\alpha\beta&\Cdots&-\alpha\beta&\frac{2n\alpha^2}{2n-1}\\
        \hline
        &&&&&\mO_{2n-2}\\
    \end{bNiceArray}
    \\& \asucceq{use Young's ineqality}
    \begin{bNiceArray}[margin,cell-space-limits=2pt]{c|ccc|c|c}
        n\gamma^2&\beta\gamma&\Cdots&\beta\gamma\\
        \hline
        \beta\gamma&\Block{3-3}{\frac{3\beta^2}{2n}\mJ_n}&&&\\
        \Vdots&&&&\\
        \beta\gamma&&&&\\
        \hline
        &&&&0\\
        \hline
        &&&&&\mO_{2n-2}\\
    \end{bNiceArray}
    \\& \asucceq{use Young's ineqality}
    \frac{n\gamma^2}{3}
    \begin{bNiceArray}[margin,cell-space-limits=2pt]{c|c|c|c}
        1\\
        \hline
        &\mO_n\\
        \hline
        &&0\\
        \hline
        &&&\mO_{2n-2}\\
    \end{bNiceArray}.
\end{align*}
where \annotate.
Now, we sum all the inequalities for the matrix $\mE\mE^\top$ obtained above with positive coefficients $\theta_1,\theta_2,\theta_3,\theta_4>0$ and obtain the following:
\begin{equation}
    \mE\mE^\top \succeq
    \begin{bNiceArray}[margin,cell-space-limits=2pt]{c|c|c|c}
        \theta_1\cdot\frac{n\gamma^2}{3}\\
        \hline
        &\theta_2\cdot\frac{\beta^2}{2}\mI_n\\
        \hline
        &&\theta_3\cdot\frac{2\alpha^2}{3n^2}\mI_{n-1}\\
        \hline
        &&&\theta_4\cdot \beta^2\mI_n
    \end{bNiceArray}.
\end{equation}
Choosing
$\theta_1 = \frac{3}{4}$,
$\theta_2 = \frac{1}{18}$,
$\theta_3 = \frac{1}{6}$, and
$\theta_4 = \frac{1}{36}$ implies the following:
\begin{equation}
    \mE\mE^\top \succeq
    \min\left\{\frac{n\gamma^2}{4},\frac{\beta^2}{36}, \frac{\alpha^2}{9n^2}\right\}\mI_{3n}.
\end{equation}
Finally, we obtain the following inequality for the matrix $\mE\mE^\top$:
\begin{align*}
    \mE\mE^\top
     & \apreceq{uses the expression for the matrix $\mE\mE^\top$ and Young's inequality}
    \begin{bNiceArray}[margin,cell-space-limits=2pt]{c|c|ccccc|c}
        2n\gamma^2\\
        \hline
        &\beta^2\mI_n + \frac{\beta^2(2n+1)}{n}\mJ_n\\
        \hline
        &&3\alpha^2&-\alpha^2\\
        &&-\alpha^2&2\alpha^2&-\alpha^2\\
        &&&\Ddots&\Ddots&\Ddots\\
        &&&&-\alpha^2&2\alpha^2&-\alpha^2\\
        &&&&&-\alpha^2&3\alpha^2\\
        \hline
        &&&&&&&\beta^2(\mI_n + 2\mJ_n)
    \end{bNiceArray}.
    \\&
    \apreceq{uses Young's inequality}
    \begin{bNiceArray}[margin,cell-space-limits=2pt]{c|c|c|c}
        2n\gamma^2\\
        \hline
        &2(n+1)\beta^2\mI_n\\
        \hline
        &&4\alpha^2\mI_{n-1}\\
        \hline
        &&&(2n+1)\beta^2\mI_n\\
    \end{bNiceArray}
    \\&\preceq
    \max\left\{2n\gamma^2,2(n+1)\beta^2,4\alpha^2\right\}\mI_{3n}.
\end{align*}
where \annotate, which concludes the proof in the case $\gamma > 0$. The remaining case $\gamma = 0$ is a trivial extension of the case $\gamma > 0$. \qed
\subsubsection{Proof of Lemma~\ref{lem:Moreau}}\label{sec:proof:Moreau}

Using the definition of the Moreau envelope, we get
\begin{equation}
    \moreau{h}{\lambda} = \argmin_{z' \in \R^d}
    \frac{1}{2\lambda}\sqn{z'-z}
    +
    \frac{\mu}{2}\sqn{z'} + \frac{L-\mu}{2}\sqn{\mF z'} - \<b,z'>.
\end{equation}
Using the first-order optimality conditions, we get
\begin{equation}
    \frac{1}{\lambda}z + b = \left(\left(\frac{1}{\lambda} + \mu\right) \mI_d + (L-\mu)\mF\mF^\top\right)z'.
\end{equation}
Hence, we obtain the following:
\begin{align*}
    \moreau{h}{\lambda}
     & \aeq{uses the definition of the Moreau envelope}
    \frac{1}{2\lambda}\sqn{z'-z}
    +
    \frac{\mu}{2}\sqn{z'} + \frac{L-\mu}{2}\sqn{\mF z'} - \<b,z'>
    \\&=
    \frac{1}{2\lambda}\sqn{z}
    +\frac{1}{2}\sqn{z'}_{\left(\frac{1}{\lambda} + \mu\right)\mI_d + (L-\mu)\mF^\top\mF}
    - \<z',\tfrac{1}{\lambda}z + b>
    \\&\aeq{uses the expression for $z'$}
    \frac{1}{2\lambda}\sqn{z}
    -\frac{1}{2}\sqn{\tfrac{1}{\lambda}z + b}_{\left(\left(\frac{1}{\lambda} + \mu\right)\mI_d + (L-\mu)\mF^\top\mF\right)^{-1}}
    \\&\aeq{uses the Woodbury matrix identity}
    \frac{1}{2\lambda}\sqn{z}
    -\frac{1}{2}\sqn{\tfrac{1}{\lambda}z + b}_{\left(\frac{1}{\lambda} + \mu\right)^{-1}\mI_d - \left(\frac{1}{\lambda} + \mu\right)^{-2}\mF^\top\left(\frac{1}{L-\mu}\mI_p + \left(\frac{1}{\lambda} + \mu\right)^{-1}\mF\mF^\top\right)^{-1}\mF}
    \\&\aeq{uses the assumption $\mF\mF^\top = \mI_p$}
    \frac{1}{2\lambda}\sqn{z}
    -\frac{1}{2}\sqn{\tfrac{1}{\lambda}z + b}_{\left(\frac{1}{\lambda} + \mu\right)^{-1}\mI_d - \left(\frac{1}{\lambda} + \mu\right)^{-2}\mF^\top\left(\frac{1}{L-\mu}\mI_p + \left(\frac{1}{\lambda} + \mu\right)^{-1}\mI_p\right)^{-1}\mF}
    \\&=
    \frac{1}{2\lambda}\sqn{z}
    -\frac{1}{2}\sqn{\tfrac{1}{\lambda}z + b}_{\frac{\lambda}{1+\lambda\mu}\mI_d - \frac{\lambda^2(L-\mu)}{(1+\lambda L)(1+\lambda\mu)}\mF^\top\mF}
    \\&\aeq{uses the definitions of $L_\lambda$ and $\mu_\lambda$}
    \frac{1}{2\lambda}\sqn{z}
    -\frac{1}{2}\sqn{\tfrac{1}{\lambda}z + b}_{(\lambda - \lambda^2\mu_\lambda)\mI_d - \lambda^2(L_\lambda - \mu_\lambda)\mF^\top\mF}
    \\&=
    \frac{\mu_\lambda}{2}\sqn{z}
    +\frac{L_\lambda - \mu_\lambda}{2}\sqn{\mF z}
    -\<z,((1 - \lambda\mu_\lambda)\mI_d - \lambda(L_\lambda - \mu_\lambda)\mF^\top\mF)b>
    \\&
    -\frac{1}{2}\sqn{b}_{(\lambda - \lambda^2\mu_\lambda)\mI_d - \lambda^2(L_\lambda - \mu_\lambda)\mF^\top\mF}
    \\&\aeq{uses the assumption $\mF^\top\mF b = 0$}
    \frac{\mu_\lambda}{2}\sqn{z}
    +\frac{L_\lambda - \mu_\lambda}{2}\sqn{\mF z}
    -(1 - \lambda\mu_\lambda)\<z,b>
    -\frac{\lambda(1 - \lambda\mu_\lambda)}{2}\sqn{b}
    \\&\aeq{uses the definitions of $B_\lambda$ and $C_\lambda$}
    \frac{\mu_\lambda}{2}\sqn{z}
    +\frac{L_\lambda - \mu_\lambda}{2}\sqn{\mF z}
    -B_\lambda\<b,z>
    -C_\lambda,
\end{align*}
where \annotate, which concludes the proof.\qed

\subsubsection{Proof of Lemma~\ref{lem:geometric}}\label{sec:proof:geometric}
One can verify that vector $\hz$ satisfies the following linear system:
\begin{equation}
    \left(2\mu \mI_d + (L-\mu)\mF^\top\mF\right)\hz = B\basis{d}{1},
\end{equation}
which can be rewritten as follows:
\begin{equation}
    \begin{bNiceArray}{cccccc}
        (L+3\mu) & -(L-\mu)\\
        -(L-\mu) & 2(L+\mu) & -(L-\mu)\\
        &\Ddots&\Ddots&\Ddots\\
        &&-(L-\mu) & 2(L+\mu) & -(L-\mu)\\
        &&&& -(L-\mu) & (L+3\mu)\\
    \end{bNiceArray}
    \begin{bNiceArray}{c}
        \hz_1\\
        \hz_2\\
        \Vdots\\
        \hz_d
    \end{bNiceArray}
    =
    \begin{bNiceArray}{c}
        2B\\
        0\\
        \Vdots\\
        0
    \end{bNiceArray}.
\end{equation}
Let $z^\circ \in \R^d$ be defined as follows:
\begin{equation}\label{eq:z_circ}
    z^\circ = \frac{2B}{(1-\rho)(L-\mu)}\cdot(\rho,\rho^2,\ldots,\rho^{d}).
\end{equation}
Then one can observe that
\begin{equation}\label{eq:misc1}
    \left(2\mu \mI_d + (L-\mu)\mF^\top\mF\right)z^\circ = B\basis{d}{1}-B\rho^d\basis{d}{d},
\end{equation}
which implies
\begin{equation}
    \left(2\mu \mI_d + (L-\mu)\mF^\top\mF\right)(z^\circ - \hz) = -B\rho^d\basis{d}{d}.
\end{equation}
Hence, we obtain the following inequality:
\begin{equation}
    \norm{z^\circ - \hz}
    \leq
    \frac{B\rho^d}{2\mu},
\end{equation}
which concludes the proof.\qed
\subsubsection{Proof of Lemma~\ref{lem:q^n}}\label{sec:proof:q^n}
We can upper-bound $q^{-n}$ as follows:
\begin{align*}
    q^{-n}
     & =
    \exp\left(n\log\left(1 + \frac{1-q}{q}\right)\right)
    \\&\aleq{uses the concavity of the logarithm}
    \exp\left(\frac{n(1-q)}{q}\right)
    \\&\aleq{uses the definition of $n$ in \cref{eq:n}}
    \exp\left(\frac{L_{xy}}{6\mu_{yx}}\cdot\frac{1-q}{q}\right)
    \\&\aleq{uses the definition of $q$ in \cref{eq:q}}
    \exp\left(\frac{L_{xy}}{6\mu_{yx}}\cdot\frac{2\sqrt{\tdelta_x\smash{\mu_y}}}{\sqrt{4\alpha^2 + \tdelta_x\smash{\mu_y}} - \sqrt{\tdelta_x\smash{\mu_y}}}\right)
    \\&=
    \exp\left(\frac{L_{xy}}{6\mu_{yx}}\cdot\frac{2\sqrt{\tdelta_x\smash{\mu_y}}\big(\sqrt{4\alpha^2 + \tdelta_x\smash{\mu_y}} + \sqrt{\tdelta_x\smash{\mu_y}}\big)}{4\alpha^2}\right)
    \\&\aleq{uses the ineqality $\sqrt{a+b} \leq \sqrt{a} + \sqrt{b}$ for $a,b>0$}
    \exp\left(\frac{L_{xy}}{6\mu_{yx}}\cdot\frac{4\sqrt{\tdelta_x\smash{\mu_y}}\big(\alpha + \sqrt{\tdelta_x\smash{\mu_y}}\big)}{4\alpha^2}\right)
    \\&\aleq{uses the definition of $\tdelta_x$ in \cref{eq:tdelta_x}}
    \exp\left(\frac{L_{xy}}{6\mu_{yx}}\cdot\frac{\sqrt{\mu_x\smash{\mu_y} + 4\mu_{xy}^2\frac{\mu_y}{L_y}}\big(\alpha + \sqrt{\mu_x\smash{\mu_y} + 4\mu_{xy}^2\frac{\mu_y}{L_y}}\big)}{\alpha^2}\right)
    \\&\aleq{uses the assumption $\mu_x\mu_y \leq \mu_{yx}^2$ and the assumption $\mu_{yx} > 0$, which, together with \Cref{ass:xy}, implies $\mu_{xy} \leq \mu_{yx}$}
    \exp\left(\frac{L_{xy}}{6\mu_{yx}}\cdot\frac{\mu_{yx}\sqrt{1 + \frac{4\mu_y}{L_y}}\big(\alpha + \mu_{yx}\sqrt{1 + \frac{4\mu_y}{L_y}}\big)}{\alpha^2}\right)
    \\&\aleq{use \Cref{ass:Pi}}
    \exp\left(\frac{L_{xy}}{6\mu_{yx}}\cdot\frac{\mu_{yx}\sqrt{2}\big(\alpha + \mu_{yx}\sqrt{2}\big)}{\alpha^2}\right)
    \\&\aleq{uses the definition of $\alpha$ in \cref{eq:alpha_beta_gamma}}
    \exp\left(\frac{\sqrt{2}}{3}\left(1 + \frac{2\sqrt{2}\mu_{yx}}{L_{xy}}\right)\right)
    \\&\aleq{use \Cref{ass:Pi}}
    \exp\left(\frac{\sqrt{2}}{3} + \frac{2}{27}\right)
    \\&\leq 2,
\end{align*}
where \annotate, which concludes the proof.\qed

\subsubsection{Proof of Lemma~\ref{lem:omega_nu}}\label{sec:proof:omega_nu}

We can upper-bound $\omega_n$ as follows:
\begin{align*}
    \omega_n
     & \aeq{uses the definition of $\omega_n$ in \cref{eq:omega_nu}}
    \frac{(1-q)(1-q^{n-2})}{(1+q^{n-1})}
    \\&\aleq{uses the fact that $q \geq 0$}
    (1-q)(1-q^{n-2})
    =
    (1-q)^2\sum_{j=0}^{n-3}q^j
    \\&\aleq{uses the fact that $q \leq 1$}
    (n-2)(1-q)^2
    \\&\aeq{uses the definition of $q$ in \cref{eq:q}}
    \frac{4(n-2)\tdelta_x\smash{\mu_y}}{\big(\sqrt{4\alpha^2 + \tdelta_x\smash{\mu_y}}+\sqrt{\tdelta_x\smash{\mu_y}}\big)^2}
    \\&\leq
    \frac{(n-2)\tdelta_x\smash{\mu_y}}{\alpha^2},
\end{align*}
where \annotate.
Next, we can upper-bound $\frac{1}{\nu_n}$ as follows:
\begin{align*}
    \frac{1}{\nu_n}
     & \aeq{uses the definition of $\nu_n$ in \cref{eq:omega_nu}}
    \frac{q(q^{1-n} - q^{n-1})}{(1-q)(1+q)}
    =
    \frac{q^{2-n}}{(1+q)}\sum_{j=0}^{2n-3}q^j
    \aleq{uses the fact that $0 \leq q \leq 1$}
    2(n-1)q^{-n}
    \aleq{uses \Cref{lem:q^n}}
    4(n-1).
\end{align*}
where \annotate, which concludes the proof.\qed
\subsubsection{Proof of Lemma~\ref{lem:++}}\label{sec:proof:++}

First, we can lower-bound $\frac{L^+}{\mu^+}$ as follows:
\begin{align*}
    \frac{L^+}{\mu^+}
     & \aeq{uses the definition of $L^+$ and $\mu^+$ in \cref{eq:+}}
    \frac{
        \frac{\beta^2L_x}{\beta^2 + L_x\mu_y}
        +\frac{\tdelta_x}{n} + \frac{\omega_n\alpha^2}{n\mu_y}
    }{
        \frac{\beta^2\tdelta_x}{\beta^2 + \tdelta_x\mu_y}
        +\frac{\tdelta_x}{n} + \frac{\omega_n\alpha^2}{n\mu_y}
    }
    =
    1 +
    \frac{
        \beta^4(L_x-\tdelta_x)
    }{
        \left(
        \beta^2\tdelta_x
        +
        (\beta^2 + \tdelta_x\mu_y)\left(\frac{\tdelta_x}{n} + \frac{\omega_n\alpha^2}{n\mu_y}\right)
        \right)
        (\beta^2 + L_x\mu_y)
    }
    \\&\ageq{uses \cref{eq:tdelta_x_upper}}
    1 +
    \frac{
        \beta^4L_x
    }{
        2\left(
        \beta^2\tdelta_x
        +
        (\beta^2 + \tdelta_x\mu_y)\left(\frac{\tdelta_x}{n} + \frac{\omega_n\alpha^2}{n\mu_y}\right)
        \right)
        (\beta^2 + L_x\mu_y)
    }
    \\&\ageq{uses \Cref{lem:omega_nu}}
    1 +
    \frac{
        \beta^4L_x
    }{
        2\tdelta_x\left(
        \beta^2
        +
        \frac{(n-1)}{n}(\beta^2 + \tdelta_x\mu_y)
        \right)
        (\beta^2 + L_x\mu_y)
    }
    \\&\aeq{uses the definition of $\tdelta_x$ in \cref{eq:tdelta_x}}
    1 +
    \frac{
        \beta^4L_x
    }{
        2\tdelta_x\left(
        \beta^2
        +
        \frac{(n-1)}{n}\left(\beta^2 + \mu_x\mu_y + \mu_{xy}^2\frac{4\mu_y}{L_y}\right)
        \right)
        (\beta^2 + L_x\mu_y)
    }
    \\&\ageq{uses the assumption $\mu_x\mu_y \leq \mu_{yx}^2$ and the assumption $\mu_{yx} > 0$, which, together with \Cref{ass:xy}, implies $\mu_{xy} \leq \mu_{yx}$}
    1 +
    \frac{
        \beta^4L_x
    }{
        2\tdelta_x\left(
        \beta^2
        +
        \frac{(n-1)}{n}\left(\beta^2 + \mu_{yx}^2(1 + \frac{4\mu_y}{L_y})\right)
        \right)
        (\beta^2 + L_x\mu_y)
    }
    \\&\ageq{uses \Cref{ass:Pi}}
    1 +
    \frac{
        \beta^4L_x
    }{
        2\tdelta_x\left(
        \beta^2
        +
        \frac{(n-1)}{n}\left(\beta^2 + 2\mu_{yx}^2\right)
        \right)
        (\beta^2 + L_x\mu_y)
    }
    \\&\aeq{uses the definition of $\beta$ in \cref{eq:alpha_beta_gamma}}
    1 +
    \frac{
        \beta^4L_x
    }{
        2\tdelta_x\beta^4\left(
        1
        +
        \frac{(n-1)}{n}\left(1 + \frac{2n^2\mu_{yx}^2}{L_{xy}^2}\right)
        \right)
        \left(1 + \frac{n^2L_x\mu_y}{L_{xy}^2}\right)
    }
    \\&\ageq{uses the definition of $n$ in \cref{eq:n}}
    1 +
    \frac{
        L_x
    }{
        2\tdelta_x\left(
        1
        +
        \frac{(n-1)}{n}\left(1 + \frac{1}{18}\right)
        \right)
        \left(1 + \frac{L_x\mu_y}{36\mu_{yx}^2}\right)
    }
    \\&\geq
    1 +
    \frac{
        9L_x\mu_{yx}^2
    }{
        37\tdelta_x
        \left(\mu_{yx}^2 + \frac{1}{36}L_x\mu_y\right)
    }
    =
    \frac{
        \frac{9}{37}L_x\mu_{yx}^2
        +
        \tdelta_x
        \left(\mu_{yx}^2 + \frac{1}{36}L_x\mu_y\right)
    }{
        \tdelta_x
        \left(\mu_{yx}^2 + \frac{1}{36}L_x\mu_y\right)
    },
\end{align*}
where \annotate. Furthermore, we can lower-bound $\frac{L^+}{\mu^+}$ as follows:
\begin{align*}
    \frac{L^+}{\mu^+}
     & \geq
    1 +
    \frac{
        9L_x\mu_{yx}^2
    }{
        37\tdelta_x
        \left(\mu_{yx}^2 + \frac{1}{36}L_x\mu_y\right)
    }
    =
    1 +
    \frac{
        9\mu_{yx}^2
    }{
        37
        \left(\frac{\tdelta_x \mu_{yx}^2}{L_x} + \frac{1}{36}\tdelta_x\mu_y\right)
    }
    \\&\ageq{uses \cref{eq:tdelta_x_upper} and the definition of $\tdelta_x$ in \cref{eq:tdelta_x}}
    1 +
    \frac{
        9\mu_{yx}^2
    }{
        37
        \left(\frac{1}{2}\mu_{yx}^2 + \frac{1}{36}\left(\mu_x\mu_y + \frac{4\mu_y\mu_{xy}^2}{L_y}\right)\right)
    }
    \\&\ageq{uses the assumption $\mu_x\mu_y \leq \mu_{yx}^2$ and the assumption $\mu_{yx} > 0$, which, together with \Cref{ass:xy}, implies $\mu_{xy} \leq \mu_{yx}$}
    1 +
    \frac{
        9\mu_{xy}^2
    }{
        37
        \left(\frac{1}{2}\mu_{yx}^2 + \frac{1}{36}\mu_{yx}^2\left(1 + \frac{4\mu_y}{L_y}\right)\right)
    }
    \\&\ageq{uses \Cref{ass:Pi}}
    1 +
    \frac{
        9\mu_{yx}^2
    }{
        37
        \left(\frac{1}{2}\mu_{yx}^2 + \frac{1}{18}\mu_{yx}^2\right)
    }
    =
    1 +
    \frac{
        81
    }{
        185
    }
    \\&\geq
    \frac{
        10
    }{
        7
    },
\end{align*}
where \annotate.
Next, we can lower-bound $\frac{L^{++}}{\mu^{++}}$ as follows:
\begin{align*}
    \frac{L^{++}}{\mu^{++}}
     & \aeq{uses the definition of $L^{++}$ and $\mu^{++}$ in \cref{eq:++}}
    \frac{
        \frac{2\nu_n\alpha^2L^+}{2\nu_n\alpha^2 + n\mu_yL^+}
    }{
        \frac{2\nu_n\alpha^2\mu^+}{2\nu_n\alpha^2 + n\mu_y\mu^+}
    }
    =
    1 +
    \frac{
        \left(L^+ - \mu^+\right)
    }{
        \mu^+
        \left(1 + \frac{n\mu_y L^+}{2\nu_n\alpha^2}\right)
    }
    \\&\ageq{uses the previously obtained inequality $\mu^+ \leq \frac{7}{10}L^+$}
    1 +
    \frac{
        3L^+
    }{
        10\mu^+
        \left(1 + \frac{n\mu_y L^+}{2\nu_n\alpha^2}\right)
    }
    \\&\ageq{use \Cref{lem:omega_nu}}
    1 +
    \frac{
        3L^+
    }{
        10\mu^+
        \left(1 + \frac{2n(n-1)\mu_y L^+}{\alpha^2}\right)
    }
    \\&\ageq{uses the definition of $\alpha$ in \cref{eq:alpha_beta_gamma}}
    1 +
    \frac{
        3L^+
    }{
        10\mu^+
        \left(1 + \frac{8n(n-1)\mu_y L^+}{L_{xy}^2}\right)
    }
    \\&\ageq{uses the definition of $n$ in \cref{eq:n}}
    1 +
    \frac{
        3L^+
    }{
        10\mu^+
        \left(1 + \frac{2\mu_y L^+}{9\mu_{yx}^2}\right)
    }
    =
    1 +
    \frac{
        3\mu_{yx}^2
    }{
        10
        \left(\frac{\mu^+}{L^+}\mu_{yx}^2 + \frac{2}{9}\mu_y\mu^+\right)
    }
    \\&\ageq{uses the definition of $\mu^+$ in \cref{eq:+}}
    1 +
    \frac{
        3\mu_{yx}^2
    }{
        10
        \left(\frac{\mu^+}{L^+}\mu_{yx}^2 + \frac{2}{9}\mu_y\left(
        \frac{\beta^2\tdelta_x}{\tdelta_x\mu_y + \beta^2}
        + \frac{\tdelta_x}{n} + \frac{\omega_n\alpha^2}{n\mu_y}
        \right)\right)
    }
    \\&\leq
    1 +
    \frac{
        3\mu_{yx}^2
    }{
        10
        \left(\frac{\mu^+}{L^+}\mu_{yx}^2 + \frac{2}{9}\mu_y\left(
        \frac{(n+1)\tdelta_x}{n} + \frac{\omega_n\alpha^2}{n\mu_y}
        \right)\right)
    }
    \\&\ageq{use \Cref{lem:omega_nu}}
    1 +
    \frac{
        3\mu_{yx}^2
    }{
        10
        \left(\frac{\mu^+}{L^+}\mu_{yx}^2 + \frac{2}{9}\mu_y\left(
        \frac{(n+1)\tdelta_x}{n} + \frac{(n-2)\tdelta_x}{n}
        \right)\right)
    }
    \geq
    1 +
    \frac{
        3\mu_{yx}^2
    }{
        10
        \left(\frac{\mu^+}{L^+}\mu_{yx}^2 + \frac{4}{9}\mu_y\tdelta_x\right)
    }
    \\&\ageq{uses the previously obtained lower bound on $\frac{L^+}{\mu^+}$}
    1 +
    \frac{
        3\mu_{yx}^2
    }{
        10
        \left(
        \frac{
            \tdelta_x\mu_{yx}^2
            \left(\mu_{yx}^2 + \frac{1}{36}L_x\mu_y\right)
        }{
            \frac{9}{37}L_x\mu_{yx}^2
            +
            \tdelta_x
            \left(\mu_{yx}^2 + \frac{1}{36}L_x\mu_y\right)
        }
        + \frac{4}{9}\mu_y\tdelta_x\right)
    }
    \\&\geq
    1 +
    \frac{
        3\mu_{yx}^2
    }{
        10
        \left(
        \frac{
            37\tdelta_x
            \left(\mu_{yx}^2 + \frac{1}{36}L_x\mu_y\right)
        }{
            9L_x
        }
        + \frac{4}{9}\mu_y\tdelta_x\right)
    }
    \geq
    1 +
    \frac{
        27\mu_{yx}^2
    }{
        370\tdelta_x
        \left(
        \mu_y
        +
        \frac{
            \mu_{yx}^2
        }{
            L_x
        }
        \right)
    }
    \\&\aeq{uses the definition of $\delta_y$ in \cref{eq:xyxy}}
    1 +
    \frac{
        27\mu_{yx}^2
    }{
        370\tdelta_x\delta_y
    }
    \ageq{uses the definition of $\tdelta_x$ in \cref{eq:tdelta_x} and the definition of $\delta_x$ in \cref{eq:xyxy}}
    1 +
    \frac{
        27\mu_{yx}^2
    }{
        1480\delta_x\delta_y
    }
    \\&\geq
    1 +
    \frac{
        \mu_{yx}^2
    }{
        55\delta_x\delta_y
    },
\end{align*}
where \annotate.
Furthermore, we can lower-bound $\frac{L^{++}}{\mu^{++}}$ as follows:
\begin{align*}
    \frac{L^{++}}{\mu^{++}}
     & \geq
    1 +
    \frac{
        \mu_{yx}^2
    }{
        55\delta_x\delta_y
    }
    \ageq{uses the definitions of $\delta_x$ and $\delta_y$ in \cref{eq:xyxy}}
    1 +
    \frac{
        \mu_{yx}^2
    }{
        55\left(\mu_x\mu_y + \frac{\mu_y}{L_y}\mu_{xy}^2 + \frac{\mu_x}{L_x}\mu_{yx}^2 + \frac{\mu_{xy}^2\mu_{yx}^2}{L_xL_y}\right)
    }
    \\&\ageq{uses the assumption $\mu_x\mu_y \leq \mu_{yx}^2$ and the assumption $\mu_{yx} > 0$, which, together with \Cref{ass:xy}, implies $\mu_{xy} \leq \mu_{yx}$}
    1 +
    \frac{
        1
    }{
        55\left(1 + \frac{\mu_y}{L_y} + \frac{\mu_x}{L_x} + \frac{\mu_{xy}^2}{L_xL_y}\right)
    }
    \ageq{uses \Cref{ass:Pi}}
    1 +
    \frac{
        16
    }{
        25\cdot 55
    }
    \geq
    1 +
    \frac{
        1
    }{
        86
    },
\end{align*}
where \annotate, which concludes the proof.\qed
\subsubsection{Proof of Lemma~\ref{lem:geometric2}}\label{sec:proof:geometric2}

First, we can upper-bound $\sqn{x' - x^*}$ as follows:
\begin{align*}
    \sqn{x' - x^*}
     & \aeq{uses the expressions for $x^*$ in \cref{eq:x1,eq:x2,eq:x3} and the definition of $x'$ in \cref{eq:x_prime}}
    \sum_{i=1}^n \sqn{\prox_{\frac{\mu_y}{\beta^2}h_1}(x_{n+1}') - \prox_{\frac{\mu_y}{\beta^2}h_1}(x_{n+1}^*)}
    \\&
    +\sqn{\prox_{\frac{n\mu_y}{2\nu_n\alpha^2} h_1^+}(z^\circ) - \prox_{\frac{n\mu_y}{2\nu_n\alpha^2} h_1^+}(z^*)}
    \\&
    +\sum_{i=n+2}^{2n-1}\sqn*{\frac{(x_{n+1}'-x_{n+1}^*)(q^{2n-i} - q^{i-2n}) + (x_{2n}'-x_{2n}^*)(q^{i-(n+1)} - q^{(n+1)-i})}{(q^{n-1} - q^{1-n})}}
    \\&
    +\sqn{\prox_{\frac{n\mu_y}{2\nu_n\alpha^2} h_2^+}(z^\circ) - \prox_{\frac{n\mu_y}{2\nu_n\alpha^2} h_2^+}(z^*)}
    \\&
    +\sum_{i=2n+1}^{3n}\sqn{\prox_{\frac{\mu_y}{\beta^2}h_2}(x_{2n}') - \prox_{\frac{\mu_y}{\beta^2}h_2}(x_{2n}^*)}
    \\&\aleq{use the nonexpansiveness of the proximal operator}
    n \sqn{x_{n+1}' - x_{n+1}^*}
    +\sqn{z^\circ - z^*}
    \\&
    +\sum_{i=n+2}^{2n-1}\sqn*{\frac{(x_{n+1}'-x_{n+1}^*)(q^{2n-i} - q^{i-2n}) + (x_{2n}'-x_{2n}^*)(q^{i-(n+1)} - q^{(n+1)-i})}{(q^{n-1} - q^{1-n})}}
    \\&
    +\sqn{z^\circ - z^*}
    +n\sqn{x_{2n}' - x_{2n}^*}
    \\&\aleq{uses Young's inequality}
    2\sqn{z^\circ - z^*}
    +n (\sqn{x_{n+1}' - x_{n+1}^*} + \sqn{x_{2n}' - x_{2n}^*})
    \\&
    +2\sum_{i=n+2}^{2n-1}\frac{(q^{2n-i} - q^{i-2n})^2}{(q^{n-1} - q^{1-n})^2}\sqn{x_{n+1}'-x_{n+1}^*}
    \\&
    +2\sum_{i=n+2}^{2n-1}\frac{(q^{i-(n+1)} - q^{(n+1)-i})^2}{(q^{n-1} - q^{1-n})^2}\sqn{x_{2n}'-x_{2n}^*}
    \\&=
    2\sqn{z^\circ - z^*}
    +\left(n + 2\sum_{i=1}^{n-2}\frac{(q^{i} - q^{-i})^2}{(q^{n-1} - q^{1-n})^2}\right) (\sqn{x_{n+1}' - x_{n+1}^*} + \sqn{x_{2n}' - x_{2n}^*})
    \\&\aeq{uses the definitions of $x_{n+1}^*$ and $x_{2n}^*$ in \cref{eq:x3} and the definitions of $x_{n+1}'$ and $x_{2n}'$ in \cref{eq:x_prime}}
    2\sqn{z^\circ - z^*}
    \\&
    +\left(n + 2\sum_{i=1}^{n-2}\frac{(q^{i} - q^{-i})^2}{(q^{n-1} - q^{1-n})^2}\right) \sqn{\prox_{\frac{n\mu_y}{2\nu_n\alpha^2} h_1^+}(z^\circ) - \prox_{\frac{n\mu_y}{2\nu_n\alpha^2} h_1^+}(z^*)}
    \\&
    +\left(n + 2\sum_{i=1}^{n-2}\frac{(q^{i} - q^{-i})^2}{(q^{n-1} - q^{1-n})^2}\right) \sqn{\prox_{\frac{n\mu_y}{2\nu_n\alpha^2} h_2^+}(z^\circ) - \prox_{\frac{n\mu_y}{2\nu_n\alpha^2} h_2^+}(z^*)}
    \\&\aleq{use the nonexpansiveness of the proximal operator}
    \left(2 + 2n + 4\sum_{i=1}^{n-2}\frac{(q^{i} - q^{-i})^2}{(q^{n-1} - q^{1-n})^2}\right) \sqn{z^\circ - z^*}
    \\&\aleq{uses \cref{eq:z_ub}}
    \left(2 + 2n + 4\sum_{i=1}^{n-2}\frac{(q^{i} - q^{-i})^2}{(q^{n-1} - q^{1-n})^2}\right)\left(\frac{A^{++}\rho^d}{2\mu^{++}}\right)^2
    \\&\aeq{uses the definition of $A^{++}$ in \cref{eq:A}}
    C_\pi'A^2\rho^{2d},
\end{align*}
where \annotate, and $C_\pi' \geq 0$ is a constant that depends on the parameters $\pi \in \Pi$.

Next, we can upper-bound $\distsol(x',y')$ as follows:
\begin{align*}
    \distsol(x',y')
     & \aeq{uses the definition of $\distsol$ in \cref{eq:distsol}}
    \delta_x\sqn{x' - x^*} + \delta_y\sqn{y' - y^*}
    \\&\aeq{uses the definition of $y'$ in \cref{eq:x_prime} and the expression for $y^*$ in \cref{eq:y}}
    \delta_x\sqn{x' - x^*} + \delta_y\sqn{\nabla g^*(\mB x') - \nabla g^*(\mB x^*)}
    \\&\aleq{uses the $(1/\mu_y)$-smoothness of function $g^*(y)$}
    \delta_x\sqn{x' - x^*} + \frac{\delta_y}{\mu_y^2}\sqn{\mB (x'-x^*)}
    \\&\aleq{uses \Cref{ass:xy}}
    \delta_x\sqn{x' - x^*} + \frac{\delta_yL_{xy}^2}{\mu_y^2}\sqn{x'-x^*}
    \\&\aleq{uses the previously obtained upper bound on $\sqn{x'-x^*}$}
    \left(\delta_x + \frac{\delta_yL_{xy}^2}{\mu_y^2}\right)C_{\pi}'A^2 \rho^{2d}
    \\&=
    C_\pi'' A^2 \rho^{2d},
\end{align*}
where \annotate, and $C_\pi'' \geq 0$ is a constant that depends on the parameters $\pi \in \Pi$.

Furthermore, we can express $x_{2n}'$ as follows:
\begin{align*}
    x_{2n}'
     & \aeq{uses the definition of $x'$ on \cref{eq:x_prime}}
    \prox_{\frac{n\mu_y}{2\nu_n\alpha^2} h_2^+}(z^\circ)
    \\&\aeq{uses the properties of the proximal operator and the definition of $h_2^{++}(z)$ in \cref{eq:h++2}}
    z^\circ - \frac{n\mu_y}{2\nu_n\alpha^2}\nabla h_2^{++}(z^\circ)
    \\&\aeq{uses the definition of $h_2^{++}(z)$ in \cref{eq:h++}}
    \left(\mI_d - \frac{n\mu_y}{2\nu_n\alpha^2}\left(\mu^{++}\mI_d + (L^{++} - \mu^{++})\mF_2^\top\mF_2\right) \right)z^\circ + A^{++}\basis{d}{1}
    \\&=
    \left(1 - \frac{n\mu_y\mu^{++}}{2\nu_n\alpha^2}\right)z^\circ - \frac{n\mu_y(L^{++} - \mu^{++})}{2\nu_n\alpha^2}\mF_2^\top\mF_2 z^\circ+ A^{++}\basis{d}{1}
    \\&\aeq{uses the adapted version of the definition of $z^\circ$ in \cref{eq:z_circ}}
    \frac{2A^{++}}{(1-\rho)(L^{++} - \mu^{++})}
    \left(
    \left(1 - \frac{n\mu_y\mu^{++}}{2\nu_n\alpha^2}\right)\mI_d - \frac{n\mu_y(L^{++} - \mu^{++})}{2\nu_n\alpha^2}\mF_2^\top\mF_2\right) (\rho,\ldots,\rho^d)
    \\&
    + A^{++}\basis{d}{1}
    \\&\aeq{uses the definition of $A^{++}$ in \cref{eq:A}}
    A(C_{\pi,2n}' \mI_d + 2C_{\pi,2n}''\mF_2^\top\mF_2)(\rho,\ldots,\rho^d) + A C_{\pi,2n}'''\basis{d}{1}
    \\&\aeq{uses the definition of matrix $\mF_2$ in \cref{eq:F}}
    A C_{\pi,2n}'''\basis{d}{1}
    +AC_{\pi,2n}'(\rho,\ldots,\rho^d)
    \\&
    +AC_{\pi,2n}''(0, \rho^2(1 - \rho),\rho^2(\rho - 1),\rho^4(1 - \rho),\rho^4(\rho-1),\ldots)
    \\&=
    A(C_{\pi,2n}''' + \rho C_{\pi,2n}')(1,0,0,\ldots)
    \\&
    +A(\rho^2C_{\pi,2n}' + \rho^2(1-\rho)C_{\pi,2n}'')(0,1,0,\rho^2,0,\rho^4,\ldots)
    \\&
    +A(\rho^3 C_{\pi,2n}' + \rho^2(\rho-1)C_{\pi,2n}'')(0,0,1,0,\rho^2,0,\rho^4,\ldots),
\end{align*}
where \annotate, and $C_{\pi,2n}',C_{\pi,2n}''\in \R$ are constants that depend on the parameters $\pi \in \Pi$. From this expression, we can conclude that the following inequality holds:
\begin{equation}
    \sqn{x_{2n}' - x_{2n}^\circ} \leq A^2 C_{\pi,2n}\rho^{2d},
\end{equation}
where $C_{\pi,2n}\geq 0$ is some constant that depends on the parameters $\pi \in \Pi$.
Similarly, we can obtain the following inequality:
\begin{equation}
    \delta_x \sqn{x' - x^\circ} + \delta_y \sqn{y' - y^\circ} \leq C_\pi''' A^2 \rho^{2d},
\end{equation}
where $C_{\pi}'''\geq 0$ is some constant that depends on the parameters $\pi \in \Pi$. Therefore, we obtain the following inequality:
\begin{align*}
    \distsol(x^\circ,y^\circ)
     & \leq
    2 \distsol(x',y') + 2\delta_x \sqn{x' - x^\circ} + 2\delta_y \sqn{y' - y^\circ}
    \\&
    \leq 2(C_\pi''' + C_\pi'')A^2\rho^{2d}
    \\&=
    C_\pi A^2\rho^{2d},
\end{align*}
which concludes the proof.\qed

\subsubsection{Proof of Lemma~\ref{lem:init_dist}}\label{sec:proof:init_dist}

Functions $f(x)$ and $g(y)$ defined in \cref{eq:f,eq:g} are quadratic. Hence, the optimality conditions~\eqref{eq:opt} can be written in the following form:
\begin{equation}
    \begin{aligned}
        \mG_x x^* + \mB^\top y^* & = Ab \\
        \mG_y y^* - \mB x^*      & = 0,
    \end{aligned}
\end{equation}
where $b = (\zeros_{2n},\ones_n)\otimes \basis{d}{1} \in \R^{n_xd}$, and $\mG_x \in \R^{n_xd\times n_xd}$ and $\mG_y \in \R^{n_yd\times n_yd}$ are symmetric matrices that can be expressed from $d\in\{1,2,\ldots\}$ and the parameters $\pi \in \Pi$ using the definitions of functions $f(x)$ and $g(y)$ in \cref{eq:f,eq:g}.
Hence, the solution $(x^*,y^*)$ can be expressed as follows:
\begin{equation}
    \begin{aligned}
        x^* & = A\cdot\mQ_x^{-1}b,             \\
        y^* & = A\cdot G_y^{-1}\mB\mQ_x^{-1}b,
    \end{aligned}
\end{equation}
where $\mQ_x \in \R^{n_xd\times n_xd}$ is defined as follows:
\begin{equation}
    \mQ_x = \mG_x + \mB^\top \mG_y^{-1}\mB.
\end{equation}
Note that matrix $\mG_x + \mB^\top \mG_y^{-1}\mB$ is invertible due to \Cref{ass:x,,ass:y,,ass:xy,,ass:xyxy}:
\begin{equation}
    \mG_x + \mB^\top \mG_y^{-1}\mB \succeq \delta_x\mI_{n_xd} \succ \mO_{n_xd}.
\end{equation}
Furthermore, using the definition of $(x^0,y^0)$ in \cref{eq:x_init}, we obtain the following:
\begin{align*}
    (x^* - x^0,y^* - y^0)
     & =
    (\mI_{n_x+n_y}\otimes(\mI_d - \mP))(x^*,y^*)
    \\&=
    (\mI_{n_x+n_y}\otimes(\mI_d - \mP))(\mQ_x^{-1}b,G_y^{-1}\mB\mQ_x^{-1}b) \cdot A.
\end{align*}
Hence, it is easy to observe that $\distsol(x^0,y^0) = B_{\pi,d}A^2$ for some constant $B_{\pi,d} \geq 0$ that depends on $d \in \{2,3,\ldots\}$ and the parameters $\pi \in \Pi$. Finally, we need to show the existence of $\bar{d} \in \{2,3,\ldots\}$ such that for all $A > 0$ the following inequality holds:
\begin{equation}
    \min_{d \in \{\hat{d},\hat{d}+1,\ldots\}} B_{\pi,d} = \min_{d \in \{\hat{d},\hat{d}+1,\ldots\}}  \frac{\distsol(x^0,y^0)}{A^2} > 0.
\end{equation}
To do this, we can lower-bound $\distsol(x^0,y^0)$ as follows:
\begin{align*}
    \distsol(x^0,y^0)
     & \geq
    \delta_x\left(\sqn{x_{n+1}^0 - x_{n+1}^*} + \sqn{x_{2n}^0 - x_{2n}^*}\right)
    \\&\aeq{uses the definition of $x^0$ in \cref{eq:x_init}}
    \delta_x\left(\sqn{(\mI_d - \mP)x_{n+1}^*} + \sqn{(\mI_d - \mP)x_{2n}^*}\right)
    \\&\ageq{uses the convexity of $\sqn{\cdot}$}
    \tfrac{1}{2}\delta_x\sqn{(\mI_d - \mP)(x_{n+1}^* + x_{2n}^*)}
    \\&\aeq{uses the expressions for $x_{n+1}^*$ and $x_{2n}^*$ in \cref{eq:x3}}
    \tfrac{1}{2}\delta_x\sqn{(\mI_d - \mP)(\prox_{\frac{n\mu_y}{2\nu_n\alpha^2} h_1^+}(z^*) + \prox_{\frac{n\mu_y}{2\nu_n\alpha^2} h_2^+}(z^*))}
    \\&\aeq{uses the properties of the proximal operator and the definition of $h_j^{++}(z)$ in \cref{eq:h++2}}
    \tfrac{1}{2}\delta_x\sqn{(\mI_d - \mP)(2z^* - \tfrac{n\mu_y}{2\nu_n\alpha^2}(\nabla h_1^{++}(z^*) + \nabla h_2^{++}(z^*)))}
    \\&\aeq{uses the definition of $z^*$ in \cref{eq:z}}
    2\delta_x\sqn{(\mI_d - \mP)z^*}
    \\&\ageq{uses Young's inequality}
    \delta_x\sqn{(\mI_d - \mP)z^\circ}
    -2\delta_x\sqn{(\mI_d - \mP)(z^* - z^\circ)}
    \\&\ageq{uses \cref{eq:z_ub}}
    \delta_x\sqn{(\mI_d - \mP)z^\circ}
    -2\delta_x\left(\frac{A^{++}\rho^d}{2\mu^{++}}\right)^2
    \\&\ageq{uses the adapted version of the definition of $z^\circ$ in \cref{eq:z_circ} in the proof of \Cref{lem:geometric} in \Cref{sec:proof:geometric}}
    \delta_x\left(\frac{2A^{++}\rho^2}{(1-\rho)(L^{++} - \mu^{++})}\right)^2
    -2\delta_x\left(\frac{A^{++}\rho^d}{2\mu^{++}}\right)^2
    \\&\geq
    (A^{++})^2\left(
    \frac{2 \delta_x\rho^{4}}{(L^{++})^2}
    -\frac{ \delta_x\rho^{2d}}{2(\mu^{++})^2}
    \right)
    \\&\ageq{uses the definition of $A^{++}$ in \cref{eq:A}}
    A^2 \cdot\left(\frac{2\nu_n\alpha^2 \beta^2}{(2\nu_n\alpha^2 + n\mu_y\mu^+)(\beta^2 + \mu_y\tdelta_x)}\right)^2\cdot
    \left(
    \frac{ 2\delta_x\rho^{4}}{(L^{++})^2}
    -\frac{ \delta_x\rho^{2d}}{2(\mu^{++})^2}
    \right),
\end{align*}
where \annotate.
It is not hard to verify that $\left(\frac{2\nu_n\alpha^2 \beta^2}{(2\nu_n\alpha^2 + n\mu_y\mu^+)(\beta^2 + \mu_y\tdelta_x)}\right)^2$ and $\frac{2 \delta_x\rho^{4}}{(L^{++})^2}$ are positive constants that do not depend on $d$. Hence, we can easily obtain the following relation:
\begin{equation}
    \liminf_{d \to + \infty}\frac{\distsol(x^0,y^0)}{A^2} > 0,
\end{equation}
which concludes the proof.\qed

\newpage

\section{Proof of Lemma~\ref{lem:lb_B2}}\label{sec:proof:lb_B2}
\newcommand{\bmu}{\bar{\mu}}
Let $\bmu_{xy} = \max\{\mu_{xy},\mu_{yx}\}$.
We consider a special instance of problem~\eqref{eq:main}, where $\sX = \sY = \R^d$, functions $f(x)$ and $g(y)$ are defined as follows:
\begin{equation}
    f(x) = \frac{\mu_x}{2}\sqn{x} - A\<\basis{d}{1},x>,\quad
    g(y) = \frac{\mu_y}{2}\sqn{y},
\end{equation}
and matrix $\mB \in \R^{d\times d}$ is defined as follows:
\begin{equation}\label{eq:B2}
    \mB = \frac{1}{2}\begin{bNiceMatrix}[r]
        \alpha & -\beta &                 \\
               & \Ddots & \Ddots          \\
               &        & \alpha & -\beta \\
               &        &        & \alpha
    \end{bNiceMatrix},
    \quad\text{where}\quad
    \alpha = L_{xy} + \bmu_{xy},\quad
    \beta = L_{xy} - \bmu_{xy}.
\end{equation}
Functions $f(x)$ and $g(y)$ obviously satisfy \Cref{ass:x,ass:y}.
Moreover, matrix $\mB$ satisfies \Cref{ass:xy} due to the following \Cref{lem:B2}, the proof is available in \Cref{sec:proof:B2}.
\begin{lemma}\label{lem:B2}
    The singular values of matrix $\mB$ defined in \cref{eq:B2} satisfy the following inequalities:
    \begin{equation}
        \bmu_{xy} \leq \smin(\mB) \leq \smax(\mB) \leq L_{xy}.
    \end{equation}
\end{lemma}

Next, we establish the following \Cref{lem:hard_sol2}, which describes the solution to the problem defined above, the proof is available in \Cref{sec:proof:hard_sol2}.
\begin{lemma}\label{lem:hard_sol2}
    For all $d \in \{2,3,\ldots\}$, the instance of problem~\eqref{eq:main} defined above has a unique solution $(x^*,y^*) \in \sX \times \sY$. Moreover, there exist vectors $x^\circ,y^\circ \in \spanset(\{(1,\rho,\ldots,\rho^{d-1})\})$ such that the following inequality holds:
    \begin{equation}
        \distsol(x^\circ,y^\circ) \leq C_\pi A^2\rho^{2d},
    \end{equation}
    where $C_\pi > 0$ is some constant that possibly depends on the parameters $\pi \in \Pi$, but does not depend on $d$,
    and $\rho \in (0,1)$ satisfies the following inequality:
    \begin{equation}
        \rho \geq \max\left\{1 - \frac{\sqrt{8}}{\sqrt{\condxya}}, \frac{4}{5}\right\}.
    \end{equation}
    In addition, the initial distance to the solution $\distsol(0,0)$ is a quadratic function of $A$, that is,
    \begin{equation}
        \distsol(0,0) = B_{\pi,d}A^2,
    \end{equation}
    where $B_{\pi,d} > 0$ is a constant that possibly depends on $d \in \{1,2,\ldots\}$ and the parameters $\pi \in \Pi$, and satisfies the inequality $\min_{d\in \{1,2,\ldots\}} B_{\pi,d} > 0$.
\end{lemma}

We can express $(x^\circ,y^\circ)$ from \Cref{lem:hard_sol2} as follows:
\begin{equation}\label{eq:x_circ2}
    x^\circ = u_x(1,\rho,\ldots,\rho^{d-1}),\quad
    y^\circ = u_y(1,\rho,\ldots,\rho^{d-1}),
    \quad\text{where}\quad
    u_x,u_y \in \R,
\end{equation}
which implies
\begin{equation}
    \sqn{x^\circ} = \frac{u_x^2(1-\rho^{2d})}{1-\rho^2},\quad
    \sqn{y^\circ} = \frac{u_y^2(1-\rho^{2d})}{1-\rho^2}.
\end{equation}
Further, we fix $k \in \rng{0}{d}$.
Using the sparse structure of matrix $\mB$ defined in \cref{eq:B2}, and using the standard arguments \citep{nesterov2013introductory,ibrahim2020linear,zhang2022near,scaman2017optimal,scaman2018optimal,kovalev2024lower}, we can show that the output vector $(x_o(\tau),y_o(\tau)) \in \sX \times \sY$ satisfies the following implication:
\begin{equation}\label{eq:span2}
    \tau \leq D \cdot\tau_\mB k
    \quad\Rightarrow\quad
    x_o(\tau),y_o(\tau) \in \begin{cases}
        \spanset(\{\basis{d}{1},\ldots,\basis{d}{k}\}) & k > 0 \\
        \{0\}                                          & k=0
    \end{cases},
\end{equation}
where $D > 0$ is a universal constant.
The right-hand side of this implication implies the following:
\begin{align*}
    \distsol(x_o(\tau), y_o(\tau))
     & =
    \delta_x \sqn{x_o(\tau) - x^*}
    +\delta_y \sqn{y_o(\tau) - y^*}
    \\&\ageq{use Young's inequality}
    \tfrac{1}{2}\delta_x \sqn{x_o(\tau) - x^\circ}
    +\tfrac{1}{2}\delta_y \sqn{y_o(\tau) - y^\circ}
    -\distsol(x^\circ,y^\circ)
    \\&\ageq{uses \cref{eq:span2} and the expression for $(x^\circ,y^\circ)$ in \cref{eq:x_circ2}}
    \tfrac{1}{2}(\delta_x u_x^2 + \delta_y u_y^2)\sum_{j=k}^{d-1}\rho^{2j}
    -\distsol(x^\circ,y^\circ)
    \\&=
    \tfrac{1}{2}(\delta_x u_x^2 + \delta_y u_y^2)\frac{\rho^{2k}-\rho^{2d}}{1-\rho^2}
    -\distsol(x^\circ,y^\circ)
    \\&=
    \tfrac{1}{2}(\delta_x \sqn{x^\circ} + \delta_y \sqn{y^\circ})\frac{\rho^{2k}-\rho^{2d}}{1-\rho^{2d}}
    -\distsol(x^\circ,y^\circ)
    \\&\ageq{use Young's inequality}
    \left(\tfrac{1}{4}\distsol(0,0) -  \tfrac{1}{2}\distsol(x^\circ,y^\circ)\right)\frac{\rho^{2k}-\rho^{2d}}{1-\rho^{2d}}
    -\distsol(x^\circ,y^\circ)
    \\&=
    \frac{\rho^{2k}-\rho^{2d}}{4(1-\rho^{2d})}\distsol(0,0)
    -\left(1+\frac{\rho^{2k}-\rho^{2d}}{2(1-\rho^{2d})}\right)\distsol(x^\circ,y^\circ)
    \\&\ageq{uses \Cref{lem:hard_sol2}}
    \frac{\rho^{2k}-\rho^{2d}}{4(1-\rho^{2d})}B_{\pi,d}A^2
    -\left(1+\frac{\rho^{2k}-\rho^{2d}}{2(1-\rho^{2d})}\right)C_\pi A^2\rho^{2d},
\end{align*}
where \annotate. Furthermore, we can choose $A = R / \sqrt{B_{\pi,d}}$ to ensure the initial distance $\distsol(0,0) = R^2$ and obtain the following:
\begin{align*}
    \distsol(x_o(\tau), y_o(\tau))
     & =
    \frac{\rho^{2k}-\rho^{2d}}{4(1-\rho^{2d})}B_{\pi,d} A^2
    -\left(1+\frac{\rho^{2k}-\rho^{2d}}{2(1-\rho^{2d})}\right)C_\pi A^2\rho^{2d}
    \\&\aeq{uses the choice of $A$ above}
    \frac{\rho^{2k}-\rho^{2d}}{4(1-\rho^{2d})}R^2
    -\left(1+\frac{\rho^{2k}-\rho^{2d}}{2(1-\rho^{2d})}\right)\frac{C_\pi R^2\rho^{2d}}{B_{\pi,d}}
    \\&\ageq{uses the inequalities $B_{\pi,d} \geq \min_{d' \in \{1,2,\ldots\}} B_{\pi,d'} > 0$ implied by \Cref{lem:hard_sol2}}
    \frac{\rho^{2k}-\rho^{2d}}{4(1-\rho^{2d})}R^2
    -\left(1+\frac{\rho^{2k}-\rho^{2d}}{2(1-\rho^{2d})}\right)\frac{C_\pi R^2\rho^{2d}}{\min_{d' \in \{1,2,\ldots\}} B_{\pi,d'}}
    \\&=
    \rho^{2k}R^2 \cdot\left(\frac{1-\rho^{2(d-k)}}{4(1-\rho^{2d})}
    -\left(1+\frac{\rho^{2k}-\rho^{2d}}{2(1-\rho^{2d})}\right)\frac{C_\pi \rho^{2(d-k)}}{\min_{d' \in \{1,2,\ldots\}} B_{\pi,d'}}\right)
    \\&\ageq{is implied by choosing a large enough value of $d \in \{1,2,\ldots\}$}
    \tfrac{1}{5}\rho^{2k}R^2,
\end{align*}
where \annotate.

Next, we consider the case $\epsilon \leq \frac{1}{5}R^2$.
Choosing $k$ as follows:
\begin{equation}\label{eq:k2}
    k =  \floor*{\frac{\log\left(\frac{1}{5}R^2/\epsilon\right)}{2\log \left(1/\rho\right)}}
\end{equation}
implies $ \distsol(x_o(\tau), y_o(\tau)) \geq \epsilon$. Therefore, by contraposition, from implication~\eqref{eq:span2}, we obtain the following:
\begin{equation}
    \distsol(0,0) < \epsilon
    \quad\Rightarrow\quad
    \tau > D \cdot \tau_\mB k.
\end{equation}
The right-hand side of this implication implies the following lower bound on the execution time $\tau$:
\begin{align*}
    \tau
     & >
    D \cdot \tau_\mB k
    \\&\aeq{uses the definition of $k$ in \cref{eq:k2}}
    D \cdot \tau_\mB \floor*{\frac{\log\left(\frac{1}{5}R^2/\epsilon\right)}{2\log \left(1/\rho\right)}}
    \\&\geq
    D \cdot \tau_\mB \cdot\frac{\log\left(\frac{1}{5}R^2/\epsilon\right) - 2\log \left(1/\rho\right)}{2\log \left(1/\rho\right)}
    \\&\ageq{use \Cref{lem:hard_sol2}}
    D \cdot \tau_\mB \cdot\frac{\log\left(\frac{1}{5}R^2/\epsilon\right) - 2\log \left(5/4\right)}{2\log \left(1/\rho\right)}
    \\&=
    \frac{D}{2} \cdot \tau_\mB \cdot\frac{\log\left(\frac{16}{125}R^2/\epsilon\right)}{\log \left(1/\rho\right)}
    \\&\ageq{uses the concavity of the logarithm}
    \frac{D}{2} \cdot \tau_\mB \cdot\frac{\log\left(\frac{16}{125}R^2/\epsilon\right)}{\left(1/\rho - 1\right)}
    \\&=
    \frac{D}{2} \cdot \frac{\tau_\mB\rho}{1-\rho}\cdot\log\left(\tfrac{16R^2}{125\epsilon}\right)
    \\&\ageq{use \Cref{lem:hard_sol2}}
    \frac{D}{5\sqrt{2}} \cdot \tau_\mB \sqrt{\condxya}\cdot\log\left(\tfrac{16R^2}{125\epsilon}\right)
    \\&=
    \Omega\left(\tau_\mB \cdot \sqrt{\condxya}\log\left(\tfrac{16R^2}{125\epsilon}\right)\right),
\end{align*}
where \annotate, which concludes the proof in the case $\epsilon \leq \frac{1}{5}R^2$. In the remaining case $\epsilon >  \frac{1}{5}R^2$, we have $\log\left(\tfrac{16R^2}{125\epsilon}\right) \leq 0$. Hence, the latter lower bound holds due to the fact that $\tau \geq 0$, which concludes the proof.
\qed

\subsection{Proofs of Auxiliary Lemmas}

\subsubsection{Proof of Lemma~\ref{lem:B2}}\label{sec:proof:B2}

We can write matrix $\mB^\top \mB$ as follows:
\begin{align*}
    \mB^\top\mB
     & =
    \frac{1}{4}\begin{bNiceMatrix}[l]
                   \alpha^2     & -\alpha\beta     &              &                  &                  \\
                   -\alpha\beta & \alpha^2+\beta^2 & -\alpha\beta &                  &                  \\
                                & \Ddots           & \Ddots       & \Ddots           &                  \\
                                &                  & -\alpha\beta & \alpha^2+\beta^2 & -\alpha\beta     \\
                                &                  &              & -\alpha\beta     & \alpha^2+\beta^2
               \end{bNiceMatrix}
    \\&=
    \frac{1}{4}\begin{bNiceMatrix}[l]
                   \alpha\beta  & -\alpha\beta &              &              &              \\
                   -\alpha\beta & 2\alpha\beta & -\alpha\beta &              &              \\
                                & \Ddots       & \Ddots       & \Ddots       &              \\
                                &              & -\alpha\beta & 2\alpha\beta & -\alpha\beta \\
                                &              &              & -\alpha\beta & 2\alpha\beta
               \end{bNiceMatrix}
    +
    \frac{1}{4}\begin{bNiceMatrix}
                   \alpha(\alpha-\beta)                            \\
                    & (\alpha-\beta)^2                             \\
                    &                  & \Ddots                    \\
                    &                  &        & (\alpha-\beta)^2
               \end{bNiceMatrix}.
\end{align*}
Therefore, we can obtain the following matrix inequalities:
\begin{equation}
    \begin{aligned}
        \mB^\top\mB & \succeq \tfrac{1}{4}\min\{(\alpha-\beta)^2,\alpha(\alpha-\beta)\}\mI_d =
        \bmu_{yx}^2 \mI_d,
        \\
        \mB^\top\mB & \preceq \tfrac{1}{4}\max\left\{(\alpha-\beta)^2 + 4\alpha\beta,\alpha(\alpha-\beta) + 2\alpha\beta\right\} = L_{xy}^2 \mI_d,
    \end{aligned}
\end{equation}
which conclude the proof.\qed

\subsubsection{Proof of Lemma~\ref{lem:hard_sol2}}\label{sec:proof:hard_sol2}

The first-order optimality conditions~\eqref{eq:opt} imply that the unique solution to the problem $(x^*,y^*)\in \sX \times \sY$ is defined by the following linear system:
\begin{equation}\label{eq:x4}
    \begin{aligned}
        (\mu_x\mu_y\mI_d  + \mB^\top\mB) x^* & = \mu_y A\basis{d}{1}, \\
        \mu_y y^* - \mB x^*                  & = 0.
    \end{aligned}
\end{equation}
Hence, we can express the initial distance $\distsol(0,0)$ as follows:
\begin{align*}
    \distsol(0,0)
     & =
    \delta_x\sqn{x^*} + \delta_y\sqn{y^*}
    \\&\aeq{uses the linear system~\eqref{eq:x4} and the fact that the matrix $(\mu_x\mu_y\mI_d + \mB^\top \mB)$ is invertible, which is implied by the matrix inequality $(\mu_x\mu_y\mI_d + \mB^\top \mB)\succeq (\mu_x\mu_y + \bmu_{xy}^2) \mI_d$ and \Cref{ass:xyxy}}
    \delta_x\sqn{\mu_yA(\mu_x\mu_y\mI_d  + \mB^\top\mB)^{-1}\basis{d}{1}} + \delta_y\sqn{A\mB(\mu_x\mu_y\mI_d  + \mB^\top\mB)^{-1}\basis{d}{1}}
    \\&=
    A^2 B_{\pi,d},
\end{align*}
where \annotate, and where $B_{\pi,d}$ is defined as follows:
\begin{equation}
    B_{\pi,d} = \delta_x\mu_y^2\sqn{(\mu_x\mu_y\mI_d  + \mB^\top\mB)^{-1}\basis{d}{1}} + \delta_y\sqn{\mB(\mu_x\mu_y\mI_d  + \mB^\top\mB)^{-1}\basis{d}{1}}.
\end{equation}
We can also lower-bound $B_{\pi,d}$ as follows:
\begin{align*}
    B_{\pi,d}
     & =
    \delta_x\mu_y^2\sqn{(\mu_x\mu_y\mI_d  + \mB^\top\mB)^{-1}\basis{d}{1}} + \delta_y\sqn{\mB(\mu_x\mu_y\mI_d  + \mB^\top\mB)^{-1}\basis{d}{1}}
    \\&\ageq{uses \Cref{lem:B2} and the fact that $\sqn{\basis{d}{1}} = 1$}
    \frac{\delta_x\mu_y^2 + \delta_y\bmu_{xy}^2}{(\mu_x\mu_y + L_{xy}^2)^2}
    \agreater{uses \Cref{ass:xyxy} and the assumption $\mu_y > 0$}
    0,
\end{align*}
where \annotate. Therefore, we obtain the desired inequality $\min_{d\in \{1,2,\ldots\}} B_{\pi,d} > 0$.

Furthermore, the first equation in \cref{eq:x4} can be written as follows:
\begin{equation}
    \begin{bNiceArray}{ccccc}
        \big(\gamma^2 - \frac{\beta}{\alpha}\big) & -1 \\
        -1 & \gamma^2 & -1 \\
        & \Ddots & \Ddots& \Ddots\\
        && -1 & \gamma^2 & -1 \\
        &&& -1 & \gamma^2\\
    \end{bNiceArray}x^* = \frac{4\mu_yA}{\alpha\beta} \basis{d}{1},
\end{equation}
where $\gamma^2$ is defined as follows:
\begin{equation}
    \gamma^2 = \frac{\alpha^2 + \beta^2 + 4\mu_x\mu_y}{\alpha\beta}
    =
    \frac{2(q+1)}{(q-1)},
\end{equation}
where $q > 1$ is defined as follows:
\begin{equation}\label{eq:q2}
    q=\frac{L_{xy}^2 + \mu_x\mu_y}{\bmu_{xy}^2 + \mu_x\mu_y},
\end{equation}
where the denominator is always positive due to \Cref{ass:xyxy}.
Let $(x^\circ,y^\circ) \in \sX \times \sY$ be defined as follows:
\begin{equation}
    x^\circ = \frac{4\mu_y A}{\beta(\alpha - \beta\rho)} \cdot (\rho,\rho^2,\ldots,\rho^d),
    \quad
    y^\circ =  \frac{2 A}{\beta} (\rho,\rho^2,\ldots,\rho^d),
\end{equation}
where $\rho \in (0,1)$ is defined as follows:
\begin{equation}\label{eq:rho2}
    \rho = \frac{\sqrt{q} - 1}{\sqrt{q} + 1}.
\end{equation}
One can verify that $(x^\circ,y^\circ)$  satisfies the following linear system:
\begin{equation}
    \begin{aligned}
        (\mu_x\mu_y\mI_d  + \mB^\top\mB) x^\circ & = \mu_y A\basis{d}{1} +  \frac{\alpha\rho^{d+1}}{\alpha - \beta\rho}\cdot\mu_y A\basis{d}{d}, \\
        \mu_y y^\circ - \mB x^\circ              & =
        -\frac{2 \rho^{d+1}}{\alpha - \beta\rho}\cdot \mu_y A\basis{d}{d},
    \end{aligned}
\end{equation}
which, together with \cref{eq:x4}, implies
\begin{equation}
    \begin{aligned}
        (\mu_x\mu_y\mI_d  + \mB^\top\mB) (x^\circ - x^*) & = \frac{\alpha\rho^{d+1}}{\alpha - \beta\rho}\cdot\mu_y A\basis{d}{d}, \\
        \mu_y(y^\circ - y^*)    - \mB (x^\circ - x^*)    & =
        -\frac{2 \rho^{d+1}}{\alpha - \beta\rho}\cdot \mu_y A\basis{d}{d},
    \end{aligned}
\end{equation}
Hence, we can upper-bound $\norm{x^\circ - x^*}$ as follows:
\begin{align*}
    \norm{x^\circ - x^*}
     & \aeq{uses the linear system~\eqref{eq:x4} and the fact that the matrix $(\mu_x\mu_y\mI_d + \mB^\top \mB)$ is invertible, which is implied by the matrix inequality $(\mu_x\mu_y\mI_d + \mB^\top \mB)\succeq (\mu_x\mu_y + \bmu_{xy}^2) \mI_d$ and \Cref{ass:xyxy}}
    \frac{\alpha\mu_y A\rho^{d+1}}{\alpha - \beta\rho}
    \cdot\norm{(\mu_x\mu_y\mI_d  + \mB^\top\mB)^{-1}\basis{d}{d}}
    \\&\aleq{uses the fact that $\mu_x\mu_y\mI_d + \mB^\top \mB \succeq (\mu_x\mu_y + \bmu_{xy}^2) \mI_d$ and $\norm{\basis{d}{d}} = 1$}
    \frac{\alpha\mu_y A\rho^{d+1}}{(\mu_x\mu_y + \bmu_{xy}^2)(\alpha - \beta\rho)},
\end{align*}
where \annotate. Furthermore, we can upper-bound $\norm{y^\circ - y^*}$ as follows:
\begin{align*}
    \norm{y^\circ - y^*}
     & \aeq{uses the linear system above}
    \frac{1}{\mu_y}\norm*{\mB (x^\circ - x^*)-\frac{2 \rho^{d+1}}{\alpha - \beta\rho}\cdot \mu_y A\basis{d}{d}}
    \\&\aleq{uses the triangle inequality}
    \frac{1}{\mu_y}\norm{\mB (x^\circ - x^*)}
    +\frac{2 A \rho^{d+1}}{\alpha - \beta\rho}\cdot \norm{\basis{d}{d}}
    \\&\aleq{uses the fact that $\norm{\mB} \leq L_{xy}$ and $\norm{\basis{d}{d}}=1$}
    \frac{L_{xy}}{\mu_y}\norm{x^\circ - x^*}
    +\frac{2 A \rho^{d+1}}{\alpha - \beta\rho}
    \\&\aleq{uses the previously obtained upper bound on $\norm{x^\circ - x^*}$}
    \frac{\alpha L_{xy} A\rho^{d+1}}{(\mu_x\mu_y + \bmu_{xy}^2)(\alpha - \beta\rho)}
    +\frac{2 A \rho^{d+1}}{\alpha - \beta\rho}
    \\&=
    \left(2 + \frac{\alpha L_{xy}}{\mu_x\mu_y + \bmu_{xy}^2}
    \right)\frac{A \rho^{d+1}}{\alpha - \beta\rho},
\end{align*}
where \annotate.
Combining the upper bounds on $\norm{x^\circ - x^*}$ and $\norm{y^\circ - y^*}$ gives the desired inequality:
\begin{equation}
    \distsol(x^\circ,y^\circ)
    =
    \delta_x \sqn{x^\circ - x^*}
    +\delta_y \sqn{y^\circ - y^*}
    \leq C_\pi A^2\rho^{2d}.
\end{equation}
It remains to lower-bound $\rho$ as follows:
\begin{align*}
    \rho
     & \aeq{uses the definition of $\rho$ in \cref{eq:rho2}}
    1 - \frac{2}{\sqrt{q} + 1}
    \aeq{uses the definition of $q$ in \cref{eq:q2}}
    1 - \frac{2}{\sqrt{\frac{L_{xy}^2 + \mu_x\mu_y}{\bmu_{xy}^2 + \mu_x\mu_y}} + 1}
    \ageq{uses \Cref{ass:Pi}}
    1 - \frac{2}{\sqrt{\frac{18^2\max\{\bmu_{xy}^2,\mu_x\mu_y\} + \mu_x\mu_y}{\bmu_{xy}^2 + \mu_x\mu_y}} + 1}
    \\&\geq
    1 - \frac{2}{\sqrt{\frac{(18^2-1)\max\{\bmu_{xy}^2,\mu_x\mu_y\}}{\bmu_{xy}^2 + \mu_x\mu_y} + 1} + 1}
    \geq
    1 - \frac{2}{\sqrt{\frac{(18^2-1)}{2} + 1} + 1}
    \geq
    \frac{4}{5},
\end{align*}
where \annotate.
In addition, we can lower-bound $\rho$ as follows:
\begin{align*}
    \rho
     & \aeq{uses the definition of $\rho$ in \cref{eq:rho2}}
    1 - \frac{2}{\sqrt{q} + 1}
    \aeq{uses the definition of $q$ in \cref{eq:q2}}
    1 - \frac{2}{\sqrt{\frac{L_{xy}^2 + \mu_x\mu_y}{\bmu_{xy}^2 + \mu_x\mu_y}} + 1}
    \ageq{uses the assumption $\mu_x\mu_y \geq \max\{\mu_{xy}^2, \mu_{yx}^2\}$}
    1 - \frac{2}{\sqrt{\frac{L_{xy}^2}{2\mu_x\mu_y}}}
    \ageq{uses the definition of $\condxya$ in \cref{eq:kappa} and the definitions of $\delta_x$ and $\delta_y$ in \cref{eq:xyxy}}
    1 - \frac{\sqrt{8}}{\sqrt{\condxya}},
\end{align*}
where \annotate, which concludes the proof.\qed

\newpage

\section{Proof of Theorem~\ref{thm:alg}}\label{sec:proof:alg}

Without loss of generality, we can assume $\mP = \mI_{d_z}$. Otherwise, we can simply make a variable change $z \to \mP^{1/2}z$. That is, we can replace functions $p_i(z)$ and operators $Q_i(z)$ with $p_i(\mP^{-1/2}z)$ and $\mP^{-1/2}Q_i(-\mP^{1/2}z)$.

We start with the following \Cref{lem:p_smooth,,lem:p_strongly_convex,,lem:Q_lipschitz} that describe the basic properties of functions $\pit{i}{z}{k}{t_1,\ldots,t_k}$ and $\pitpr{i}{z}{k}{t_1,\ldots,t_k}$, and operators $Q_i(z)$. The proofs of these lemmas is a trivial utilization of the definitions of functions $\pit{i}{z}{k}{t_1,\ldots,t_k}$ and $\pitpr{i}{z}{k}{t_1,\ldots,t_k}$ on \cref{line:p,line:p_prime} of \Cref{alg}, and the definitions of constants $\Lit{i}{k}{t_1,\ldots,t_k}$, $\Mit{i}{k}{t_1,\ldots,t_k}$, and $\Hit{i}{k}{t_1,\ldots,t_k}$ on \cref{line:L,,line:M,,line:H} of \Cref{alg}. We omit these proofs due to their simplicity.
\begin{lemma}\label{lem:p_smooth}
    For all $k \leq i$, function $\pitpr{i}{z}{k}{t_1,\ldots,t_k}$ is $\Lit{i}{k}{t_1,\ldots,t_k}$-smooth.
\end{lemma}
\begin{lemma}\label{lem:p_strongly_convex}
    For all $k \geq i$, function $\pit{i}{z}{k}{t_1,\ldots,t_k}$ is $\Hit{i}{k}{t_1,\ldots,t_k}$-strongly convex.
\end{lemma}
\begin{lemma}\label{lem:Q_lipschitz}
    For all $1 \leq i \leq n$, operator $Q_i(z)$ is $\Mit{i}{k}{t_1,\ldots,t_k}$-Lipschitz.
\end{lemma}

Next, we establish the key \Cref{lem:alg_main} that describes the convergence properties of \Cref{alg}. The proof is available in \Cref{sec:proof:alg_main}.
\newcommand{\Hsum}[1]{H_{\Sigma}^{#1}}
\newcommand{\rit}[4]{r_{#1}^{#3;#4}(#2)}
\newcommand{\Dit}[4]{D_{#1;#4}^{#2;#3}}
\begin{lemma}\label{lem:alg_main}
    For all $1 \leq k \leq n$ and $\hz \in \Zcon$, the following inequality holds:
    \begin{equation}\label{eq:alg_main}
        \begin{aligned}
            0
             & \geq
            \sum_{i=1}^n\left(
            \rit{i}{\zit{k}{t_1,\ldots,t_{k-1},t_{k}+1/2}}{k}{t_1,\ldots,t_k}
            -\rit{i}{\hz}{k}{t_1,\ldots,t_k}
            \right)
            +\frac{\Hsum{k}}{2}\sqn{\zit{k}{t_1,\ldots,t_{k-1},t_{k}+1/2} - \hz}
            \\&
            +\sum_{i=k+1}^{n} \left(
            \Dit{i}{k}{t_1,\ldots,t_k}{t_{k}+1} - \Dit{i}{k}{t_1,\ldots,t_k}{t_{k}}
            \right),
        \end{aligned}
    \end{equation}
    where function $\rit{i}{z}{k}{t_1,\ldots,t_k} \colon \sZ \to \R$ is defined as follows:
    \begin{equation}\label{eq:r}
        \rit{i}{z}{k}{t_1,\ldots,t_k} = \pit{i}{z}{k}{t_1,\ldots,t_k}
        +\begin{cases}
            \<z,Q_i(\hz)> & i > k    \\
            0             & i \leq k
        \end{cases},
    \end{equation}
    constant $\Hsum{k} \geq 0$ is defined as follows:
    \begin{equation}\label{eq:Hsum}
        \Hsum{k} = \sum_{i=1}^k \Hit{i}{i}{t_1,\ldots,t_i},
    \end{equation}
    and $\Dit{i}{k}{t_1,\ldots,t_k}{t} \geq 0$ is defined as follows:
    \begin{equation}\label{eq:D}
        \Dit{i}{k}{t_1,\ldots,t_k}{t}
        =
        \frac{
            \Lit{i}{k}{t_1\ldots,t_k}\prod_{j=k+1}^{i}\alpha_{T_j-1}^2
            +\Mit{i}{k}{t_1\ldots,t_k}\prod_{j=k+1}^{i}\alpha_{T_j-1}
        }{2}
        \sqn{\zit{i}{t_1,\ldots,t_{k-1},t,0\ldots0} - \hz}.
    \end{equation}
\end{lemma}

Now, we are ready to prove \Cref{thm:alg}. Using the initialization steps on \cref{line:z_in,line:p_in} of \Cref{alg}, the definition of the output $\zout\in \sZ$ on \cref{line:z_out} of \Cref{alg}, and the arguments that are identical to the proof of \Cref{lem:alg_main}, we obtain the following simplified version of \cref{eq:alg_main} for the case $k=0$:
\begin{equation}
    \sum_{i=1}^n \left(r_i^0(\zout) - r_i(\hz)\right)
    \leq
    \sum_{i=1}^n\frac{L_i\prod_{j=1}^{i}\alpha_{T_j-1}^2 + M_i\prod_{j=1}^{i}\alpha_{T_j-1}}{2}\cdot\sqn{z^0 - z}.
\end{equation}
Furthermore, we prove the following \Cref{lem:Zcon} in \Cref{sec:proof:Zcon}.
\begin{lemma}\label{lem:Zcon}
    The output $\zout$ of \Cref{alg} satisfies the inclusion $\zout \in \Zcon$.
\end{lemma}
It remains to upper-bound $\alpha_t$. First, we lower-bound $\alpha_t^{-1}$ as follows:
\begin{align*}
    \alpha_{t}^{-1}
    \aeq{use the definition of $\alpha_t$ in \cref{eq:alpha_t}}
    \tfrac{1}{2} + \tfrac{1}{2}\sqrt{1 + 4\alpha_{t-1}^{-2}}
    \geq
    \tfrac{1}{2} + \alpha_{t-1}^{-1}
    \geq \cdots \geq
    \tfrac{1}{2}t+\alpha_0^{-1}
    \aeq{use the definition of $\alpha_t$ in \cref{eq:alpha_t}}
    \tfrac{1}{2}(t+2),
\end{align*}
where \annotate. This lower bound on $\alpha_{t}^{-1}$ implies
the following upper bound on $\alpha_{T-1}$:
\begin{equation}
    \alpha_{T-1} \leq \frac{2}{T+1},
\end{equation}
which concludes the proof.\qed

\newpage

\subsection{Proofs of Auxiliary Lemmas}

\subsubsection{Proof of Lemma~\ref{lem:alg_main}}\label{sec:proof:alg_main}
We prove this lemma by induction for $k = n,\ldots,1$.

    {\bf Base case ($k=n$).}
Using \cref{line:argmin,line:call} of \Cref{alg}, we get the following relation:
\begin{equation}
    \zit{n}{t_1,\ldots,t_{n-1},t_n+1/2} = \argmin_{z \in \Zcon}\sum_{i=1}^n \pit{i}{z}{n}{t_1,\ldots,t_n}.
\end{equation}
Moreover, function $\sum_{i=1}^n \pit{i}{z}{n}{t_1,\ldots,t_n}$ is $\Hsum{n}$-strongly convex due to \Cref{lem:p_strongly_convex} and the definition of $\Hsum{n}$ in \cref{eq:Hsum}. Hence, using the fact that $\hz \in \Zcon$ and the fact that
\begin{equation}
    \pit{i}{\hz}{n}{t_1,\ldots,t_n} = \rit{i}{\hz}{n}{t_1,\ldots,t_n},
\end{equation}
which is implied by which is implied by \cref{eq:r}, we obtain the following inequality:
\begin{equation}
    \sum_{i=1}^n \rit{i}{\hz}{n}{t_1,\ldots,t_n}
    \geq
    \sum_{i=1}^n \rit{i}{\zit{n}{t_1,\ldots,t_{n-1},t_n+1/2}}{n}{t_1,\ldots,t_n}
    +\frac{\Hsum{n}}{2}\sqn{\zit{n}{t_1,\ldots,t_{n-1},t_n+1/2} - \hz},
\end{equation}
which is nothing else but the desired \cref{eq:alg_main} in the case $k = n$.

    {\bf Introduction step ($k \to k-1$).}
We assume that \cref{eq:alg_main} holds for $2\leq k \leq n-1$, which implies the following:
\begin{align*}
    0
     & \geq
    \sum_{i=k+1}^{n} \left(
    \Dit{i}{k}{t_1,\ldots,t_k}{t_{k}+1} - \Dit{i}{k}{t_1,\ldots,t_k}{t_{k}}
    \right)
    +\frac{\Hsum{k}}{2}\sqn{\zit{k}{t_1,\ldots,t_{k-1},t_{k}+1/2} - \hz}
    \\&
    +\sum_{i=1}^n\left(
    \rit{i}{\zit{k}{t_1,\ldots,t_{k-1},t_{k}+1/2}}{k}{t_1,\ldots,t_k}
    -\rit{i}{\hz}{k}{t_1,\ldots,t_k}
    \right)
    \\&\aeq{uses the definition of functions $\rit{i}{z}{k}{t_1,\ldots,t_k}$ in \cref{eq:r} and the definition of functions $\pit{i}{z}{k}{t_1,\ldots,t_k}$ on \cref{line:p} of \Cref{alg}}
    \sum_{i=k+1}^{n} \left(
    \Dit{i}{k}{t_1,\ldots,t_k}{t_{k}+1} - \Dit{i}{k}{t_1,\ldots,t_k}{t_{k}}
    \right)
    +\frac{\Hsum{k}}{2}\sqn{\zit{k}{t_1,\ldots,t_{k-1},t_{k}+1/2} - \hz}
    \\&
    +\sum_{i=1,i\neq k}^n\left(
    \pitpr{i}{\zit{k}{t_1,\ldots,t_{k-1},t_{k}+1/2}}{k}{t_1,\ldots,t_k}
    -\pitpr{i}{\hz}{k}{t_1,\ldots,t_k}
    \right)
    +\sum_{i=k+1}^n \<Q_i(\hz),\zit{k}{t_1,\ldots,t_{k-1},t_{k}+1/2} - \hz>
    \\&
    +\frac{\Hit{k}{k}{t_1,\ldots,t_k}}{2}\sqn{\zit{k}{t_1,\ldots,t_{k-1},t_{k}+1/2} - \zit{k}{t_1,\ldots,t_k}}
    -\frac{\Hit{k}{k}{t_1,\ldots,t_k}}{2}\sqn{\hz - \zit{k}{t_1,\ldots,t_k}}
    \\&
    +\<\zit{k}{t_1,\ldots,t_{k-1},t_{k}+1/2} - \hz, \git{k}{k}{t_1,\ldots,t_k}>
    \\&\aeq{uses the definition of $\git{k}{k}{t_1,\ldots,t_k}$ on \cref{line:grad} of \Cref{alg}}
    \sum_{i=k+1}^{n} \left(
    \Dit{i}{k}{t_1,\ldots,t_k}{t_{k}+1} - \Dit{i}{k}{t_1,\ldots,t_k}{t_{k}}
    \right)
    +\frac{\Hsum{k}}{2}\sqn{\zit{k}{t_1,\ldots,t_{k-1},t_{k}+1/2} - \hz}
    \\&
    +\sum_{i=1,i\neq k}^n\left(
    \pitpr{i}{\zit{k}{t_1,\ldots,t_{k-1},t_{k}+1/2}}{k}{t_1,\ldots,t_k}
    -\pitpr{i}{\hz}{k}{t_1,\ldots,t_k}
    \right)
    +\sum_{i=k+1}^n \<Q_i(\hz),\zit{k}{t_1,\ldots,t_{k-1},t_{k}+1/2} - \hz>
    \\&
    +\<\zit{k}{t_1,\ldots,t_{k-1},t_{k}+1/2} - \hz, \nabla \pitpr{i}{\zit{k}{t_1,\ldots,t_k}}{k}{t_1,\ldots,t_k}>
    +\frac{\Hit{k}{k}{t_1,\ldots,t_k}}{2}\sqn{\zit{k}{t_1,\ldots,t_{k-1},t_{k}+1/2} - \zit{k}{t_1,\ldots,t_k}}
    \\&
    +\<\zit{k}{t_1,\ldots,t_{k-1},t_{k}+1/2} - \hz, Q_k(\zit{k}{t_1,\ldots,t_k})>
    -\frac{\Hit{k}{k}{t_1,\ldots,t_k}}{2}\sqn{\hz - \zit{k}{t_1,\ldots,t_k}}
    \\&\ageq{uses the convexity and $\Lit{k}{k}{t_1,\ldots,t_k}$-smoothness of function $\pitpr{i}{\zit{k}{t_1,\ldots,t_k}}{k}{t_1,\ldots,t_k}$, where the smoothness property is implied by \Cref{lem:p_smooth}}
    \sum_{i=k+1}^{n} \left(
    \Dit{i}{k}{t_1,\ldots,t_k}{t_{k}+1} - \Dit{i}{k}{t_1,\ldots,t_k}{t_{k}}
    \right)
    +\frac{\Hsum{k}}{2}\sqn{\zit{k}{t_1,\ldots,t_{k-1},t_{k}+1/2} - \hz}
    \\&
    +\sum_{i=1}^n\left(
    \pitpr{i}{\zit{k}{t_1,\ldots,t_{k-1},t_{k}+1/2}}{k}{t_1,\ldots,t_k}
    -\pitpr{i}{\hz}{k}{t_1,\ldots,t_k}
    \right)
    +\sum_{i=k+1}^n \<Q_i(\hz),\zit{k}{t_1,\ldots,t_{k-1},t_{k}+1/2} - \hz>
    \\&
    +\frac{\Hit{k}{k}{t_1,\ldots,t_k} - \Lit{k}{k}{t_1,\ldots,t_k}}{2}\sqn{\zit{k}{t_1,\ldots,t_{k-1},t_{k}+1/2} - \zit{k}{t_1,\ldots,t_k}}
    -\frac{\Hit{k}{k}{t_1,\ldots,t_k}}{2}\sqn{\hz - \zit{k}{t_1,\ldots,t_k}}
    \\&
    +\<\zit{k}{t_1,\ldots,t_{k-1},t_{k}+1/2} - \hz, Q_k(\zit{k}{t_1,\ldots,t_k})>
    \\&\aeq{uses the definition of $\gitq{k}{t_1,\ldots,t_k}$ on \cref{line:delta_q} of \Cref{alg}}
    \sum_{i=k+1}^{n} \left(
    \Dit{i}{k}{t_1,\ldots,t_k}{t_{k}+1} - \Dit{i}{k}{t_1,\ldots,t_k}{t_{k}}
    \right)
    +\frac{\Hsum{k}}{2}\sqn{\zit{k}{t_1,\ldots,t_{k-1},t_{k}+1/2} - \hz}
    \\&
    +\sum_{i=1}^n\left(
    \pitpr{i}{\zit{k}{t_1,\ldots,t_{k-1},t_{k}+1/2}}{k}{t_1,\ldots,t_k}
    -\pitpr{i}{\hz}{k}{t_1,\ldots,t_k}
    \right)
    +\sum_{i=k+1}^n \<Q_i(\hz),\zit{k}{t_1,\ldots,t_{k-1},t_{k}+1/2} - \hz>
    \\&
    +\frac{\Hit{k}{k}{t_1,\ldots,t_k} - \Lit{k}{k}{t_1,\ldots,t_k}}{2}\sqn{\zit{k}{t_1,\ldots,t_{k-1},t_{k}+1/2} - \zit{k}{t_1,\ldots,t_k}}
    -\frac{\Hit{k}{k}{t_1,\ldots,t_k}}{2}\sqn{\hz - \zit{k}{t_1,\ldots,t_k}}
    \\&
    +\<\zit{k}{t_1,\ldots,t_{k-1},t_{k}+1/2} - \hz, Q_k(\zit{k}{t_1,\ldots,t_{k-1},t_k+1/2}) + \gitq{k}{t_1,\ldots,t_k}>
    \\&\ageq{uses the monotonicity of operator $Q_k(z)$}
    \sum_{i=k+1}^{n} \left(
    \Dit{i}{k}{t_1,\ldots,t_k}{t_{k}+1} - \Dit{i}{k}{t_1,\ldots,t_k}{t_{k}}
    \right)
    +\frac{\Hsum{k}}{2}\sqn{\zit{k}{t_1,\ldots,t_{k-1},t_{k}+1/2} - \hz}
    \\&
    +\sum_{i=1}^n\left(
    \pitpr{i}{\zit{k}{t_1,\ldots,t_{k-1},t_{k}+1/2}}{k}{t_1,\ldots,t_k}
    -\pitpr{i}{\hz}{k}{t_1,\ldots,t_k}
    \right)
    +\sum_{i=k}^n \<Q_i(\hz),\zit{k}{t_1,\ldots,t_{k-1},t_{k}+1/2} - \hz>
    \\&
    +\frac{\Hit{k}{k}{t_1,\ldots,t_k} - \Lit{k}{k}{t_1,\ldots,t_k}}{2}\sqn{\zit{k}{t_1,\ldots,t_{k-1},t_{k}+1/2} - \zit{k}{t_1,\ldots,t_k}}
    -\frac{\Hit{k}{k}{t_1,\ldots,t_k}}{2}\sqn{\hz - \zit{k}{t_1,\ldots,t_k}}
    \\&
    +\<\zit{k}{t_1,\ldots,t_{k-1},t_{k}+1/2} - \hz, \gitq{k}{t_1,\ldots,t_k}>
    \\&\aeq{uses the parallelogram rule of the form $\<a,b> = \frac{c}{2}\sqn{a + b/c} - \frac{c}{2}\sqn{a} - \frac{1}{2c}\sqn{b}$}
    \sum_{i=k+1}^{n} \left(
    \Dit{i}{k}{t_1,\ldots,t_k}{t_{k}+1} - \Dit{i}{k}{t_1,\ldots,t_k}{t_{k}}
    \right)
    +\frac{\Hsum{k}}{2}\sqn{\zit{k}{t_1,\ldots,t_{k-1},t_{k}+1/2} - \hz}
    \\&
    +\sum_{i=1}^n\left(
    \pitpr{i}{\zit{k}{t_1,\ldots,t_{k-1},t_{k}+1/2}}{k}{t_1,\ldots,t_k}
    -\pitpr{i}{\hz}{k}{t_1,\ldots,t_k}
    \right)
    +\sum_{i=k}^n \<Q_i(\hz),\zit{k}{t_1,\ldots,t_{k-1},t_{k}+1/2} - \hz>
    \\&
    +\frac{\Hit{k}{k}{t_1,\ldots,t_k} - \Lit{k}{k}{t_1,\ldots,t_k}}{2}\sqn{\zit{k}{t_1,\ldots,t_{k-1},t_{k}+1/2} - \zit{k}{t_1,\ldots,t_k}}
    -\frac{\Hit{k}{k}{t_1,\ldots,t_k}}{2}\sqn{\hz - \zit{k}{t_1,\ldots,t_k}}
    \\&
    +\frac{\Hit{k}{k}{t_1,\ldots,t_k}}{2}\sqn{\zit{k}{t_1,\ldots,t_{k-1},t_{k}+1/2} + (\Hit{k}{k}{t_1,\ldots,t_k})^{-1}\gitq{k}{t_1,\ldots,t_k} - \hz}
    \\&
    -\frac{\Hit{k}{k}{t_1,\ldots,t_k}}{2}\sqn{\zit{k}{t_1,\ldots,t_{k-1},t_{k}+1/2} - \hz}
    -\frac{1}{2\Hit{k}{k}{t_1,\ldots,t_k}}\sqn{\gitq{k}{t_1,\ldots,t_k}}
    \\&\aeq{uses \cref{line:z} of \Cref{alg}, the definition of $\Hsum{k}$ in \cref{eq:Hsum}, and the definition of $\gitq{k}{t_1,\ldots,t_k}$ on \cref{line:delta_q} of \Cref{alg}}
    \sum_{i=k+1}^{n} \left(
    \Dit{i}{k}{t_1,\ldots,t_k}{t_{k}+1} - \Dit{i}{k}{t_1,\ldots,t_k}{t_{k}}
    \right)
    +\frac{\Hsum{k-1}}{2}\sqn{\zit{k}{t_1,\ldots,t_{k-1},t_{k}+1/2} - \hz}
    \\&
    +\sum_{i=1}^n\left(
    \pitpr{i}{\zit{k}{t_1,\ldots,t_{k-1},t_{k}+1/2}}{k}{t_1,\ldots,t_k}
    -\pitpr{i}{\hz}{k}{t_1,\ldots,t_k}
    \right)
    +\sum_{i=k}^n \<Q_i(\hz),\zit{k}{t_1,\ldots,t_{k-1},t_{k}+1/2} - \hz>
    \\&
    +\frac{\Hit{k}{k}{t_1,\ldots,t_k} - \Lit{k}{k}{t_1,\ldots,t_k}}{2}\sqn{\zit{k}{t_1,\ldots,t_{k-1},t_{k}+1/2} - \zit{k}{t_1,\ldots,t_k}}
    -\frac{\Hit{k}{k}{t_1,\ldots,t_k}}{2}\sqn{\hz - \zit{k}{t_1,\ldots,t_k}}
    \\&
    +\frac{\Hit{k}{k}{t_1,\ldots,t_k}}{2}\sqn{\zit{k}{t_1,\ldots,t_{k-1},t_{k}+1} - \hz}
    -\frac{1}{2\Hit{k}{k}{t_1,\ldots,t_k}}\sqn{Q_k(\zit{k}{t_1,\ldots,t_k}) - Q_k(\zit{k}{t_1,\ldots,t_{k-1},t_k+1/2})}
    \\&\ageq{uses the $\Mit{k}{k}{t_1,\ldots,t_k}$-Lipszhitzness of operator $Q_k(z)$, which is implied by \Cref{lem:Q_lipschitz}}
    \sum_{i=k+1}^{n} \left(
    \Dit{i}{k}{t_1,\ldots,t_k}{t_{k}+1} - \Dit{i}{k}{t_1,\ldots,t_k}{t_{k}}
    \right)
    +\frac{\Hsum{k-1}}{2}\sqn{\zit{k}{t_1,\ldots,t_{k-1},t_{k}+1/2} - \hz}
    \\&
    +\sum_{i=1}^n\left(
    \pitpr{i}{\zit{k}{t_1,\ldots,t_{k-1},t_{k}+1/2}}{k}{t_1,\ldots,t_k}
    -\pitpr{i}{\hz}{k}{t_1,\ldots,t_k}
    \right)
    +\sum_{i=k}^n \<Q_i(\hz),\zit{k}{t_1,\ldots,t_{k-1},t_{k}+1/2} - \hz>
    \\&
    +\frac{(\Hit{k}{k}{t_1,\ldots,t_k})^2 - \Lit{k}{k}{t_1,\ldots,t_k}\Hit{k}{k}{t_1,\ldots,t_k} - (\Mit{k}{k}{t_1,\ldots,t_k})^2}{2\Hit{k}{k}{t_1,\ldots,t_k}}\sqn{\zit{k}{t_1,\ldots,t_{k-1},t_{k}+1/2} - \zit{k}{t_1,\ldots,t_k}}
    \\&
    -\frac{\Hit{k}{k}{t_1,\ldots,t_k}}{2}\sqn{\hz - \zit{k}{t_1,\ldots,t_k}}
    +\frac{\Hit{k}{k}{t_1,\ldots,t_k}}{2}\sqn{\zit{k}{t_1,\ldots,t_{k-1},t_{k}+1} - \hz}
    \\&\ageq{uses the definition of $\Hit{k}{k}{t_1,\ldots,t_k}$ on \cref{line:H} of \Cref{alg}}
    \sum_{i=k+1}^{n} \left(
    \Dit{i}{k}{t_1,\ldots,t_k}{t_{k}+1} - \Dit{i}{k}{t_1,\ldots,t_k}{t_{k}}
    \right)
    +\frac{\Hsum{k-1}}{2}\sqn{\zit{k}{t_1,\ldots,t_{k-1},t_{k}+1/2} - \hz}
    \\&
    +\sum_{i=1}^n\left(
    \pitpr{i}{\zit{k}{t_1,\ldots,t_{k-1},t_{k}+1/2}}{k}{t_1,\ldots,t_k}
    -\pitpr{i}{\hz}{k}{t_1,\ldots,t_k}
    \right)
    +\sum_{i=k}^n \<Q_i(\hz),\zit{k}{t_1,\ldots,t_{k-1},t_{k}+1/2} - \hz>
    \\&
    -\frac{\Hit{k}{k}{t_1,\ldots,t_k}}{2}\sqn{\hz - \zit{k}{t_1,\ldots,t_k}}
    +\frac{\Hit{k}{k}{t_1,\ldots,t_k}}{2}\sqn{\zit{k}{t_1,\ldots,t_{k-1},t_{k}+1} - \hz},
\end{align*}
where \annotate.
Furthermore, we obtain the following:
\begin{align*}
    0
     & \ageq{uses the previous inequality and the definition of function $\pitpr{i}{z}{k}{t_1,\ldots,t_k}$ on \cref{line:p_prime} of \Cref{alg}}
    \sum_{i=k+1}^{n} \left(
    \Dit{i}{k}{t_1,\ldots,t_k}{t_{k}+1} - \Dit{i}{k}{t_1,\ldots,t_k}{t_{k}}
    \right)
    +\frac{\Hsum{k-1}}{2}\sqn{\zit{k}{t_1,\ldots,t_{k-1},t_{k}+1/2} - \hz}
    \\&
    +\sum_{i=1}^{k-1}\left(
    \pit{i}{\zit{k}{t_1,\ldots,t_{k-1},t_{k}+1/2}}{k-1}{t_1,\ldots,t_{k-1}}
    -\pit{i}{\hz}{k-1}{t_1,\ldots,t_{k-1}}
    \right)
    \\&
    +\frac{1}{\alpha_{t_k}}\sum_{i=k}^n
    \pit{i}{\alpha_{t_k}\zit{k}{t_1,\ldots,t_{k-1},t_{k}+1/2} + (1-\alpha_{t_k})\hzit{k}{t_k}}{k-1}{t_1,\ldots,t_{k-1}}
    \\&
    -\frac{1}{\alpha_{t_k}}\sum_{i=k}^n
    \pit{i}{\alpha_{t_k}\hz + (1-\alpha_{t_k})\hzit{k}{t_k}}{k-1}{t_1,\ldots,t_{k-1}}
    +\sum_{i=k}^n \<Q_i(\hz),\zit{k}{t_1,\ldots,t_{k-1},t_{k}+1/2} - \hz>
    \\&
    -\frac{\Hit{k}{k}{t_1,\ldots,t_k}}{2}\sqn{\hz - \zit{k}{t_1,\ldots,t_k}}
    +\frac{\Hit{k}{k}{t_1,\ldots,t_k}}{2}\sqn{\zit{k}{t_1,\ldots,t_{k-1},t_{k}+1} - \hz}
    \\&\ageq{uses the definition of $\hzit{k}{t_k+1}$ on \cref{line:hz} of \Cref{alg} and the convexity of functions $\pit{i}{z}{k-1}{t_1,\ldots,t_{k-1}}$}
    \sum_{i=k+1}^{n} \left(
    \Dit{i}{k}{t_1,\ldots,t_k}{t_{k}+1} - \Dit{i}{k}{t_1,\ldots,t_k}{t_{k}}
    \right)
    +\frac{\Hsum{k-1}}{2}\sqn{\zit{k}{t_1,\ldots,t_{k-1},t_{k}+1/2} - \hz}
    \\&
    +\sum_{i=1}^{k-1}\left(
    \pit{i}{\zit{k}{t_1,\ldots,t_{k-1},t_{k}+1/2}}{k-1}{t_1,\ldots,t_{k-1}}
    -\pit{i}{\hz}{k-1}{t_1,\ldots,t_{k-1}}
    \right)
    \\&
    +\frac{1}{\alpha_{t_k}}\sum_{i=k}^n
    \pit{i}{\hzit{k}{t_k+1}}{k-1}{t_1,\ldots,t_{k-1}}
    -\frac{1-\alpha_{t_k}}{\alpha_{t_k}}\sum_{i=k}^n
    \pit{i}{\hzit{k}{t_k}}{k-1}{t_1,\ldots,t_{k-1}}
    \\&
    -\sum_{i=k}^n
    \pit{i}{\hz}{k-1}{t_1,\ldots,t_{k-1}}
    +\sum_{i=k}^n \<Q_i(\hz),\zit{k}{t_1,\ldots,t_{k-1},t_{k}+1/2} - \hz>
    \\&
    -\frac{\Hit{k}{k}{t_1,\ldots,t_k}}{2}\sqn{\hz - \zit{k}{t_1,\ldots,t_k}}
    +\frac{\Hit{k}{k}{t_1,\ldots,t_k}}{2}\sqn{\zit{k}{t_1,\ldots,t_{k-1},t_{k}+1} - \hz}
    \\&\aeq{uses the definition of $\hzit{k}{t_k+1}$ on \cref{line:hz} of \Cref{alg}}
    \sum_{i=k+1}^{n} \left(
    \Dit{i}{k}{t_1,\ldots,t_k}{t_{k}+1} - \Dit{i}{k}{t_1,\ldots,t_k}{t_{k}}
    \right)
    +\frac{\Hsum{k-1}}{2}\sqn{\alpha_{t_k}^{-1}\hzit{k}{t_k + 1} - (1-\alpha_{t_k})\alpha_{t_k}^{-1}\hzit{k}{t_k} - \hz}
    \\&
    +\sum_{i=1}^{k-1}
    \pit{i}{\alpha_{t_k}^{-1}\hzit{k}{t_k + 1} - (1-\alpha_{t_k})\alpha_{t_k}^{-1}\hzit{k}{t_k}}{k-1}{t_1,\ldots,t_{k-1}}
    \\&
    +\sum_{i=k}^n \<Q_i(\hz),\alpha_{t_k}^{-1}\hzit{k}{t_k + 1} - (1-\alpha_{t_k})\alpha_{t_k}^{-1}\hzit{k}{t_k} - \hz>
    \\&
    +\frac{1}{\alpha_{t_k}}\sum_{i=k}^n
    \pit{i}{\hzit{k}{t_k+1}}{k-1}{t_1,\ldots,t_{k-1}}
    -\frac{1-\alpha_{t_k}}{\alpha_{t_k}}\sum_{i=k}^n
    \pit{i}{\hzit{k}{t_k}}{k-1}{t_1,\ldots,t_{k-1}}
    -\sum_{i=1}^n
    \pit{i}{\hz}{k-1}{t_1,\ldots,t_{k-1}}
    \\&
    -\frac{\Hit{k}{k}{t_1,\ldots,t_k}}{2}\sqn{\hz - \zit{k}{t_1,\ldots,t_k}}
    +\frac{\Hit{k}{k}{t_1,\ldots,t_k}}{2}\sqn{\zit{k}{t_1,\ldots,t_{k-1},t_{k}+1} - \hz}
    \\&\ageq{uses the convexity of function $\pit{i}{z}{k-1}{t_1,\ldots,t_{k-1}}$ and $\sqn{\cdot}$}
    \sum_{i=k+1}^{n} \left(
    \Dit{i}{k}{t_1,\ldots,t_k}{t_{k}+1} - \Dit{i}{k}{t_1,\ldots,t_k}{t_{k}}
    \right)
    \\&
    +\frac{1}{\alpha_{t_k}}\left(
    \frac{\Hsum{k-1}}{2}\sqn{\hzit{k}{t_k + 1} - \hz}
    +\sum_{i=k}^n \<Q_i(\hz),\hzit{k}{t_k + 1} - \hz>
    \right)
    \\&
    -\frac{1-\alpha_{t_k}}{\alpha_{t_k}}\left(
    \frac{\Hsum{k-1}}{2}\sqn{\hzit{k}{t_k} - \hz}
    +\sum_{i=k}^n \<Q_i(\hz),\hzit{k}{t_k} - \hz>
    \right)
    \\&
    +\frac{1}{\alpha_{t_k}}\sum_{i=1}^n
    \pit{i}{\hzit{k}{t_k+1}}{k-1}{t_1,\ldots,t_{k-1}}
    -\frac{1-\alpha_{t_k}}{\alpha_{t_k}}\sum_{i=1}^n
    \pit{i}{\hzit{k}{t_k}}{k-1}{t_1,\ldots,t_{k-1}}
    -\sum_{i=1}^n
    \pit{i}{\hz}{k-1}{t_1,\ldots,t_{k-1}}
    \\&
    -\frac{\Hit{k}{k}{t_1,\ldots,t_k}}{2}\sqn{\hz - \zit{k}{t_1,\ldots,t_k}}
    +\frac{\Hit{k}{k}{t_1,\ldots,t_k}}{2}\sqn{\zit{k}{t_1,\ldots,t_{k-1},t_{k}+1} - \hz}
    \\&\aeq{uses the definition of functions $\rit{i}{z}{k-1}{t_1,\ldots,t_{k-1}}$ in \cref{eq:r}}
    \frac{1}{\alpha_{t_k}}\left(
    \frac{\Hsum{k-1}}{2}\sqn{\hzit{k}{t_k + 1} - \hz}
    +\sum_{i=1}^n \left(
    \rit{i}{\hzit{k}{t_k+1}}{k-1}{t_1,\ldots,t_{k-1}}
    -\rit{i}{\hz}{k-1}{t_1,\ldots,t_{k-1}}
    \right)
    \right)
    \\&
    -\frac{1-\alpha_{t_k}}{\alpha_{t_k}}\left(
    \frac{\Hsum{k-1}}{2}\sqn{\hzit{k}{t_k} - \hz}
    +\sum_{i=1}^n \left(
    \rit{i}{\hzit{k}{t_k}}{k-1}{t_1,\ldots,t_{k-1}}
    -\rit{i}{\hz}{k-1}{t_1,\ldots,t_{k-1}}
    \right)
    \right)
    \\&
    -\frac{\Hit{k}{k}{t_1,\ldots,t_k}}{2}\sqn{\hz - \zit{k}{t_1,\ldots,t_k}}
    +\frac{\Hit{k}{k}{t_1,\ldots,t_k}}{2}\sqn{\zit{k}{t_1,\ldots,t_{k-1},t_{k}+1} - \hz}
    \\&
    +\sum_{i=k+1}^{n} \left(
    \Dit{i}{k}{t_1,\ldots,t_k}{t_{k}+1} - \Dit{i}{k}{t_1,\ldots,t_k}{t_{k}}
    \right),
\end{align*}
where \annotate.
Next, we divide both sides of the inequality by $\alpha_{t_k}$ and obtain the following:
\begin{align*}
    0
     & \geq
    \frac{1}{\alpha_{t_k}^2}\left(
    \frac{\Hsum{k-1}}{2}\sqn{\hzit{k}{t_k + 1} - \hz}
    +\sum_{i=1}^n \left(
    \rit{i}{\hzit{k}{t_k+1}}{k-1}{t_1,\ldots,t_{k-1}}
    -\rit{i}{\hz}{k-1}{t_1,\ldots,t_{k-1}}
    \right)
    \right)
    \\&
    -\frac{1-\alpha_{t_k}}{\alpha_{t_k}^2}\left(
    \frac{\Hsum{k-1}}{2}\sqn{\hzit{k}{t_k} - \hz}
    +\sum_{i=1}^n \left(
    \rit{i}{\hzit{k}{t_k}}{k-1}{t_1,\ldots,t_{k-1}}
    -\rit{i}{\hz}{k-1}{t_1,\ldots,t_{k-1}}
    \right)
    \right)
    \\&
    -\frac{\Hit{k}{k}{t_1,\ldots,t_k}}{2\alpha_{t_k}}\sqn{\hz - \zit{k}{t_1,\ldots,t_k}}
    +\frac{\Hit{k}{k}{t_1,\ldots,t_k}}{2\alpha_{t_k}}\sqn{\zit{k}{t_1,\ldots,t_{k-1},t_{k}+1} - \hz}
    \\&
    +\sum_{i=k+1}^{n} \frac{1}{\alpha_{t_k}}\left(
    \Dit{i}{k}{t_1,\ldots,t_k}{t_{k}+1} - \Dit{i}{k}{t_1,\ldots,t_k}{t_{k}}
    \right)
    \\&\aeq{uses the definition of $\Hit{k}{k}{t_1,\ldots,t_k}$ on \cref{line:H} of \Cref{alg}}
    \frac{1}{\alpha_{t_k}^2}\left(
    \frac{\Hsum{k-1}}{2}\sqn{\hzit{k}{t_k + 1} - \hz}
    +\sum_{i=1}^n \left(
    \rit{i}{\hzit{k}{t_k+1}}{k-1}{t_1,\ldots,t_{k-1}}
    -\rit{i}{\hz}{k-1}{t_1,\ldots,t_{k-1}}
    \right)
    \right)
    \\&
    -\frac{1-\alpha_{t_k}}{\alpha_{t_k}^2}\left(
    \frac{\Hsum{k-1}}{2}\sqn{\hzit{k}{t_k} - \hz}
    +\sum_{i=1}^n \left(
    \rit{i}{\hzit{k}{t_k}}{k-1}{t_1,\ldots,t_{k-1}}
    -\rit{i}{\hz}{k-1}{t_1,\ldots,t_{k-1}}
    \right)
    \right)
    \\&
    -\frac{\Lit{k}{k}{t_1,\ldots,t_k} + \Mit{k}{k}{t_1,\ldots,t_k}}{2\alpha_{t_k}}\sqn{\hz - \zit{k}{t_1,\ldots,t_k}}
    +\frac{\Lit{k}{k}{t_1,\ldots,t_k} + \Mit{k}{k}{t_1,\ldots,t_k}}{2\alpha_{t_k}}\sqn{\zit{k}{t_1,\ldots,t_{k-1},t_{k}+1} - \hz}
    \\&
    +\sum_{i=k+1}^{n} \frac{1}{\alpha_{t_k}}\left(
    \Dit{i}{k}{t_1,\ldots,t_k}{t_{k}+1} - \Dit{i}{k}{t_1,\ldots,t_k}{t_{k}}
    \right)
    \\&\aeq{uses the definitions of $\Lit{k}{k}{t_1,\ldots,t_k}$ and $\Mit{k}{k}{t_1,\ldots,t_k}$ on \cref{line:L,line:M} of \Cref{alg}}
    \frac{1}{\alpha_{t_k}^2}\left(
    \frac{\Hsum{k-1}}{2}\sqn{\hzit{k}{t_k + 1} - \hz}
    +\sum_{i=1}^n \left(
    \rit{i}{\hzit{k}{t_k+1}}{k-1}{t_1,\ldots,t_{k-1}}
    -\rit{i}{\hz}{k-1}{t_1,\ldots,t_{k-1}}
    \right)
    \right)
    \\&
    -\frac{1-\alpha_{t_k}}{\alpha_{t_k}^2}\left(
    \frac{\Hsum{k-1}}{2}\sqn{\hzit{k}{t_k} - \hz}
    +\sum_{i=1}^n \left(
    \rit{i}{\hzit{k}{t_k}}{k-1}{t_1,\ldots,t_{k-1}}
    -\rit{i}{\hz}{k-1}{t_1,\ldots,t_{k-1}}
    \right)
    \right)
    \\&
    +\frac{\Lit{k}{k-1}{t_1,\ldots,t_{k-1}}}{2}
    \left(
    \sqn{\zit{k}{t_1,\ldots,t_{k-1},t_k+1} - \hz}
    -\sqn{\zit{k}{t_1,\ldots,t_k} - \hz}
    \right)
    \\&
    +\frac{\Mit{k}{k-1}{t_1,\ldots,t_{k-1}}}{2\alpha_{T_k-1}}
    \left(
    \sqn{\zit{k}{t_1,\ldots,t_{k-1},t_k+1} - \hz}
    -\sqn{\zit{k}{t_1,\ldots,t_k} - \hz}
    \right)
    \\&
    +\sum_{i=k+1}^{n} \frac{1}{\alpha_{t_k}}\left(
    \Dit{i}{k}{t_1,\ldots,t_k}{t_{k}+1} - \Dit{i}{k}{t_1,\ldots,t_k}{t_{k}}
    \right),
\end{align*}
where \annotate.
Using the definition of $\Dit{i}{k}{t_1,\ldots,t_k}{t}$ in \cref{eq:D}, one can verify that $\frac{1}{\alpha_{t_k}}\Dit{i}{k}{t_1,\ldots,t_k}{t}$ does not depend on $t_k$. Moreover, using the definition of $\alpha_{t_k}$ in \cref{eq:alpha_t}, one can verify that $\alpha_0 = 1$ and $\frac{1-\alpha_{t_k}}{\alpha_{t_k}^2} = \frac{1}{\alpha_{t_k-1}^2}$ for $t_k \geq 1$. Hence, we can do the telescoping and obtain for arbitrary $t_k \in \rng{0}{T_k-1}$ the following:
\begin{align*}
    0
     & \geq
    \frac{1}{\alpha_{T_k-1}^2}\left(
    \frac{\Hsum{k-1}}{2}\sqn{\hzit{k}{T_k} - \hz}
    +\sum_{i=1}^n \left(
    \rit{i}{\hzit{k}{T_k}}{k-1}{t_1,\ldots,t_{k-1}}
    -\rit{i}{\hz}{k-1}{t_1,\ldots,t_{k-1}}
    \right)
    \right)
    \\&
    +\frac{\Lit{k}{k-1}{t_1,\ldots,t_{k-1}}}{2}
    \left(
    \sqn{\zit{k}{t_1,\ldots,t_{k-1},T_k} - \hz}
    -\sqn{\zit{k}{t_1,\ldots,t_{k-1},0} - \hz}
    \right)
    \\&
    +\frac{\Mit{k}{k-1}{t_1,\ldots,t_{k-1}}}{2\alpha_{T_k-1}}
    \left(
    \sqn{\zit{k}{t_1,\ldots,t_{k-1},T_k} - \hz}
    -\sqn{\zit{k}{t_1,\ldots,t_{k-1},0} - \hz}
    \right)
    \\&
    +\sum_{i=k+1}^{n} \frac{1}{\alpha_{t_k}}\left(
    \Dit{i}{k}{t_1,\ldots,t_k}{T_k} - \Dit{i}{k}{t_1,\ldots,t_k}{0}
    \right).
\end{align*}
After multiplying both sides of the inequality by $\alpha_{T_k-1}^2$, we obtain for arbitrary $t_k \in \rng{0}{T_k-1}$ the following:
\begin{align*}
    0
     & \geq
    \sum_{i=1}^n \left(
    \rit{i}{\hzit{k}{T_k}}{k-1}{t_1,\ldots,t_{k-1}}
    -\rit{i}{\hz}{k-1}{t_1,\ldots,t_{k-1}}
    \right)
    +\frac{\Hsum{k-1}}{2}\sqn{\hzit{k}{T_k} - \hz}
    \\&
    +\frac{\Lit{k}{k-1}{t_1,\ldots,t_{k-1}}\cdot\alpha_{T_k-1}^2}{2}
    \left(
    \sqn{\zit{k}{t_1,\ldots,t_{k-1},T_k} - \hz}
    -\sqn{\zit{k}{t_1,\ldots,t_{k-1},0} - \hz}
    \right)
    \\&
    +\frac{\Mit{k}{k-1}{t_1,\ldots,t_{k-1}}\cdot\alpha_{T_k-1}}{2}
    \left(
    \sqn{\zit{k}{t_1,\ldots,t_{k-1},T_k} - \hz}
    -\sqn{\zit{k}{t_1,\ldots,t_{k-1},0} - \hz}
    \right)
    \\&
    +\sum_{i=k+1}^{n} \frac{\alpha_{T_k-1}^2}{\alpha_{t_k}}\left(
    \Dit{i}{k}{t_1,\ldots,t_k}{T_k} - \Dit{i}{k}{t_1,\ldots,t_k}{0}
    \right)
    \\&\aeq{use the definition of $\Dit{i}{k}{t_1,\ldots,t_k}{t}$ in \cref{eq:D}}
    \sum_{i=1}^n \left(
    \rit{i}{\hzit{k}{T_k}}{k-1}{t_1,\ldots,t_{k-1}}
    -\rit{i}{\hz}{k-1}{t_1,\ldots,t_{k-1}}
    \right)
    +\frac{\Hsum{k-1}}{2}\sqn{\hzit{k}{T_k} - \hz}
    \\&
    +\frac{\Lit{k}{k-1}{t_1,\ldots,t_{k-1}}\cdot\alpha_{T_k-1}^2}{2}
    \left(
    \sqn{\zit{k}{t_1,\ldots,t_{k-1},T_k} - \hz}
    -\sqn{\zit{k}{t_1,\ldots,t_{k-1},0} - \hz}
    \right)
    \\&
    +\frac{\Mit{k}{k-1}{t_1,\ldots,t_{k-1}}\cdot\alpha_{T_k-1}}{2}
    \left(
    \sqn{\zit{k}{t_1,\ldots,t_{k-1},T_k} - \hz}
    -\sqn{\zit{k}{t_1,\ldots,t_{k-1},0} - \hz}
    \right)
    \\&
    +\frac{\alpha_{T_k-1}^2}{\alpha_{t_k}}\sum_{i=k+1}^{n}
    \frac{
        \Lit{i}{k}{t_1\ldots,t_k}\prod_{j=k+1}^{i}\alpha_{T_j-1}^2
        +\Mit{i}{k}{t_1\ldots,t_k}\prod_{j=k+1}^{i}\alpha_{T_j-1}
    }{2}
    \sqn{\zit{i}{t_1,\ldots,t_{k-1},T_k,0\ldots0} - \hz}
    \\&
    -\frac{\alpha_{T_k-1}^2}{\alpha_{t_k}}\sum_{i=k+1}^{n}
    \frac{
        \Lit{i}{k}{t_1\ldots,t_k}\prod_{j=k+1}^{i}\alpha_{T_j-1}^2
        +\Mit{i}{k}{t_1\ldots,t_k}\prod_{j=k+1}^{i}\alpha_{T_j-1}
    }{2}
    \sqn{\zit{i}{t_1,\ldots,t_{k-1},0\ldots0} - \hz}
    \\&\aeq{uses the definitions of $\Lit{i}{k}{t_1,\ldots,t_k}$ and $\Mit{i}{k}{t_1,\ldots,t_k}$ on \cref{line:L,line:M} of \Cref{alg}}
    \sum_{i=1}^n \left(
    \rit{i}{\hzit{k}{T_k}}{k-1}{t_1,\ldots,t_{k-1}}
    -\rit{i}{\hz}{k-1}{t_1,\ldots,t_{k-1}}
    \right)
    +\frac{\Hsum{k-1}}{2}\sqn{\hzit{k}{T_k} - \hz}
    \\&
    +\frac{\Lit{k}{k-1}{t_1,\ldots,t_{k-1}}\cdot\alpha_{T_k-1}^2}{2}
    \left(
    \sqn{\zit{k}{t_1,\ldots,t_{k-1},T_k} - \hz}
    -\sqn{\zit{k}{t_1,\ldots,t_{k-1},0} - \hz}
    \right)
    \\&
    +\frac{\Mit{k}{k-1}{t_1,\ldots,t_{k-1}}\cdot\alpha_{T_k-1}}{2}
    \left(
    \sqn{\zit{k}{t_1,\ldots,t_{k-1},T_k} - \hz}
    -\sqn{\zit{k}{t_1,\ldots,t_{k-1},0} - \hz}
    \right)
    \\&
    +\sum_{i=k+1}^{n}
    \frac{
    \Lit{i}{k-1}{t_1\ldots,t_{k-1}}\prod_{j=k}^{i}\alpha_{T_j-1}^2
    +\Mit{i}{k-1}{t_1\ldots,t_{k-1}}\prod_{j=k}^{i}\alpha_{T_j-1}
    }{2}
    \sqn{\zit{i}{t_1,\ldots,t_{k-1},T_k,0\ldots0} - \hz}
    \\&
    -\sum_{i=k+1}^{n}
    \frac{
    \Lit{i}{k-1}{t_1\ldots,t_{k-1}}\prod_{j=k}^{i}\alpha_{T_j-1}^2
    +\Mit{i}{k-1}{t_1\ldots,t_{k-1}}\prod_{j=k}^{i}\alpha_{T_j-1}
    }{2}
    \sqn{\zit{i}{t_1,\ldots,t_{k-1},0\ldots0} - \hz}
    \\&\aeq{uses \cref{line:z_table} of \Cref{alg}}
    \sum_{i=1}^n \left(
    \rit{i}{\hzit{k}{T_k}}{k-1}{t_1,\ldots,t_{k-1}}
    -\rit{i}{\hz}{k-1}{t_1,\ldots,t_{k-1}}
    \right)
    +\frac{\Hsum{k-1}}{2}\sqn{\hzit{k}{T_k} - \hz}
    \\&
    +\sum_{i=k}^{n}
    \frac{
    \Lit{i}{k-1}{t_1\ldots,t_{k-1}}\prod_{j=k}^{i}\alpha_{T_j-1}^2
    +\Mit{i}{k-1}{t_1\ldots,t_{k-1}}\prod_{j=k}^{i}\alpha_{T_j-1}
    }{2}
    \sqn{\zit{i}{t_1,\ldots,t_{k-2},t_{k-1}+1,0\ldots0} - \hz}
    \\&
    -\sum_{i=k}^{n}
    \frac{
    \Lit{i}{k-1}{t_1\ldots,t_{k-1}}\prod_{j=k}^{i}\alpha_{T_j-1}^2
    +\Mit{i}{k-1}{t_1\ldots,t_{k-1}}\prod_{j=k}^{i}\alpha_{T_j-1}
    }{2}
    \sqn{\zit{i}{t_1,\ldots,t_{k-1},0\ldots0} - \hz}
    \\&\aeq{use the definition of $\Dit{i}{k}{t_1,\ldots,t_k}{t}$ in \cref{eq:D}}
    \sum_{i=1}^n \left(
    \rit{i}{\hzit{k}{T_k}}{k-1}{t_1,\ldots,t_{k-1}}
    -\rit{i}{\hz}{k-1}{t_1,\ldots,t_{k-1}}
    \right)
    +\frac{\Hsum{k-1}}{2}\sqn{\hzit{k}{T_k} - \hz}
    \\&
    +\sum_{i=k}^{n}\left(
    \Dit{i}{k-1}{t_1,\ldots,t_{k-1}}{t_{k-1}+1}
    -\Dit{i}{k-1}{t_1,\ldots,t_{k-1}}{t_{k-1}}
    \right)
    \\&\aeq{uses \cref{line:return,line:call} of \Cref{alg}}
    \sum_{i=1}^n \left(
    \rit{i}{\zit{k-1}{t_1,\ldots,t_{k-2},t_{k-1}+1/2}}{k-1}{t_1,\ldots,t_{k-1}}
    -\rit{i}{\hz}{k-1}{t_1,\ldots,t_{k-1}}
    \right)
    \\&
    +\frac{\Hsum{k-1}}{2}\sqn{\zit{k-1}{t_1,\ldots,t_{k-2},t_{k-1}+1/2} - \hz}
    +\sum_{i=k}^{n}\left(
    \Dit{i}{k-1}{t_1,\ldots,t_{k-1}}{t_{k-1}+1}
    -\Dit{i}{k-1}{t_1,\ldots,t_{k-1}}{t_{k-1}}
    \right),
\end{align*}
where \annotate. The latter inequality is nothing else but the desired \cref{eq:alg_main} for $k-1$, which concludes the proof.\qed

\newpage
\subsubsection{Proof of Lemma~\ref{lem:Zcon}}\label{sec:proof:Zcon}
We prove by induction that for all $1\leq k\leq n$, the following inclusion holds:
\begin{equation}\label{eq:zit_con}
    \zit{k}{t_1,\ldots,t_{k-1},t_k+1/2} \in \Zcon.
\end{equation}
Indeed, in the base case $k = n$, \cref{eq:zit_con} holds due to \cref{line:argmin,line:call} of \Cref{alg}. Next, we assume the inclusion~\eqref{eq:zit_con} for a fixed $k$ satisfying $1 \leq k \leq n$. We have $\hzit{k}{1} = \zit{k}{t_1,\ldots,t_{k-1},1/2} \in \Zcon$ due to the definition on \cref{line:hz} of \Cref{alg} and the fact that $\alpha_0 = 1$, which is implied by \cref{eq:alpha_t}. Furthermore, for $t_k \geq 1$, we have $\hzit{k}{t_k + 1} \in \Zcon$ due to the definition on \cref{line:hz} of \Cref{alg} and the inclusions $\hzit{k}{t_k}, \zit{k}{t_1,\ldots,t_{k-1},t_k+1/2} \in \Zcon$, and $\alpha_{t_k} \in (0,1)$. Hence, we obtain $\hzit{k}{T_k} \in \Zcon$. This, together with \cref{line:return,line:call} of \Cref{alg}, implies \cref{eq:zit_con} for $k-1$, if $k \geq 2$, or $\zout \in \Zcon$, if $k=1$, which concludes the proof.\qed

\section{Proof of Corollary~\ref{cor:alg}}\label{sec:proof:alg2}

Without loss of generality, we can assume that
\begin{equation}
    \max\left\{\sqrt{{L_{i+1}}/{\epsilon}}, {M_{i+1}}/{\epsilon},1\right\} \geq \max\left\{\sqrt{{L_i}/{\epsilon}}, {M_i}/{\epsilon},1\right\}
    \quad\text{for all}\quad
    i \in \rng{1}{n-1}.
\end{equation}
Otherwise, we can simply reshuffle the pairs $(p_i(z), Q_i(z))$.
To ensure the desired inequality~\eqref{eq:vi_gap}, it is sufficient to choose $\{T_i\}_{i=1}^n$ as follows:
\begin{equation}
    \begin{aligned}
        T_{1}   & = \ceil*{2\cdot\max\left\{\sqrt{{L_1}/{\epsilon}}, {M_1}/{\epsilon}, 1\right\}}, \\
        T_{i+1} & =
        \ceil*{2\cdot \frac{\max\left\{\sqrt{{L_{i+1}}/{\epsilon}}, {M_{i+1}}/{\epsilon}, 1\right\}}{\max\left\{\sqrt{{L_i}/{\epsilon}}, {M_i}/{\epsilon}, 1\right\}}}\quad\text{for}\quad
        i \in \rng{1}{n-1}.
    \end{aligned}
\end{equation}
Indeed, we can lower-bound $\prod_{j=1}^i T_j$ as follows:
\begin{align*}
    \prod_{j=1}^{i} T_j
     & =
    T_1\prod_{j=2}^{i} T_j
    \geq
    2\max\left\{\sqrt{{L_1}/{\epsilon}}, {M_1}/{\epsilon},1\right\}
    \prod_{j=1}^{i-1}
    \frac{2\max\left\{\sqrt{{L_{i+1}}/{\epsilon}}, {M_{i+1}}/{\epsilon},1\right\}}{\max\left\{\sqrt{{L_i}/{\epsilon}}, {M_i}/{\epsilon},1\right\}}
    \\&=
    2^i\max\left\{\sqrt{{L_i}/{\epsilon}}, {M_i}/{\epsilon},1\right\},
\end{align*}
which, together with \cref{eq:alg}, implies \cref{eq:vi_gap}.
Similarly, we can upper-bound $\prod_{j=1}^i T_j$ as follows:
\begin{align*}
    \prod_{j=1}^{i} T_j
     & =
    T_1\prod_{j=2}^{i} T_j
    \\&\leq
    \left(2\max\left\{\sqrt{{L_1}/{\epsilon}}, {M_1}/{\epsilon},1\right\}+1\right)
    \prod_{j=1}^{i-1}
    \left(\frac{2\max\left\{\sqrt{{L_{i+1}}/{\epsilon}}, {M_{i+1}}/{\epsilon},1\right\}}{\max\left\{\sqrt{{L_i}/{\epsilon}}, {M_i}/{\epsilon},1\right\}} + 1\right)
    \\&\leq
    3\max\left\{\sqrt{{L_1}/{\epsilon}}, {M_1}/{\epsilon},1\right\}
    \prod_{j=1}^{i-1}
    \frac{3\max\left\{\sqrt{{L_{i+1}}/{\epsilon}}, {M_{i+1}}/{\epsilon},1\right\}}{\max\left\{\sqrt{{L_i}/{\epsilon}}, {M_i}/{\epsilon},1\right\}}
    \\&=
    3^i\max\left\{\sqrt{{L_i}/{\epsilon}}, {M_i}/{\epsilon},1\right\}
    \\&\leq
    3^n\max\left\{\sqrt{{L_i}/{\epsilon}}, {M_i}/{\epsilon},1\right\},
\end{align*}
Finally, it is easy to verify that \Cref{alg} performs $\prod_{j=1}^{i}T_j$ computations of the gradient $\nabla p_i(z)$ on \cref{line:grad} and $2^i\cdot\prod_{j=1}^{i}T_j$ computations of operator $Q_i(z)$ on \cref{line:grad,line:delta_q}, which concludes the proof.\qed

\newpage
\section{Proof of Theorem~\ref{thm:opt_alg}}\label{sec:proof:opt_alg}
Let $z = (x,y)\in \sZ = \sX \times \sY$ and $z' = (x',y')\in \sZ = \sX \times \sY$. Then we can upper-bound $\norm{\nabla p_1(z) - \nabla p_1(z')}_{\mP^{-1}}$ as follows:
\begin{align*}
    \norm{\nabla p_1(z) - \nabla p_1(z')}_{\mP^{-1}}
     & \aeq{uses the definition of function $p_1(z)$ in \cref{eq:p} and the definition of matrix $\mP$ in \cref{eq:P}}
    \delta_x^{-1/2}\cdot \norm{\nabla f(x) - \nabla f(x')}
    \\&\aleq{uses the smoothness property in \Cref{ass:x}}
    \delta_x^{-1/2}\cdot L_x\norm{x - x'}
    \aeq{uses the definition of matrix $\mP$ in \cref{eq:P}}
    \frac{L_x}{\delta_x}\norm{x - x'}_\mP
    \aeq{uses the definition of $\condx$ in \cref{eq:kappa}}
    \condx\norm{x - x'}_\mP,
\end{align*}
where \annotate. Hence, we can choose $L_1 = \condx$, and similarly, we can choose $L_2 = \condy$.
In addition, we can upper-bound $\norm{\nabla p_3(z) - \nabla p_3(z')}_{\mP^{-1}}$ as follows:
\begin{align*}
    \norm{\nabla p_3(z) - \nabla p_3(z')}_{\mP^{-1}}
     & \aeq{uses the definition of function $p_3(z)$ in \cref{eq:p} and the definition of matrix $\mP$ in \cref{eq:P}}
    \sqrt{
        \frac{\beta_x^2}{\delta_x}\sqn{\mB^\top\mB(x-x')}
        +\frac{\beta_y^2}{\delta_y}\sqn{\mB\mB^\top(y-y')}
    }
    \\&\aleq{uses \Cref{ass:xy}}
    \sqrt{
    \frac{\beta_x^2L_{xy}^4}{\delta_x}\sqn{x-x'}
    +\frac{\beta_y^2L_{xy}^4}{\delta_y}\sqn{y-y'}
    }
    \\&\aeq{uses the definition of matrix $\mP$ in \cref{eq:P}}
    \max\left\{\frac{\beta_x L_{xy}^2}{\delta_x},\frac{\beta_y L_{xy}^2}{\delta_y}\right\}
    \sqn{z-z'}_{\mP}
    \\&\aeq{uses the definition of $\condxya$ in \cref{eq:kappa}}
    \condxya\cdot\max\left\{\beta_x\delta_y,\beta_y\delta_x\right\}
    \norm{z-z'}_{\mP}
    \\&\aleq{uses the definitions of $\beta_x$ and $\beta_y$ in \cref{eq:beta_xy}}
    \condxya\cdot\max\left\{\frac{\delta_y}{4L_y},\frac{\delta_x}{4L_x}\right\}
    \norm{z-z'}_{\mP}
    \\&\aeq{uses the definitions of $\condx$ and $\condy$ in \cref{eq:kappa}}
    \condxya\cdot\max\left\{\frac{1}{4\condx},\frac{1}{4\condy}\right\}
    \norm{z-z'}_{\mP}
    \\&\aleq{uses the fact that $\condx,\condy \geq 1$, which is implied by \Cref{ass:Pi}}
    \tfrac{1}{4}\condxya
    \norm{z-z'}_{\mP},
\end{align*}
where \annotate. Hence, we can choose $L_3 = \condxya$.
Furthermore, we can upper-bound $\norm{Q_3(z) - Q_3(z')}_{\mP^{-1}}$ as follows:
\begin{align*}
    \norm{Q_3(z) - Q_3(z')}_{\mP^{-1}}
     & \aeq{uses the definition of operator $Q_3(z)$ in \cref{eq:Q}}
    \norm*{
        \begin{bNiceMatrix}[c]
            \mO_{d_x} & \mB^\top  \\
            -\mB      & \mO_{d_y}
        \end{bNiceMatrix}
        \begin{bNiceMatrix}
            x-x' \\ y-y'
        \end{bNiceMatrix}
    }_{\mP^{-1}}
    \\&\aeq{use the definition of matrix $\mP$ in \cref{eq:P}}
    \sqrt{\delta_x^{-1}\sqn{\mB^\top(y-y')} + \delta_y^{-1}\sqn{\mB(x-x')}}
    \\&\aleq{uses \Cref{ass:xy}}
    \sqrt{\frac{L_{xy}^2}{\delta_x}\sqn{y-y'} + \frac{L_{xy}^2}{\delta_y}\sqn{x-x'}}
    \\&\aeq{uses the definition of $\condxya$ in \cref{eq:kappa}}
    \sqrt{\condxya\delta_y\sqn{y-y'} + \condxya\delta_x\sqn{x-x'}}
    \\&\aeq{use the definition of matrix $\mP$ in \cref{eq:P}}
    \sqrt{\condxya}\norm{z-z'}_\mP,
\end{align*}
where \annotate. Hence, we can choose $M_3 = \sqrt{\condxya}$.

Next, we can obtain the following:
\begin{align*}
    \epsilon n \sqn{\zin - z^*}_\mP
     & \ageq{uses \Cref{cor:alg}}
    p(\zout) - p(z^*) + \<Q(z),\zout - z^*>
    \\&\aeq{the definitions of functions $p_i(z)$ and operators $Q_i(z)$ in \cref{eq:p,eq:Q}}
    f(\xout) - f(x^*) + g(\yout) - g(y^*) + \<\mB^\top y^*,\xout - x^*> - \<\mB x^*, \yout - y^*>
    \\&
    +\frac{\beta_x}{2}\sqn{\mB \xout - \nabla g(\yin)} + \frac{\beta_y}{2}\sqn{\mB^\top \yout + \nabla f(\xin)}
    \\&
    -\frac{\beta_x}{2}\sqn{\mB x^* - \nabla g(\yin)} - \frac{\beta_y}{2}\sqn{\mB^\top y^* + \nabla f(\xin)}
    \\&\ageq{uses Young's inequality}
    f(\xout) - f(x^*) + g(\yout) - g(y^*) + \<\mB^\top y^*,\xout - x^*> - \<\mB x^*, \yout - y^*>
    \\&
    +\frac{\beta_x}{4}\sqn{\mB (\xout - x^*)} + \frac{\beta_y}{4}\sqn{\mB^\top (\yout - y^*)}
    \\&
    -\beta_x\sqn{\mB x^* - \nabla g(\yin)} - \beta_y\sqn{\mB^\top y^* + \nabla f(\xin)}
    \\&\aeq{uses the optimality conditions~\eqref{eq:opt}}
    \bg_f(\xout,x^*) + \bg_g(\yout,y^*)
    +\frac{\beta_x}{4}\sqn{\mB (\xout - x^*)} + \frac{\beta_y}{4}\sqn{\mB^\top (\yout - y^*)}
    \\&
    -\beta_x\sqn{\nabla g(y^*) - \nabla g(\yin)} - \beta_y\sqn{\nabla f(x^*) - \nabla f(\xin)}
    \\&\ageq{uses the smoothness properties in \Cref{ass:x,ass:y}}
    \bg_f(\xout,x^*) + \bg_g(\yout,y^*)
    +\frac{\beta_x}{4}\sqn{\mB (\xout - x^*)} + \frac{\beta_y}{4}\sqn{\mB^\top (\yout - y^*)}
    \\&
    -2\beta_x \bg_f(\xin,x^*) - 2\beta_y\bg_g(\yin,y^*)
    \\&\ageq{uses the definitions of $\beta_x$ and $\beta_y$ in \cref{eq:beta_xy}}
    \bg_f(\xout,x^*) + \bg_g(\yout,y^*)
    +\frac{\beta_x}{4}\sqn{\mB (\xout - x^*)} + \frac{\beta_y}{4}\sqn{\mB^\top (\yout - y^*)}
    \\&
    -\tfrac{1}{2} \bg_f(\xin,x^*) - \tfrac{1}{2}\bg_g(\yin,y^*),
\end{align*}
where \annotate.

Using the definitions of functions $p_i(z)$ and operators $Q_i(z)$ in \cref{eq:p,eq:Q}, \Cref{ass:xy,ass:xyxy}, and \Cref{lem:sol}, we can conclude that $\proj_\cS(\zin) = \proj_\cS(\zout)$. Hence, we get the following:
\begin{equation}
    \beta_x\sqn{\mB (\xout - x^*)} + \beta_y\sqn{\mB^\top (\yout - y^*)}
    \geq
    \beta_x\mu_{xy}^2\sqn{\xout - x^*} + \beta_y\mu_{yx}^2\sqn{\yout - y^*},
\end{equation}
which implies the following:
\begin{align*}
    \epsilon n \sqn{\zin - z^*}_\mP
     & \geq
    \bg_f(\xout,x^*) + \bg_g(\yout,y^*)
    -\tfrac{1}{2} \bg_f(\xin,x^*) - \tfrac{1}{2}\bg_g(\yin,y^*)
    \\&
    +\tfrac{1}{4}\beta_x\mu_{xy}^2\sqn{\xout - x^*} + \tfrac{1}{4}\beta_y\mu_{yx}^2\sqn{\yout - y^*}
    \\&\ageq{uses the strong convexity properties in \Cref{ass:x,ass:y}}
    \tfrac{3}{4}\bg_f(\xout,x^*) + \tfrac{3}{4}\bg_g(\yout,y^*)
    -\tfrac{1}{2} \bg_f(\xin,x^*) - \tfrac{1}{2}\bg_g(\yin,y^*)
    \\&
    +\left(\tfrac{1}{8}\mu_x + \tfrac{1}{4}\beta_x\mu_{xy}^2\right)\sqn{\xout - x^*}
    +\left(\tfrac{1}{8}\mu_y + \tfrac{1}{4}\beta_y\mu_{yx}^2\right)\sqn{\yout - y^*}
    \\&\ageq{uses the definitions of $\beta_x$ and $\beta_y$ in \cref{eq:beta_xy} and the definitions of $\delta_x$ and $\delta_y$ in \cref{eq:xyxy}}
    \tfrac{3}{4}\bg_f(\xout,x^*) + \tfrac{3}{4}\bg_g(\yout,y^*)
    -\tfrac{1}{2} \bg_f(\xin,x^*) - \tfrac{1}{2}\bg_g(\yin,y^*)
    \\&
    +\tfrac{1}{16}\delta_x\sqn{\xout - x^*}
    +\tfrac{1}{16}\delta_y\sqn{\yout - y^*}
    \\&\aeq{uses the definition of $\mP$ in \cref{eq:P}}
    \tfrac{3}{4}\bg_f(\xout,x^*) + \tfrac{3}{4}\bg_g(\yout,y^*)
    -\tfrac{1}{2} \bg_f(\xin,x^*) - \tfrac{1}{2}\bg_g(\yin,y^*)
    \\&
    +\tfrac{1}{16}\sqn{\zout - z^*}_{\mP},
\end{align*}
where \annotate. Furthermore, we have $\sqn{\zin - z^*}_\mP = \distsol(\xin,\yin)$ and $\sqn{\zout - z^*}_\mP = \distsol(\xout,\yout)$ due to the definition of $z^*$. Hence, we obtain the following inequality:
\begin{equation}
    \begin{multlined}
        \distsol(\xout,\yout) + 12\bg_f(\xout,x^*) + 12\bg_g(\yout,y^*)\\
        \leq
        16\epsilon n\distsol(\xin,\yin) + 8\bg_f(\xin,x^*) + 8\bg_g(\yin,y^*).
    \end{multlined}
\end{equation}
Choosing $\epsilon = \frac{1}{24n} = \frac{1}{72}$ concludes the proof.\qed

\newpage
\section{Proof of Corollary~\ref{cor:opt_alg}}\label{sec:proof:opt_alg2}
We use a restarted version of \Cref{alg}. That is, we apply \Cref{alg} $T$ times and use the output at each run as the input for the next run. Formally, by $\zin^t = (\xin^t,\yin^t)$ and $\zout^t = (\xout^t,\yout^t)$ we denote the input and the output of \Cref{alg} at $t$-th run, where $t \in \rng{0}{T-1}$. Then we have $\zin^0 = 0$ and $\zin^{t+1} = \zout^{t}$ for all $t \in \rng{0}{T-1}$. Hence, we can upper-bound $\distsol(\xin^T,\yin^T)$ as follows:
\begin{align*}
    \distsol(\xin^T,\yin^T)
     & \aleq{use the definition of $\Psi(z)$ in \cref{eq:Psi}}
    \Psi(\zin^T)
    \\&\aleq{uses \Cref{thm:opt_alg}}
    \left(\tfrac{2}{3}\right)^T\Psi(\zin^0)
    \\&\aeq{use the definition of $\Psi(z)$ in \cref{eq:Psi}, where $(x^*,y^*) = \proj_{\cS}(\zin^0)$}
    \left(\tfrac{2}{3}\right)^T\left(\distsol(\xin^0,\yin^0) + 12 \bg_f(\xin^0,x^*) + 12 \bg_g(\yin^0,y^*)\right)
    \\&\aleq{uses the smoothness properties in \Cref{ass:x,ass:y}}
    \left(\tfrac{2}{3}\right)^T\left(\distsol(\xin^0,\yin^0) + 12 L_x\sqn{\xin^0-x^*} + 12 L_y\sqn{\yin^0-y^*}\right)
    \\&\aeq{uses the definitions of $\condx$ and $\condy$ in \cref{eq:kappa} and the definition of $\distsol$ in \cref{eq:distsol}}
    \left(\tfrac{2}{3}\right)^T\left(1 + 12\condx + 12\condy\right)\distsol(\xin^0,\yin^0)
    \\&\aeq{uses the definitions $\zin^0 = 0$, $R^2 = \distsol(0,0)$ and $c = 1 + 12\condx + 12\condy$}
    \left(\tfrac{2}{3}\right)^TcR^2,
\end{align*}
where \annotate.
Next, we choose $T$ as follows:
\begin{equation}
    T = \ceil*{\frac{\log (cR^2/\epsilon)}{\log (3/2)}},
\end{equation}
which implies $\distsol(\xin^T,\yin^T) \leq \epsilon$. Note that $T \geq 0$ due to the fact that $\epsilon \leq R^2$ and $c \geq 1$. In addition, we can upper-bound $T$ as follows:
\begin{align*}
    T
     & \aleq{uses the properties of $\ceil{\cdot}$}
    \frac{\log (cR^2/\epsilon)}{\log (3/2)} + 1
    \\&=
    \left(\frac{1}{\log (3/2)} + \frac{1}{\log (cR^2/\epsilon)}\right)\log (cR^2/\epsilon)
    \\&\aleq{uses the assumption $\epsilon \leq R^2$}
    \left(\frac{1}{\log (3/2)} + \frac{1}{\log (c)}\right)\log (cR^2/\epsilon)
    \\&\aeq{uses the deifnition $c = 1 + 12\condx + 12\condy$}
    \left(\frac{1}{\log (3/2)} + \frac{1}{\log (1 + 12\condx + 12\condy)}\right)\log (cR^2/\epsilon)
    \\&\aleq{uses the fact that $\kappa_x,\kappa_y \geq 1$, which is implied by \Cref{ass:Pi}}
    \left(\frac{1}{\log (3/2)} + \frac{1}{\log (25)}\right)\log (cR^2/\epsilon)
    \\&=
    \cO\left(\log \frac{cR^2}{\epsilon}\right),
\end{align*}
where \annotate.

It remains to combine \cref{eq:LM} in \Cref{thm:opt_alg} and multiply $T$ by the appropriate number of computations of the gradients $\nabla p_i(z)$ and operators $Q_i(z)$, which are provided by \Cref{cor:alg}. Note that the computation of the gradients $\nabla p_1(z)$ and $\nabla p_2(z)$ is equivalent to the computation of the gradients $\nabla f(x)$ and $\nabla g(y)$, respectively. The computation of the gradient $\nabla p_3(z)$ and operator $Q_3(z)$ requires $\cO(1)$ matrix-vector multiplications with the matrices $\mB$ and $\mB^\top$, as well as a single computation of the gradients $\nabla f(\xin)$ and $\nabla g(\yin)$ at the beginning of the algorithm, which concludes the proof. \qed

\end{document}